\documentclass[12 pt]{article}
\counterwithin*{equation}{section}
\usepackage[all]{xy}
\usepackage{pdflscape}
\usepackage[hypcap=false]{caption}
\usepackage{amsmath,amssymb,amsfonts,latexsym,float,graphics,epsfig} 
\usepackage{setspace}
\usepackage{xcolor}
\usepackage[top=2cm, bottom=2cm, left=2cm, right=2cm]{geometry}
\usepackage[colorlinks]{hyperref}
\usepackage{tikz}
\usetikzlibrary{cd}

\title{On Non-autonomous Hamiltonian Dynamics, Dual Spaces, and Kinetic Lifts}  

\begin{document}

\maketitle

\begin{center}
Begüm Ateşli\footnote{E-mail: 
\href{mailto:b.atesli@gtu.edu.tr}{b.atesli@gtu.edu.tr}, corresponding author}\\
Department of Mathematics, \\ Gebze Technical University, 41400 Gebze,
Kocaeli, Turkey.

\bigskip

O\u{g}ul Esen\footnote{E-mail: 
\href{mailto:oesen@gtu.edu.tr}{oesen@gtu.edu.tr}}\\
Department of Mathematics, \\ Gebze Technical University, 41400 Gebze,
Kocaeli, Turkey.

\bigskip

Manuel de Le\'on\footnote{E-mail: \href{mailto:mdeleon@icmat.es}{mdeleon@icmat.es}}
\\ Instituto de Ciencias Matem\'aticas \\
 Consejo Superior de Investigaciones Cient\'ificas
 \\
C/ Nicol\'as Cabrera, 13--15, 28049, Madrid, Spain
\\
and
\\
Real Academia de Ciencias de Espa\~na
\\
C/ Valverde, 22, 28004 Madrid, Spain.

\bigskip 

Cristina Sard\'on\footnote{E-mail: \href{mailto:mariacristina.sardon@upm.es}{mariacristina.sardon@upm.es}}
\\ Department of Applied Mathematics 
\\ Universidad Polit\'ecnica de Madrid 
\\ C/ Jos\'e Guti\'errez Abascal, 2, 28006, Madrid. Spain.\footnote{2nd e-mail: \href{mailto:mcsardon@icade.comillas.edu}{mcsardon@icade.comillas.edu }}
\\ Faculty of Economics and Bussiness Administration
\\ ICADE-Comillas, Madrid.
\\ C/ Alberto Aguilera, 23, 28015, Madrid. Spain.

\begin{abstract}
Vlasov kinetic theory is the dynamics of a bunch of particles flowing according to symplectic Hamiltonian dynamics. More recently, this geometry has been extended to contact Hamiltonian dynamics. In this paper, we introduce geometric kinetic theories within the framework of cosymplectic and cocontact manifolds to extend the present literature to time-dependent dynamics. The cosymplectic and the cocontact kinetic theories are obtained in terms of both momentum variables and density functions. These alternative realizations are linked via Poisson/momentum maps. Furthermore, in cocontact geometry, we introduce a hierarchical analysis of nine distinct dynamical motions as various manifestations of Hamiltonian, evolution, and gradient flows. 

\smallskip

\smallskip


\end{abstract}

\end{center}
\tableofcontents
\onehalfspacing

\numberwithin{equation}{section}

\setlength{\parindent}{2em}
\setlength{\parskip}{3ex}

\section{Introduction}

It is well-established that symplectic manifolds provide a proper geometric framework for the classical mechanics of autonomous Hamiltonian systems \cite{AbMa78, Arn,guil90,leon89, holm2009geometric, Marsden1999}. In symplectic Hamiltonian dynamics, the symplectic volume is preserved as a manifestation of the Liouville theorem, and the energy (the Hamiltonian function itself) is conserved along the motion. These characteristics ensure the non-dissipative nature of symplectic Hamiltonian dynamics. Symplectic manifolds are particular instances of Poisson manifolds. The operation $H \mapsto X_H^s$, which maps a Hamiltonian function to the (symplectic) Hamiltonian vector field (generating the dynamics), has a kernel consisting only of constant functions.

\noindent
\textbf{Cosymplectic Dynamics.} If the differential equations governing a physical system have explicit time dependence, the geometric analysis of this system goes beyond the realm of symplectic geometry. One way to deal with the geometries of nonautonomous dynamics is to employ cosymplectic geometry \cite{BazzGoer15, Cape, LeonSar2, LeonTuyn}. Cosymplectic manifolds are odd-dimensional counterparts of symplectic manifolds, where an extra dimension (apart from the position and the momenta) can be interpreted as the time variable. In cosymplectic Hamiltonian dynamics, both the cosymplectic volume and the Hamiltonian are conserved by the Hamiltonian vector field, but not by the so-called evolution vector field, the latter providing the true dynamics. Thus, cosymplectic Hamiltonian dynamics is reversible as well. Cosymplectic manifolds are also particular instances of Poisson manifolds. However, in this case, the operation $H\mapsto X_H^{cs}$, which maps a Hamiltonian function to the (cosymplectic) Hamiltonian dynamics, has a kernel larger than constant functions. This non-trivial kernel brings an advantage by allowing three different types of vector fields (namely, Hamiltonian, evolution, and gradient) used to model various physical structures in the literature. We list all these vector fields in Table \ref{cosymplectic-table}.

\noindent
\textbf{Contact Dynamics.} 
Another odd-dimensional counterpart of symplectic manifolds is contact manifolds. In this case, the extra dimension is used to record the action term induced by the Hamiltonian function \cite{Arn,Liber87}. Geometrically, contact manifolds differ from symplectic and cosymplectic manifolds, since they are Jacobi manifolds. In terms of Hamiltonian dynamics, this leads to the violation of the conservation laws (of volume and energy) that exist in the Poisson setting \cite{BrCrTa17,de2019contact,esen2021contact}. This, physically, provide a framework that is proper for irreversible dynamical motions such as thermodynamically processes \cite{Bravettithermo,Goto15,Grmela14,Mrugala,Houches}. On contact manifolds, one can define three different types of dynamical systems, namely contact Hamiltonian vector field, contact evolution vector field (in this work, we refer to it as energy vector field to avoid confusion with the evolution vector field in cosymplectic geometry), and quantomorphic (strict contact) vector field. Table \ref{contact-table} is listing all these fields. An important fact in contact Hamiltonian systems is that here the energy dissipates at a known rate, but this is not a consequence of the variation with respect to time.

\noindent
\textbf{Cocontact Dynamics.} 
Cocontact dynamics has emerged as a compelling formalism for studying nonautonomous dissipative systems, offering an extension beyond conventional contact dynamics \cite{LeGaGrMuRi23}. In this framework, two additional dimensions, time and action, intertwine with the traditional position and momenta variables, providing a richer geometric structure. Cocontact dynamics possesses characteristics from both cosymplectic and contact realizations. Accordingly, it is creating a unique theory of dissipative and time-dependent dynamics. Cocontact dynamics has been gathering interest in various research domains. For instance, \cite{Leon-cocontact23} has explored its application to Hamilton-Jacobi theory, while  \cite{LeGaMuRiRo23, Ri23} has investigated its relevance in field theory. \cite{azuaje2023canonical} has exhibited its connection with diffeomorphisms groups, and \cite{MaNo17} has explored cocontact geometry within higher-order theories. Additionally, a unified formalism has been proposed in \cite{RiTo23}.  

\noindent
\textbf{Goal 1.} In the present work, one of the novel results is to determine nine different dynamical motions on cocontact manifolds and to present them in a hierarchic order in Section \ref{sec-coco-ham}, Section \ref{sec-coco-evo}, and Section \ref{sec-coco-grad}. For the summary of all these fields see Table  \ref{cocontact-table} provided in the summary. Before stating the second goal let us first examine the symplectic geometry underlying Vlasov kinetic theory.

\noindent
\textbf{Vlasov Dynamics.} 
Another goal of this work is to extend the geometric Vlasov dynamics of plasma motion to the cosymplectic and cocontact frameworks. First, let us examine the geometric Vlasov dynamics in the upcoming paragraph and then provide more concrete comments on the novel results (in the direction of kinetic theories) presented in the paper. Consider a non-relativistic and collisionless plasma  resting in $Q\subset \mathbb{R}^3$. The dynamics is determined by plasma density function $f$ defined on the mometum phase space $T^*Q$ equipped with the Darboux' coordinates $(q^i,p_i)$. The equation of motions is a coupled integro-differential system 
\begin{equation}\label{PV}
\begin{split}
\frac{\partial f}{\partial t}+\frac{1}{m}p_{i}\frac{\partial f}{\partial
q^{i}}-e\frac{\partial ^{2}\phi }{\partial q^{i}}\frac{\partial f}{\partial
p_{i}}=0.\\ 
\nabla _{q}^{2}\phi _{f}(q)=-e\int f(q,p)d^{3}p
\end{split}
\end{equation}
called Vlasov-Poisson equations. Here, $e$ is the charge and $\phi$ is the potential. The geometric analysis of plasma dynamics is available in various literature for example see the following list \cite{marsden1983hamiltonian,marsden83b,marsden1982hamiltonian,morrison80,morrison1981hamiltonian,morrison1982poisson} for the early treatments. In those studies, it is well established that the Vlasov–Poisson system \eqref{PV} admits a Hamiltonian formulation. More precisely, the Vlasov equation is fitting the Lie-Poisson (a Poisson framework available on the dual of a Lie algebra) picture whereas the Poisson equation is realized as the momentum map due to the gauge invariance of Hamiltonian dynamics.

\noindent
\textbf{Momentum-Vlasov Dynamics.} A  discussion was conducted on the Lie-Poisson formulation of the Vlasov equation \cite{EsGu12,Gu10}. This exploration led to the derivation of intermediate equations concerning the space of one-forms on $T^*Q$. Specifically, a pure Lie-Poisson formulation was developed, written on the dual of the Lie algebra of Hamiltonian vector fields. A crucial step in this process involved the concrete determination of the dual space. Once the dual space was accurately identified, the dynamics of a dual element, namely a one-form denoted as $\Pi$, became governed by the Hamiltonian vector field $X^s_H$, as expressed by what we call momentum-Vlasov equation
\begin{equation}\label{MV-intro}
\dot{\Pi}=-\mathcal{L}_{X^s_H}\Pi.
\end{equation} 
The connection between the Vlasov equation (referred to as the former equation in equation \eqref{PV}) and the momentum-Vlasov equation \eqref{MV-intro} was established through the dual mapping of the Lie algebra homomorphism $H\mapsto X^s_H$. This mapping was computed to yield the relationship:
\begin{equation}\label{density}
f=\mathrm{div} \sharp^s(\Pi),
 \end{equation}
Here, $\sharp^s$ denotes the  musical isomorphism induced from the symplectic two-form. It is from the space of one-forms to the space of vector fields. In this context, $\mathrm{div}$ signifies divergence with respect to the symplectic volume. The relation expressed in the density equation emerged as a result of the duality of the Lie algebra homomorphism, thus functioning as both a Poisson and a momentum map. Within this geometric framework, the Poisson equation (the latter equation in \eqref{PV}) is the momentum map due to the gauge invariance of Hamiltonian dynamics.

One interesting aspect of the momentum-Vlasov equation \eqref{MV-intro} lies in its purely geometric derivation. To be more precise, starting from the Hamiltonian vector field $X^s_H$ and lifting it to the cotangent bundle, followed by taking the vertical representative, a generalized vector field $V\widehat{X}^s_H$ is obtained. This vector field governs the motion of the sections $\Pi$ of the cotangent fibrations over $T^*Q$. This geometric pathway precisely corresponds to the momentum-Vlasov equation \eqref{MV-intro}. Such a geometrization served as motivation for further analysis conducted on fluid dynamics \cite{EsGrMiGu19,EsGu11}. The algebraic structure underlying the dynamics of the momentum-Vlasov equation received more detailed examination in \cite{EsSu21}. 

\noindent
\textbf{Goal 2.} The second goal of this work, is to carry Vlasov kinetic theory to more general geometric constructions. More precisely, we relax the condition of being autonomous by adding the time variable. In Subsection \ref{sec-kin-cos-mom}, this is achieved on the cosymplectic geometry for momentum variables, and in Subsection \ref{sec-kin-cos-den} the geometry of time-dependent Vlasov theory is derived in terms of the (plasma) density function. The relationship between these two realizations is established through a Poisson/momentum map which is obtained by dualizing a Lie algebra homomorphism.

\noindent
\textbf{Goal 3.} The kinetic equation of particles undergoing contact Hamiltonian dynamics has been studied in \cite{EsGu11} for 3D and in \cite{esen2023conformal} for arbitrary dimensions. Our goal here is to extend this formalism to cocontact geometry to permit explicit time dependence. This is done in Section \ref{sec-kin-coco-mom} for momentum variables and \ref{sec-kin-coco-den} in terms of the plasma density function. As in the cosymplectic case, the connection between these two realizations is established by means of a Poisson/momentum map. 

We create Table \ref{kinetic-algebra-table} for the results of goals 2 and 3 at the algebraic level and the concrete definitions of the dual spaces, and Table \ref{kinetic-table} for the kinetic equations obtained in the present work.

\noindent
\textbf{Contents.} The main body of the paper consists of two sections. In the following section, we shall start with basic expositions of symplectic and cosymplectic manifolds, as well as Hamiltonian, evolution, and gradient vector fields in cosymplectic geometry. Then, we shall present kinetic theories for both momentum and density variables in cosymplectic geometry permitting time dependence. In Section \ref{Sec-Coco}, we shall first present contact geometry and dynamical vector fields. Then, we will review the basics of cocontact geometry. Referring to the dual character of cocontact geometry, we shall present various vector fields and the relationships between them. This section will conclude with geometric kinetic theory on the cocontact framework. Additionally, Appendix \ref{Appapp} is included for the completeness of the work, where one may find necessary background information on Poisson and Jacobi manifolds, as well as geometric kinetic theory in a more abstract framework.
 

\section{Kinetic Theory on Cosymplectic Manifolds}
 
\subsection{Symplectic Manifolds and Dynamics}

A symplectic manifold is a pair $(M,\Omega)$ consisting of a  manifold $M$, and a closed nondegenerate (namely symplectic) two-form $\Omega$. Symplectic manifolds are even-dimensional. For a $2n$-dimensional symplectic manifold the non-degeneracy of the symplectic two-form $\Omega$ is given by
\begin{equation}
    \Omega^n\neq 0. 
\end{equation}
If the symplectic two-form admits a potential one-form as $\Omega=-d\Theta$, then $(M,\Omega)$ is called an exact symplectic manifold. A symplectic manifold is a volume manifold with 
the symplectic volume defined referring to the $n$-th power (assuming that the dimension of $M$ is $2n$) of the symplectic two-form $\Omega$ as
 \begin{equation}\label{symp-volume}
\mu=\frac{(-1)^{n(n-1)/2}}{n!} \Omega^n.
 \end{equation}

\noindent
\textbf{The Canonical Symplectic Structure.} The cotangent manifold $T^*Q$ of a manifold $Q$ is a symplectic manifold equipped with the canonical symplectic two-form $\Omega_Q=-d\Theta_Q$. Here, $\Theta_Q$ is the canonical (Liouville) one-form on $T^*Q$ whose value on 
a vector field defined as 
\begin{equation} \label{canonicaloneform}
\Theta_Q( X) = \langle \tau
_{T^{\ast }Q}\left( X\right) ,T\pi _{Q}\left( X
\right)  \rangle , 
\end{equation}%
for a vector field $X$ on $T^{\ast }Q$.

There exist the Darboux coordinates $(q^i,p_i)$ on $T^{\ast }Q$ (just the natural bundle coordinates) 
\begin{equation}
(q^i,p_i):T^*Q\longrightarrow \mathbb{R}^{2n}
\end{equation}
such that the canonical one form $\Theta_Q$ and the canonical symplectic two-form $\Omega _{Q}$ become
\begin{equation}\label{ss}
\Theta _{Q}=p_i dq^i, \qquad \Omega _{Q}=dq^i \wedge dp_i,
\end{equation}
respectively. the Darboux theorem states that this local structure is generic for all symplectic manifolds.  

\noindent
\textbf{Musical Mappings.} Let $(M,\Omega)$ be a symplectic manifold. 
The non-degenerate character of a symplectic form $\Omega$ gives the following musical flat isomorphism
\begin{equation}\label{flat-map-symp}
    \flat^{s}: \mathfrak{X}(M)\longrightarrow \Gamma^1(M),\qquad X \mapsto \flat^{s}(X)=\iota_X \Omega .
\end{equation}
The inverse is the musical sharp isomorphism $\sharp^{s}:=(\flat^{s})^{-1}$.  
The canonical Poisson bivector field and the musical isomorphism induced from this field are
\begin{equation}\label{canonical-bivector}
\Lambda^{s}(\alpha,\beta):= \Omega(\sharp^{s}(\alpha),\sharp^{s}(\beta)), \qquad 
 \sharp_{\Lambda^{s}}(\alpha): =\Lambda^{s}(\bullet,\alpha)= \sharp^{s}(\alpha),
\end{equation}
respectively. This leads us to define the following canonical Poisson bracket 
\begin{equation}\label{poissonsymplectic-}
 \{ F,H\}^{s} =  \Lambda^{s}(dF,dH).
\end{equation}
In the Darboux coordinates $(q^i,p_i) $ on the symplectic manifold $M$, we compute the musical mapping $\sharp^{s}$ and consequently the musical mapping $\sharp_{\Lambda^{s}}$ in \eqref{canonical-bivector} as
\begin{equation}\label{musical-local-s}
\sharp^{s} = \sharp_{\Lambda^{s}} :\alpha_i dq^i + \alpha^i dp^i 
\mapsto 
\alpha^i \frac{\partial}{\partial q^i}-\alpha_i \frac{\partial}{\partial p_i}.
\end{equation}
The canonical Poisson bivector \eqref{canonical-bivector} and the canonical Poisson bracket \eqref{poissonsymplectic-} are
\begin{equation}\label{poisson}
\Lambda^s=\frac{\partial}{\partial q^i} \wedge \frac{\partial}{\partial p_i},\qquad  
\left\{ F,H\right\}^{s} =\frac{\partial F }{\partial q^{i}}\frac{\partial H}{\partial p_{i}}-\frac{\partial F}{\partial p_{i}}\frac{\partial H}{\partial q^{i}}, 
\end{equation} 
respectively.

\noindent
\textbf{Canonical Diffeomorphisms and Hamiltonian Dynamics.} 
Let $(M,\Omega)$ be a symplectic manifold. The Lie (pseudo) group of all symplectic diffeomorphisms is defined to be 
\begin{equation}
\mathrm{Diff}_{\mathrm{s}}(M)=\{\varphi \in \mathrm{Diff}(M) ~:~ \varphi^*(\Omega)=\Omega\}.
\end{equation} 
We refer \cite{banyaga97,McSa17,RaSc81}  for more details on groups of diffeomorphisms. A vector field is called a symplectic 
 vector field if its flows consists of symplectic diffeomorphisms. These vector fields determine the Lie algebra (of the symplectic diffeomorphisms group) given by 
\begin{equation}
\mathfrak{X}_{\mathrm{s}}(M)=\{X \in \mathfrak{X}(M) ~:~\mathcal{L}_X \Omega = 0\}.
\end{equation} 
For a given Hamiltonian function $H$ on a symplectic manifold $( M,\Omega ) $, the Hamiltonian vector field $X_{H}^{s}$  is defined through  
\begin{equation}
\iota_{X_{H}^{s}}\Omega =dH,  \label{Hamvf}
\end{equation}
where $\iota_ {X^s_H}$ is the contraction operator (the interior derivative). In terms of the musical mappings defined in \eqref{flat-map-symp}, the Hamilton's equation \eqref{Hamvf} can be written as 
\begin{equation}
X_{H}^{s}=\sharp ^{s}(dH).
\end{equation}
In the Darboux coordinates $(q^i,p_i)$, the Hamiltonian vector field on $M$   is computed to be
\begin{equation}\label{X-ham}
X_{H}^{s}=\frac{\partial H}{\partial p_i}\frac{\partial }{\partial q^i }- \frac{\partial H}{\partial q^i }\frac{\partial }{\partial p_i}
\end{equation}
so that the dynamics governed by the Hamiltonian vector field, called Hamilton's equation, is 
\begin{equation}
\dot{q}^i =\frac{\partial H}{\partial p_i},\qquad \dot{p}_i=-\frac{
\partial H}{\partial q^i }.
\end{equation}

Obviously, a Hamiltonian vector field $X_H^s$ is an element of $\mathfrak{X}_{\mathrm{s}}(M)$, and we are interested in the space of Hamiltonian vector fields written as
\begin{equation}\label{alg-sym-ham}
\mathfrak{X}_{\mathrm{s-ham}}(M): =\{X_H^s \in \mathfrak{X}(M) :~  \iota_{X_{H}^{s}}\Omega =dH \}
\end{equation} 
which coincides with $\mathfrak{X}_{\mathrm{s}}(M)$ if the first cohomology class is trivial. 
In this realization, the canonical Poisson bracket is
\begin{equation}\label{poissonsymplectic}
\{F,H\}^{s}= \Omega \left(X_{F}^s,X_{H}^{s}\right)  ,
\end{equation}
where $X_{F}^s$ and $X_{H}^s$ are Hamiltonian vector fields in the sense of (\ref{Hamvf}). Notice that we can write the Poisson bracket via  
\begin{equation}
X_{H}^{s} ( F) = \{ F,H\}^{s}
\end{equation}
as well. 
The following identity 
\begin{equation}\label{Poissonham}
\left[ X_{H}^s,X_{F}^s\right] =-X_{\{ H,F\} ^{s}}.
\end{equation}
manifests that the space of Hamiltonian vector fields $\mathfrak{X}_{\mathrm{s-ham}}(M)$ is a Lie subalgebra of the space $\mathfrak{X}_{\mathrm{s}}(M)$ of symplectic vector fields. .
Further, \eqref{Poissonham} gives that the mapping  
\begin{equation} \label{Phi}
\Phi: C^\infty(M)\longrightarrow \mathfrak{X}_{\mathrm{s-ham}}(M), \qquad H\mapsto X_{H}^{s}
\end{equation} 
is a Lie algebra homomorphism (assuming the Lie bracket on $\mathfrak{X}_{\mathrm{s-ham}}(M)$ is the minus of the Jacobi-Lie bracket of vector fields). This mapping is far from being injective since its kernel is consisting of constant functions.

\subsection{Cosymplectic Manifolds}

A cosymplectic manifold is a triple $(\bar{M},\tau,\Omega)$ consisting of a $(2n + 1)$-dimensional manifold $\bar{M}$, a closed one-form $\tau$, and a closed two-form $\Omega$ satisfying the non-degeneracy condition \cite{Cape,leon89,Liber87}
\begin{equation}
    \tau\wedge\Omega^n\neq 0.
\end{equation} 
There exists the Reeb vector field satisfing the following contractions
\begin{equation}\label{Reeb-cosymp}
\iota_{\mathcal{R}^\tau}\tau=1,\qquad \iota_{\mathcal{R}^\tau}\Omega=0.
\end{equation}
We define the volume form as
\begin{equation}\label{volume-cos}
    \bar{\mu} = \tau\wedge\mu = \frac{(-1)^{n(n-1)/2}}{n!} \tau\wedge\Omega^n.
\end{equation}
 See that $\mu$ is similar to the symplectic volume given in \eqref{symp-volume}.

\noindent
 \textbf{Cosymplectization and Symplectization.} Generic examples of cosymplectic manifolds are extended symplectic manifolds $\mathbb{R}\times M$. We define a two-form $\Omega$ by pulling back the symplectic two-form on $M$ and a one-form $\tau$ by pulling back the differential form $dt$ on $\mathbb{R}$. This reads cosymplectization $(\bar{M},\tau,\Omega)$ of the symplectic manifold $(M,\Omega)$, where $\bar{M}=\mathbb{R}\times M$. This geometrization suggests the Darboux coordinates $(t,q^i,p_i)$ on cosymplectic manifold such that
\begin{equation}\label{canonical-cosym}
\Omega=dq^i\wedge dp_i,\qquad \tau=dt.
\end{equation}
In this local realization, the Reeb vector field is written as
\begin{equation}
    \mathcal{R}^{\tau}=\frac{\partial}{\partial t}.
\end{equation}
Let $(\bar{M},\tau,\Omega)$ be a cosymplectic manifold and consider its extension $\bar{\bar{M}}=\mathbb{R}\times \bar{M}$ and define the two-form \cite{Leon-CosympReduction}
\begin{equation}	\label{Omega-bar}
	\omega=\pi^*\Omega+\pi^*\tau \wedge ds,
\end{equation}	
	 where $\pi$ is the canonical projection from $\bar{\bar{M}}$ to $\bar{M}$ whereas $s$ is the coordinate function on $\mathbb{R}$. Then, a simple calculation shows that $(\bar{\bar{M}}, \omega)$ is a symplectic manifold.
  Notice that $\pi$ is a Poisson morphism.

\noindent
\textbf{Musical Mappings.}
For a cosymplectic manifold, there is an isomorphism
\begin{equation}\label{flat-map}
    \flat^{cs}: \mathfrak{X}(\bar{M})\longrightarrow \Gamma^1(\bar{M}),\qquad X \mapsto \flat^{cs}(X )=\iota_X \Omega + \langle \tau,X \rangle \tau.
\end{equation}
We denote the inverse of this isomorphism by $\sharp^{cs}$. Notice that the inverse image of the one-form $\tau$ is the Reeb field $\mathcal{R}^\tau$.
 A cosymplectic manifold $(\bar{M},\tau,\Omega)$ determines a Poisson bivector field
 \begin{equation}\label{biV-cs}
 \Lambda^{cs} (\alpha,\beta) = \Omega (\sharp^{cs}(\alpha),\sharp^{cs}(\beta))
\end{equation}
 Notice that, in this case, we have 
  \begin{equation}\label{sharp-L-cs}
 \sharp_{\Lambda^{cs}}(\alpha): =\Lambda^{cs}(\bullet,\alpha)= \sharp^{cs}(\alpha) - \langle \alpha, \mathcal{R}^\tau \rangle  \mathcal{R}^\tau.
\end{equation}
The Poisson bracket is then computed to be 
 \begin{equation}\label{Poisson-cs}
    \{F,H\}^{cs}= \Omega (\sharp^{cs}(dF),\sharp^{cs}(dH)).
\end{equation}
 In terms of the Darboux coordinates $(t,q^i,p_i)$, we compute the musical mappings $\sharp^{cs}$ and  $\sharp_{\Lambda^{cs}}$   as
\begin{equation}\label{musical-local-cs}
\begin{split}
\sharp^{cs} &:\alpha_i dq^i + \alpha^i dp^i + 
u dt\mapsto 
\alpha^i \frac{\partial}{\partial q^i}-\alpha_i \frac{\partial}{\partial p_i}  + u \frac{\partial}{\partial t}. 
\\
\sharp_{\Lambda^{cs}} &:\alpha_i dq^i + \alpha^i dp^i + 
u dt\mapsto 
\alpha^i \frac{\partial}{\partial q^i}- \alpha_i  \frac{\partial}{\partial p_i} . 
\end{split}
\end{equation}
The bivector $\Lambda^{cs}$ in \eqref{biV-cs} and the Poisson bracket \eqref{Poisson-cs} 
are computed to be 
\begin{equation}\label{Local-biV-Poisson-cs}
\Lambda^{cs}= \frac{\partial  }{\partial q^{i}}\wedge \frac{\partial  }{\partial p_{i}} , \qquad 
\{F,H\}^{cs}(t,q,p)=\frac{\partial H}{\partial p_i} 
\frac{\partial F}{\partial q^i} 
-
\frac{\partial H}{\partial q^i} 
\frac{\partial F}{\partial p_i}.
\end{equation} 
Even though the cosymplectic Poisson bracket $\{\bullet,\bullet\}^{cs}$ looks the same as the symplectic Poisson bracket in \eqref{poissonsymplectic}, they are not the same. The bracket  \eqref{poissonsymplectic} is non-degenerate while the bracket  \eqref{Poisson-cs} is degenerate by admitting non constant Casimir functions. See that a function depending only on the time variable is a Casimir function for the algebra \eqref{Poisson-cs}.

\subsection{Dynamics  on Cosymplectic Manifolds}

Consider a Hamiltonian function $H$ defined on a cosymplectic manifold $(M,\tau,\Omega)$. We define three vector fields related to the Hamiltonian function. 

We begin defining the group of diffeomorphisms.

\noindent
\textbf{(Weakly) Cosymplectic Diffeomorphisms.}
Let $(\bar{M},\tau,\Omega)$ be a cosymplectic manifold. We define the group of (weakly) cosymplectic diffeomorphisms as \cite{BazzGoer15} 
\begin{equation}\label{Diff-cos}
{\rm Diff}_{\mathrm{cs}} ( \bar{M} ) =\left\{ \varphi \in {\rm Diff} ( \bar{M}   ) :\varphi^{\ast }\Omega=\Omega-dh_\varphi\wedge \tau, ~ \varphi^{\ast }\tau=\tau, \quad h_\varphi\in C^\infty(\bar{M} )\right\} .
\end{equation}
The Lie algebra of this group is the space of (weakly)  cosymplectic vector fields defined to be 
\begin{equation}\label{alg-cos}
\mathfrak{X}_{\mathrm{cs}} (\bar{M}) =\left\{ X\in \mathfrak{X} ( 
\bar{M}) :\mathcal{L}_{X}\Omega =- dh_X\wedge \tau, ~  \mathcal{L}_{X}\tau=0, \quad h_X\in C^\infty(\bar{M} ) \right\} .  
\end{equation}

\noindent
\textbf{(Weakly) Cosymplectic Hamiltonian Vector Fields.}
Consider the Poisson structure $(\bar{M},\Lambda^{cs})$ induced by the cosymplectic manifold $(\bar{M},\tau,\Omega)$. Accordingly, for a given Hamiltonian function $H$ on $\bar{M}$, we define (weakly) cosymplectic Hamiltonian vector field $X^{cs}_H$ as 
  \begin{equation}\label{Poisson-bra-coi-1}
 \{F,H\}^{cs}=\Omega(X^{cs}_F,X^{cs}_H)=X^{cs}_H(F). 
 \end{equation}
where $\{\bullet,\bullet\}^{cs}$ is the cosymplectic Poisson bracket defined in \eqref{Poisson-cs}. In terms of the musical mappings, we compute the (weakly) cosymplectic Hamiltonian vector field as
\begin{equation}\label{Poisson-bra-coi-2--} 
X^{cs}_H  = \sharp_{\Lambda^{cs}}(dH)=\sharp^{cs}(dH)-\langle dH, \mathcal{R}^\tau \rangle \mathcal{R}^\tau. 
\end{equation}
In the Darboux coordinates $(t,q^i,p_i)$, the (weakly) cosymplectic Hamiltonian vector field \eqref{Poisson-bra-coi-2--} is 
  \begin{equation}
X^{cs}_H(t,q,p)=\frac{\partial H}{\partial p_i} 
\frac{\partial }{\partial q^i} 
-
\frac{\partial H}{\partial q^i} 
\frac{\partial }{\partial p_i}.
 \end{equation}
 Notice that, if the real variable $t$ is identified with the time (we will be more clear about this discussion in the upcoming paragraphs where we present evolution vector fields) then  the dynamics determined by the (weakly), the a cosymplectic Hamiltonian vector field represents time dependent Hamiltonian dynamics 
   \begin{equation}
 {\dot q}^i =\frac{\partial H}{\partial p_i}(t,q,p), \qquad
 {\dot p}_i =-\frac{\partial H}{\partial q^i}(t,q,p).
 \end{equation}

It is possible to establish that the following equalities can also be used to determine a (weakly) cosymplectic Hamiltonian vector field uniquely
\begin{equation}\label{Cos-X-H}
\iota_{X^{cs}_H}\Omega =dH- \langle dH,\mathcal{R}^{\tau} \rangle\tau, \qquad \iota_{X^{cs}_H}\tau=0.
\end{equation}
Notice that a (weakly) cosymplectic Hamiltonian vector field is an element of the algebra $\mathfrak{X}_{\mathrm{cs}} (\bar{M})$ of (weakly) cosymplectic vector fields defined in \eqref{alg-cos}. More precisely, for a given Hamiltonian function $H$ the conformal factor $h_X$ in the definition of (weakly) cosymplectic vector field is computed to be $\langle dH,\mathcal{R}^{\tau} \rangle$.  
 In order to see the conservation of the differential forms $(\tau,\Omega)$, we compute the Lie derivatives 
 \begin{equation}
\mathcal{L}_{X^{cs}_H}\tau=0 ,\qquad \mathcal{L}_{X^{cs}_H}\Omega=-d\langle dH, \mathcal{R}^\tau \rangle\wedge \tau.
 \end{equation}
This says that even though the flows of $X^{cs}_H$ do not preserve the two-form $\Omega$, they preserve the volume-form $\bar{\mu}$ in \eqref{volume-cos}. That is, $X^{cs}_H$ is divergence-free: 
  \begin{equation}
  \mathcal{L}_{X^{cs}_H}(\tau\wedge \Omega^n)=  \mathcal{L}_{X^{cs}_H}(\tau) \wedge \Omega^n+ \tau\wedge \mathcal{L}_{X^{cs}_H}(\Omega^n)=-n \tau\wedge d\langle dH, \mathcal{R}^\tau \rangle \wedge\tau\wedge \Omega^{n-1}=0.
   \end{equation}
The Hamiltonian function $H$ is conserved under the flow generated by the Hamiltonian vector field $X^{cs}_H$, that is, $X^{cs}_H(H)=0$.

Furthermore, we have that the Jacobi-Lie bracket of two (weakly) cosymplectic Hamiltonian vector fields is again a (weakly) cosymplectic Hamiltonian vector field by satisfying
 \begin{equation}\label{eksi}
[X^{cs}_F,X^{cs}_H]=-X^{cs}_{\{F,H\}^{cs}}.
\end{equation} 
So, the space of (weakly)  cosymplectic Hamiltonian vector fields turns out to be a Lie subalgebra of $\mathfrak{X}_{\mathrm{cs}} (\bar{M})$. We denote this subalgebra as
\begin{equation}\label{X-cs}
\mathfrak{X}_{\mathrm{cs-ham}} (\bar{M})=\big\{X^{cs}_H\in \mathfrak{X}(\bar{M}): \iota_{X^{cs}_H}\Omega =dH- \langle dH,\mathcal{R}^{\tau} \rangle\tau, \quad \iota_{X^{cs}_H}\tau=0 \big\}.
\end{equation}
We note that the algebraic equation \eqref{eksi} lets us define a Lie algebra homomorphism from the space $C^\infty(\bar{M})$ of smooth functions on the cosymplectic manifold $\bar{M}$ to the space of (weakly) cosymplectic Hamiltonian vector fields as
\begin{equation}\label{Psi}
\bar{\Phi}:C^\infty(\bar{M})\longrightarrow \mathfrak{X}_{\mathrm{cs-ham}}(\bar{M}),\qquad H\mapsto X_H^{cs}.
\end{equation}

\noindent
\textbf{Gradient Vector Field.} Recall that, for a given cosymplectic manifold $(\bar{M},\tau,\Omega)$, we have two musical isomorphisms, namely $\sharp_{\Lambda^{cs}}$ and $\sharp^{cs}$. The former maps the exact differential $dH$ of a Hamiltonian function to the (weakly) cosymplectic Hamiltonian vector field  $X_H^{cs}$. The latter maps $dH$ to the gradient vector field, that is \cite{Leon-GradOnCosymp}
\begin{equation}
\operatorname{grad} H:=\sharp^{cs}(dH).
\end{equation} 
The relationship between the gradient vector field $\operatorname{grad} H$ and the Hamiltonian vector field $X^{cs}_H$ is the following
\begin{equation}\label{Poisson-bra-coi-2} 
\operatorname{grad} H  = X^{cs}_H + \sharp^{cs} (\langle dH, \mathcal{R}^\tau \rangle \tau)  .
\end{equation}
Notice that, alternatively, we can define the gradient vector field as 
\begin{equation}\label{Cos-grad-H}
\iota_{\operatorname{grad} H}\Omega =dH- \langle dH,\mathcal{R}^{\tau} \rangle\tau, \qquad \iota_{\operatorname{grad} H}\tau=\langle dH,\mathcal{R}^{\tau} \rangle.
\end{equation} 
So, we see that the gradient vector field $\operatorname{grad} H$ fails to be an element of the space $\mathfrak{X}_{\mathrm{cs}} (\bar{M}) $ of (weakly) cosymplectic vector fields.  
In the Darboux coordinates $(t,q^i,p_i)$, the gradient vector field is computed to be 
\begin{equation}
\operatorname{grad} H(t,q,p)=\frac{\partial H}{\partial t} 
\frac{\partial }{\partial t} +\frac{\partial H}{\partial p_i} 
\frac{\partial }{\partial q^i} 
-
\frac{\partial H}{\partial q^i} 
\frac{\partial }{\partial p_i},
\end{equation} 
To investigate the conservation properties of the gradient vector field, first we compute the following Lie derivatives
\begin{equation}
 \mathcal{L}_{\operatorname{grad} H}\Omega=-d\langle dH, \mathcal{R}^\tau \rangle\wedge \tau,\qquad \mathcal{L}_{\operatorname{grad} H}\tau=d\langle dH, \mathcal{R}^\tau \rangle.
 \end{equation}
This says that the flows of $\operatorname{grad} H$ do not preserve the one-form $\tau$ and the two-form $\Omega$. Then, we see that
    \begin{equation}
  \mathcal{L}_{\operatorname{grad} H}(\tau\wedge \Omega^n)=  (\mathcal{R}^\tau)^2(H)\tau\wedge \Omega^n.
    \end{equation}
So, the flows of $\operatorname{grad} H$ do not preserve the volume form $\bar{\mu}$ in \eqref{volume-cos} and
\begin{equation}
    \mathrm{div}\operatorname{grad} H = (\mathcal{R}^\tau)^2(H).
\end{equation}
Furthermore, we obtain that
\begin{equation}
   \mathcal{L}_{\operatorname{grad} H} H = \left(\mathcal{R}^\tau(H)\right)^2.
\end{equation}
Therefore, the Hamiltonian function $H$ is not conserved under the flow generated by the vector field $\operatorname{grad} H$.

 \noindent
\textbf{Evolution Vector Field.} It is possible to find other vector fields than the (weakly) cosymplectic Hamiltonian vector field belonging to the space  $\mathfrak{X}_{\mathrm{cs}} (\bar{M}) $ of (weakly) cosymplectic vector fields. For this, we simply integrate the second condition in \eqref{alg-cos} and arrive at a particular local solution  $\iota_{X}\tau=1$. Motivating this, for a Hamiltonian function $H$, we define a evolution vector field $E_H$ by the following identities
\begin{equation}\label{evolution-H}
	\iota_{E_H}\Omega =dH- \langle dH,\mathcal{R}^{\tau} \rangle\tau, \qquad  \iota_{E_H}\tau=1. 
\end{equation}
Notice that, in terms of the musical isomorphism $\sharp^{cs}$, and the Reeb field, the evolution vector field is written 
as 
\begin{equation}\label{evolution-H-alt}
	E_H=\sharp^{cs}(dH) + \mathcal{R}^{\tau} - \langle dH,\mathcal{R}^{\tau} \rangle\mathcal{R}^{\tau}. 
\end{equation}
Also, in this formulation, the first term on the right hand side is the the gradient field $\operatorname{grad} H$, while the first and the third terms read the Hamiltonian vector field $X^{cs}_H$. 
In the Darboux coordinates $(t,q^i,p_i)$, one has that
\begin{equation}
 E_H(t,q,p)=\frac{\partial}{\partial t} + \frac{\partial H}{\partial p_i} 
\frac{\partial }{\partial q^i} 
-
\frac{\partial H}{\partial q^i} 
\frac{\partial }{\partial p_i} 
\end{equation}
and the dynamic generated by the evolution vector field $E_H$ is precisely  
\begin{equation}\label{hamileq22}
 {\dot t} =1,\qquad {\dot q}^i =\frac{\partial H}{\partial p_i}(t,q,p), \qquad
 {\dot p}_i =-\frac{\partial H}{\partial q^i}(t,q,p).
 \end{equation} 
Since $\dot{t}=1$, we can consider $t$ as a time-parameter up to an affine term. So, we can argue that cosymplectic Hamiltonian dynamics is the convenient geometric formulation for time-dependent Hamiltonians.

We compute the Lie derivative of the forms in the direction of the evolution vector field $E_H$
    \begin{equation}
 \mathcal{L}_{E_H}\Omega=-d\langle dH, \mathcal{R}^\tau \rangle\wedge \tau,\qquad \mathcal{L}_{E_H}\tau=0.
    \end{equation}
This says that the flows of $E_H$ do not preserve the two-form $\Omega$, but they preserve the one-form $\tau$. We compute
\begin{equation}
    \mathcal{L}_{E_H}(\tau\wedge \Omega^n)= 0,
\end{equation}
that is,  $E_H$ is divergence-free. On the other hand, we compute that
\begin{equation}
     \mathcal{L}_{E_H}(H) = \mathcal{R}^\tau(H).
\end{equation}
 Thus, the flows of $E_H$ do not preserve the Hamiltonian function $H$.

 We refer the reader to Table \ref{cosymplectic-table} to compare the dynamics created by these three vector fields.

\subsection{Kinetic Lift  Cosymplectic Hamiltonian Dynamics: Momentum Formulation}
\label{sec-kin-cos-mom}

Consider a cosymplectic manifold $(\bar{M},\tau,\Omega)$ and then the space $\mathfrak{X}_{\mathrm{cs-ham}}(\bar{M})$ of  (weakly) cosymplectic Hamiltonian vector fields given in \eqref{X-cs}.     
The dual space $\mathfrak{X}_{\mathrm{cs-ham}}^*(\bar{M})$ of cosymplectic Hamiltonian vector fields  $\mathfrak{X}_{\mathrm{cs-ham}}(\bar{M})$ consists of one-form densities in $\Gamma ^{1}(\bar{M})\otimes \mathrm{Den}(\bar{M})$. To have a more precise definition, we now compute the following pairing 
\begin{equation} \label{dual-calc-C-1}
\begin{split}
\left\langle \bar{\Pi} \otimes \bar{\mu},X^{cs}_H \right\rangle _{L_2} &= \int_{\bar{M}} \langle \bar{\Pi}, X^{cs}_H 
 \rangle   \bar{\mu} =\int_{\bar{M}} \big\langle \bar{\Pi}, \sharp_{\Lambda^{cs}}(dH) \big\rangle \bar{\mu} 
=-\int_{\bar{M}} \big\langle \sharp_{\Lambda^{cs}}(\bar{\Pi}), dH \big\rangle \bar{\mu} 
 \\&
 =-\int_{\bar{M}} \big\langle \sharp^{cs}\bar{\Pi}- \langle \bar{\Pi}, \mathcal{R}^\tau\rangle \mathcal{R}^\tau,  dH  \big\rangle \bar{\mu} 
  \\&
 =  - \int_{\bar{M}}  \langle \sharp^{cs}\bar{\Pi},  dH   \rangle \bar{\mu} + \int_{\bar{M}} \langle \bar{\Pi}, \mathcal{R}^\tau\rangle \langle dH , \mathcal{R}^\tau\rangle \bar{\mu} .
\end{split}
\end{equation} 
We analyse the integrals on the right-hand side one by one. The first integral on the right-hand side can be computed to be
\begin{equation}\label{dual-calc-C-2} 
\begin{split}
-\int_{\bar{M}}  \langle \sharp^{cs}\bar{\Pi},  dH   \rangle \bar{\mu} 
& =-\int_{\bar{M}} \iota_{\sharp^{cs}\bar{\Pi}} dH   \bar{\mu}    =-\int_{\bar{M}} \iota_{\sharp^{cs}\bar{\Pi}} (dH\wedge \bar{\mu} ) -\int_{\bar{M}}  dH\wedge \iota_{\sharp^{cs}\bar{\Pi}} ( \bar{\mu} )
 \\ & =-\int_{\bar{M}} dH \wedge \iota_{\sharp^{cs}\bar{\Pi}}
(\bar{\mu} )
=-\int_{\bar{M}}  d\left( H  \iota_{\sharp^{cs}\bar{\Pi}}\bar{\mu} \right) + \int_{\bar{M}}   H ~ d \iota_{\sharp^{cs}\bar{\Pi}}\bar{\mu} 
 \\ & =
\int_{\bar{M}} H~\mathcal{L}_{\sharp^{cs}\bar{\Pi}} \bar{\mu}   = \int_{\bar{M}} H~ \mathrm{div} \sharp^{cs}\bar{\Pi}  \bar{\mu}
\end{split}
\end{equation} 
up to total divergence. 
The second integral on the right hand side of \eqref{dual-calc-C-1} gives 
\begin{equation}\label{dual-calc-C-3} 
\begin{split}
\int_{\bar{M}} \langle \bar{\Pi}, \mathcal{R}^\tau\rangle \langle dH , \mathcal{R}^\tau\rangle \bar{\mu} &=\int_{\bar{M}} \iota_{\mathcal{R}^\tau}(dH)   \langle \bar{\Pi}, \mathcal{R}^\tau \rangle   \bar{\mu} = \int_{\bar{M}} \langle \bar{\Pi}, \mathcal{R}^\tau \rangle dH \wedge     \iota_{\mathcal{R}^\tau}\bar{\mu}
\\& = - \int_{\bar{M}} H d\langle \bar{\Pi}, \mathcal{R}^\tau \rangle\wedge     \iota_ {\mathcal{R}^\tau}\bar{\mu} -
\int_{\bar{M}} H \langle \bar{\Pi}, \mathcal{R}^\tau \rangle    d \iota_ {\mathcal{R}^\tau}\bar{\mu}
\\& = - \int_{\bar{M}} H \iota_ {\mathcal{R}^\tau }d\langle \bar{\Pi}, \mathcal{R}^\tau \rangle  \bar{\mu} -
\int_{\bar{M}} H \langle \bar{\Pi}, \mathcal{R}^\tau \rangle    \mathcal{L}_ {\mathcal{R}^\tau}\bar{\mu}
\\& = - \int_{\bar{M}} H \mathcal{L}_ {\mathcal{R}^\tau }\langle \bar{\Pi}, \mathcal{R}^\tau \rangle    \bar{\mu} -
\int_{\bar{M}} H \langle \bar{\Pi}, \mathcal{R}^\tau \rangle    \mathrm{div}(\mathcal{R}^\tau)\bar{\mu}
\\& = -\int_{\bar{M}} H \left(\mathcal{L}_ {\mathcal{R}^\tau }\langle \bar{\Pi}, \mathcal{R}^\tau \rangle+
\langle \bar{\Pi}, \mathcal{R}^\tau \rangle    \mathrm{div}(\mathcal{R}^\tau)\right)\bar{\mu},
\end{split}
\end{equation}
up to total divergence. 
A direct computation proves that the vector field $\mathcal{R}^\tau$ is divergence-free. This gives that the second term in the last line of the calculation \eqref{dual-calc-C-3} is zero. So that, combining the calculations (\ref{dual-calc-C-2}) and \eqref{dual-calc-C-3}, we arrive at the following equation
\begin{equation} \label{heye}
\left\langle \bar{\Pi} \otimes \bar{\mu},X^{cs}_H \right\rangle _{L_2} = \int_{\bar{M}} H \left( \mathrm{div} (\sharp^{cs} \bar{\Pi})   - \mathcal{L}_ {\mathcal{R}^\tau }  \langle \bar{\Pi},\mathcal{R}^\tau\rangle     \right)  \bar{\mu}.
\end{equation}
We look for the non-degeneracy of the pairing between the space $\mathfrak{X}_{\mathrm{cs-ham}}(\bar{M})$ of cosymplectic Hamiltonian vector fields and its dual space  $\mathfrak{X}_{\mathrm{cs-ham}}^*(\bar{M})$. So, we ask that the integral on the right hand side of \eqref{heye} should not be zero for an arbitrary Hamiltonian function $H$. This observation gives a more precise definition of the dual space as 
\begin{equation}\label{X-cs-ham}
\mathfrak{X}_{cs-\mathrm{ham}}^*(\bar{M}): =\Big\{\bar{\Pi}
\otimes \bar{\mu} \in \Gamma ^{1}(\bar{M})\otimes \mathrm{Den}(\bar{M}): \mathrm{div} (\sharp^{cs} \bar{\Pi}) -     \mathcal{L}_ {\mathcal{R}^\tau  } \langle \bar{\Pi},\mathcal{R}^\tau\rangle   \neq 0  \Big\} \cup \{0\}.
\end{equation}

Furthermore, one has the dual (which is also Poisson and momentum) mapping 
\begin{equation}\label{den-cos}
\mathfrak{X}_{\mathrm{cs-ham}}^*(\bar{M})\longrightarrow \mathrm{Den}(\bar{M}),\qquad \bar{\Pi} \mapsto f\bar{\mu}=\left(\mathrm{div} (\sharp^{cs} \bar{\Pi})   -\mathcal{L}_{\mathcal{R}^\tau }  \langle \bar{\Pi},\mathcal{R}^\tau\rangle\right) \bar{\mu}.
\end{equation}
In terms of the local coordinates, this map reads 
\begin{equation}
\bar{\Pi}=\bar{\Pi}_t dt + \bar{\Pi}_i dq^i + \bar{\Pi}^i dp_i\mapsto  
   f(t,q,p)=  \frac{\partial \bar{\Pi}^i}{\partial q^i} - \frac{\partial \bar{\Pi}_i}{\partial p_i} .
    \end{equation}
    Fixing the cosymplectic volume, we determine a dual element only by a one-form $\bar{\Pi}$. Since a Hamiltonian vector field is divergence-free, referring to the divergence-free character of Hamiltonian vector fields and the discussion done in Appendix \ref{App-LP-Kinetic}, especially the coadjoint action \eqref{coad-gen}, the coadjoint action of $X^{cs}_H $ on a dual element $\bar{\Pi}$ is computed to be
\begin{equation}
ad^*_{X^{cs}_H }\bar{\Pi}=\mathcal{L}_{X^{cs}_H }\bar{\Pi}.
\end{equation}
So that, the Lie-Poisson dynamics is 
\begin{equation}\label{MV-cos}
\dot{\bar{\Pi}}=-\mathcal{L}_{X^{cs}_H }\bar{\Pi}.
\end{equation}
We call this equation a time-dependent momentum-Vlasov equation. 

\noindent \textbf{Symplectic (Autonomous) Kinetic Dynamics in Momentum Formulation.} As a particular case, we consider now that no explicit time appears in the tensor fields and functions. In this case, we have a symplectic manifold $(M,\Omega)$ and the space of Hamiltonian vector fields $\mathfrak{X}_{\mathrm{ham}}(M)$ defined in \eqref{alg-sym-ham}. The dual space 
$\mathfrak{X}_{\mathrm{ham}}^*(M)$ is \cite{Gu10}
\begin{equation}\label{dual-symp-ham}
\mathfrak{X}_{\mathrm{ham}}^*(M): =\{\Pi
\otimes \mu \in \Gamma ^{1}(M)\otimes \mathrm{Den}(M)):\mathrm{div}\sharp^s(\Pi)\neq 0 \} \cup {0}.
\end{equation}
Accordingly, the dual operation to the Lie algebra homomorphism \eqref{Phi} is  
\begin{equation} \label{phi-*}
\Phi^* : \mathfrak{X}_{\mathrm{ham}}^*(M) \longrightarrow  \mathrm{Den}(M),\qquad \Pi
\otimes \mu\mapsto f  \mu =\mathrm{div}\sharp^s(\Pi)   \mu.
\end{equation}
Being the dual of a Lie algebra homomorphism, $\Phi^*$ is a momentum and a Poisson mapping. In terms of the local coordinates we have that 
\begin{equation}
\Pi=\Pi_i dq^i + \Pi^i dp_i\mapsto     f(q,p)=  \frac{\partial \Pi^i}{\partial q^i} - \frac{\partial \Pi_i}{\partial p_i} .
    \end{equation}

Since we determine the dual space properly, we can determine the Lie-Poisson equations on $\mathfrak{X}_{s-\mathrm{ham}}^*(M)$. A Hamiltonian vector field is divergence-free so, from the calculation \eqref{LP}, we have that 
\begin{equation}\label{MV}
\dot{\Pi}= - ad^*_{X^s_H} \Pi = - \mathcal{L}_{X^s_H} \Pi.
\end{equation}
This is called as momentum-Vlasov equation in the literature \cite{EsGu12,EsSu21,Gu10}.

\subsection{Kinetic Lift  Cosymplectic Hamiltonian Dynamics: Density Formulation} \label{sec-kin-cos-den}

We consider a cosymplectic manifold $(\bar{M},\tau,\Omega)$, and then we start with the following observation  
\begin{equation}\label{associ}
\int_{\bar{M}} \{F,H\}^{cs}K \bar{\mu} = \int_{\bar{M}} F\{H,K\}^{cs} \bar{\mu} 
\end{equation}
where $\bar{\mu}$ is the volume form on the cosymplectic manifold $(\bar{M},\tau,\Omega)$. 
Let us prove this as follows
\begin{equation}\label{symp-bra-comt}
\begin{split}
\int_{\bar{M}} \{F,H\}^{cs}K \bar{\mu} &= \int_{\bar{M}}  X_H(F) K\bar{\mu} =\int_{\bar{M}}  K \left(\iota_{X_H} dF \right)    \wedge  \bar{\mu}
\\
&=\int_{\bar{M}}  K dF \wedge  \iota_{X_H} \bar{\mu}=- \int_{\bar{M}}  F dK  \wedge  \iota_{X_H} \bar{\mu} -
\int_{\bar{M}}  F K  d\iota_{X_H} \bar{\mu}
\\
&= - \int_{\bar{M}}  F \left( \iota_{X_H} dK \right)    \bar{\mu} -
\int_{\bar{M}}  F K  \mathcal{L}_{X_H} \bar{\mu} =  - \int_{\bar{M}}  F X_H (K)    \bar{\mu}
\\
&= -
\int_{\bar{M}} F\{K,H\}^{cs} \bar{\mu} = \int_{\bar{M}} F\{H,K\}^{cs} \bar{\mu}.
\end{split}
\end{equation}

Assuming that the adjoint action on the space of smooth functions $C^\infty(\bar{M})$ is the Poisson bracket $\{\bullet,\bullet\}^{cs}$, then the coadjoint action of the Lie algebra $C^\infty(\bar{M})$ on its dual space $\mathrm{Den}(\bar{M})\simeq C^\infty(\bar{M})$ is 
\begin{equation}
ad^*_Hf= \{f,H\}^{cs}.
\end{equation}
Then, the Lie Poisson equation on the function level is
\begin{equation}\label{LP-cos-density}
\dot{f}=\{H,f\}^{cs}.
\end{equation}
In terms of the Darboux' coordinates $(t,q^i,p_i)$, this dynamics determines the following non-autonomous partial differential equation
\begin{equation}\label{non-Vlasov}
\frac{\partial f}{\partial t}(t,q,p)-\frac{\partial H}{\partial q^i}\frac{\partial f}{\partial p_i}(t,q,p)+ \frac{\partial H}{\partial p_i}\frac{\partial f}{\partial q^i}(t,q,p)=0.
\end{equation}
As in the symplectic case, if the Hamiltonian function $H=p^2/2m+e\phi(q)$ (where $e$ is the charge and $\phi$ is the potential) is the total energy of the single particle in the plasma, then we arrive at the Vlasov equation which permits explicit dependence of time for the plasma density function, that is,   
\begin{equation}\label{V} 
\frac{\partial f}{\partial t}(t,q,p)+\frac{1}{m}p_{i}\frac{\partial f}{\partial
q^{i}}(t,q,p)-e\frac{\partial  \phi }{\partial q^{i}}\frac{\partial f}{\partial
p_{i}}(t,q,p)=0 
\end{equation}

\textbf{From Momentum to Density Formulations.} 
Let us link the density formulation of the dynamics with the dynamics \eqref{MV-cos} in terms of one-forms. We start with the following calculation  
\begin{equation}\label{cos-ss}
\begin{split}
\int_{\bar{M}}\langle ad^*_{{X^{cs}_H}}\bar{\Pi},X^{cs} _K\rangle \bar{\mu} &= \int_{\bar{M}}\langle\bar{\Pi}, ad_{{X^{cs}_H}}X^{cs} _K\rangle \bar{\mu}=-\int_{\bar{M}}\langle\bar{\Pi}, [{X^{cs}_H},X^{cs} _K]\rangle \bar{\mu}
\\
&= \int_{\bar{M}}\langle\bar{\Pi}, X^{cs} _{\{H,K\}^{cs}}\rangle \bar{\mu} 
\\
&= \int \left( \mathrm{div} (\sharp^{cs} \bar{\Pi})   - \mathcal{L}_ {\mathcal{R}^\tau }  \langle \bar{\Pi},\mathcal{R}^\tau\rangle     \right) \{H,K\} ^{cs} \bar{\mu}
\\
&=\int \big\{ \left( \mathrm{div} (\sharp^{cs} \bar{\Pi})   - \mathcal{L}_ {\mathcal{R}^\tau }  \langle \bar{\Pi},\mathcal{R}^\tau\rangle     \right), H\big\}^{cs} K \bar{\mu}
\\
&=\int  \{f,H \}^{cs} K \bar{\mu}
\end{split}
\end{equation}
where we have employed the identity \eqref{eksi} in the second line,  we have used \eqref{heye} in the third line, and we applied \eqref{associ} in the fourth line. Notice that the last line is simply substitution of the density function \eqref{den-cos}.  
On the other hand, by starting with the left hand side of \eqref{cos-ss}, and referring directly to \eqref{heye}, one can deduce
\begin{equation}\label{hess-1}
\begin{split}
\int_{\bar{M}}\langle ad^*_{{X^{cs}_H}}\bar{\Pi},X^{cs} _K\rangle \bar{\mu} &= \int  \mathrm{div} (\sharp^{cs} ad^*_{{X^{cs}_H}}\bar{\Pi})  -\mathcal{L}_ {\mathcal{R}^\tau }  \langle ad^*_{{X^{cs}_H}}\bar{\Pi},\mathcal{R}^\tau\rangle   K \bar{\mu}.
\end{split}
\end{equation}
Collecting up these observations, we have that
\begin{equation}\label{hess-2}
  \mathrm{div} (\sharp^{cs} ad^*_{{X^{cs}_H}}\bar{\Pi})  -\mathcal{L}_ {\mathcal{R}^\tau }  \langle ad^*_{{X^{cs}_H}}\bar{\Pi},\mathcal{R}^\tau\rangle  =
\left \{\left( \mathrm{div} (\sharp^{cs} \bar{\Pi})   -\mathcal{L}_ {\mathcal{R}^\tau }  \langle \bar{\Pi},\mathcal{R}^\tau\rangle     \right),H\right \}^{cs} =\{ f,H\}^{cs}.
\end{equation}
So, one has 
\begin{equation}
\dot{f}=\mathrm{div} (\sharp^{cs} \dot{\bar{\Pi}})   - \mathcal{L}_ {\mathcal{R}^\tau }  \langle  \dot{\bar{\Pi}},\mathcal{R}^\tau\rangle = -\mathrm{div} (\sharp^{cs} \mathcal{L}_{{X^{cs}_H}}\bar{\Pi})   + \mathcal{L}_ {\mathcal{R}^\tau }  \langle  \mathcal{L}_{{X^{cs}_H}}\bar{\Pi},\mathcal{R}^\tau\rangle =\{H,f \}^{cs}.
\end{equation}
As a result, we obtain the nonautonomous Vlasov equation  \eqref{LP-cos-density}. 

 \noindent \textbf{Symplectic (Autonomous) Kinetic Dynamics in Density Formulation.} Once again, we try to investigate the case where there is no explicit time dependence. To have this, we first recall that the following identity holds for the canonical Poisson brackets 
\begin{equation}\label{flip-flow}
\int_M \{F,H\}^{s}K \mu = \int_M F\{H,K\}^{s} \mu 
\end{equation}
for all smooth functions $F$, $H$ and $K$ defined on the symplectic manifold $M$. Here, $\mu$ is the symplectic volume. Assuming that the adjoint action on the space of smooth functions $C^\infty(M)$ is the canonical Poisson bracket, then the coadjoint action of the Lie algebra $C^\infty(M)$ on its dual space $\mathrm{Den}(M)$ is
 \begin{equation}
ad^*_H f =\{f,H\}^{s}
\end{equation}
Let us examine the coadjoint flow in the light of the dual mapping \eqref{phi-*} 
\begin{equation}
\begin{split}
 \langle ad^*_{X_H^{s}}\Pi \otimes \mu ,X_K^{s}\rangle &= \int_M\langle\Pi, ad_{X_H^{s}}X_K^{s}\rangle \mu
=-\int_M\langle\Pi, [X_H^{s},X_K^{s}]  \rangle \mu
\\
&=\int_M\langle\Pi, X^{s}_{\{H,K\}^{s}}\rangle \mu 
=\int \mathrm{div}\sharp^s(\Pi)  \{H,K\}^{s} \mu
\\
&=\int \{ \mathrm{div}\sharp^s(\Pi), H\} ^{s} K \mu.
\end{split}
\end{equation} 
On the other hand, one has that 
\begin{equation}
\begin{split}
\int_M\langle ad^*_{X^{s}_H}\Pi,X^{s}_K\rangle \mu &= \int  \mathrm{div}\sharp^s(ad^*_{X^{s}_H}\Pi)   K \mu.
\end{split}
\end{equation}
Accordingly, we have that
\begin{equation}
\mathrm{div}\left(\sharp^s(\mathcal{L}_{X^{s}_H}\Pi)\right) =
\big \{ \mathrm{div}\sharp^s(\Pi), H \big \}^{s}.
\end{equation}
The relationship between the coadjoint flows reads
\begin{equation}
\dot{f}=\mathrm{div}\sharp^s(\dot{\Pi})=-\mathrm{div}\left(\sharp^s(\mathcal{L}_{X_H^{s}}\Pi)\right)= \left \{H, \mathrm{div}\sharp^s(\Pi)\right \}^{s} = \{H,f\}^s.
\end{equation}
As a result, we have the Lie-Poisson dynamics as 
\begin{equation}
\dot{f} = \{H,f\}^s.
\end{equation}
If the Hamiltonian function $H=p^2/2m+e\phi(q)$  is the total energy of the single particle in the plasma, then we arrive at the Vlasov equation 
\begin{equation}\label{Vla} 
\frac{\partial f}{\partial t}+\frac{1}{m}p_{i}\frac{\partial f}{\partial
q^{i}}-e\frac{\partial  \phi }{\partial q^{i}}\frac{\partial f}{\partial
p_{i}}=0 
\end{equation}
governing collisionless non-relativistic plasma motion. In the literature, one couples the Vlasov equation with Poisson equation
\begin{equation} \label{P} 
\nabla _{q}^{2}\phi _{f}(q)=-e\int f(q,p)d^{3}p 
\end{equation}
to arrive at the integro--differential system \eqref{V}  and \eqref{P}, so called Vlasov-Poisson equations. 
 
\newpage

\section{Kinetic Theory on Cocontact Manifolds}\label{Sec-Coco}

\subsection{Contact Manifolds and Dynamics}

Let $N$ be an odd, say $(2n + 1)$, dimensional manifold. A contact structure on $N$ is a maximally non-integrable smooth distribution of codimension one. A locally contact structure can be determined  by the kernel of a one-form $\eta$ satisfying the non-integrability condition
\begin{equation}  
    d\eta^n
\wedge \eta \neq 0.
\end{equation}
Such a one-form $\eta$ is called a (local) contact form \cite{Arn,Liber87}. In the present paper, we shall consider the existence of a global contact one-form. A (global) contact one-form for a given contact structure is not unique. Indeed, if $\eta$ is a contact one-form for a fixed contact structure, then $\lambda \eta$ also defines the same contact structure for any non-zero real-valued function $\lambda$ defined on $N$. Nevertheless, we call a  $(2n+1)$-dimensional manifold $N$ as a contact manifold if it is equipped with a contact one-form $\eta$, and we shall denote a contact manifold by $(N, \eta)$. Given a contact one-form $\eta$, the vector field $\mathcal{R}^\eta$ satisfying 
\begin{equation}
\iota_{\mathcal{R}^\eta}\eta =1,\qquad \iota_{\mathcal{R}^\eta}d\eta =0
\end{equation}
is unique, and it is called the Reeb vector field. 

\noindent
\textbf{Contactization and Symplectization.} Let us consider an exact symplectic manifold $(M,\Omega=-d\Theta)$ and consider  the principal circle bundle
\begin{equation}\label{contactization}
S^{1}\rightsquigarrow (N,\eta )\overset{pr}{\longrightarrow }(M,\Omega=-d\Theta),
\end{equation}%
called the quantization bundle. For a local coordinate system $z$ on the circle (which may locally be considered to be $\mathbb{R}$), and the Darboux coordinates $(q^i,p_i)$ on the symplectic manifold, the contact manifold admits the Darboux coordinates $(q^i,p_i,z)$ on $N$. In this realization, the contact one-form and the Reeb field are 
\begin{equation}\label{eta-Q}
    \eta=dz-\Theta=dz-p_idq^i,\qquad \mathcal{R}^\eta=\frac{\partial}{\partial z},
\end{equation}
where $\Theta$ (by the abuse of notation) is the pullback of the potential one-form on $M$. One may substitute the symplectic manifold $M$ with the cotangent bundle $T^*Q$. In this case, the contact manifold $N$ locally turns out to be the extended cotangent bundle $T^*Q\times \mathbb{R}$. This also suggests the coordinates $(q^i,p_i,z)$ on the contact manifold $N$. In this case, the volume form $\nu$ on the contact manifold is computed to be 
\begin{equation}\label{contact-volume}
\nu= \mu \wedge\eta
\end{equation}
where $\mu$ (by the abuse of notation) is the pullback of the symplectic volume in \eqref{symp-volume}. 
Accordingly, a generic example of a contact manifold is obtained by the so-called contactization of a symplectic manifold. 
On the other hand, it is possible to determine a symplectic manifold starting with a contact manifold. Concretely, the symplectization of a contact manifold $(N,\eta)$ is the
symplectic manifold 
\begin{equation}\label{con-to-symp}
  (\mathbb{R} \times N ,d(e^t\eta)=e^{t}(d\eta + \tau \wedge 
 \eta )), 
\end{equation}
where $t$ denotes the standard coordinate on $\mathbb{R}$
factor. We remark that, if \eqref{con-to-symp} is a symplectic manifold, then $(M,\eta)$  turns out to be contact.

\textbf{Musical Mappings.} 
For a contact manifold $(N, \eta)$, there is a musical isomorphism 
\begin{equation}\label{flat-map-cont}
\flat^{c}:\mathfrak{X}(N)\longrightarrow \Gamma^1(N),\qquad X\mapsto \iota_Xd\eta+\langle \eta, X \rangle \eta.
\end{equation} 
We denote the inverse of \eqref{flat-map} by $\sharp^{c}$. 
It is worth noting that the image of the Reeb vector field $\mathcal{R}^\eta$ under the musical mapping $\flat^{c}$ is the contact one-form $\eta$. 
Referring to the musical mapping $\sharp^{c}$, we define a bivector field $\Lambda^c$ on  $N$ as
\begin{equation}\label{Lambda}
\Lambda^c(\alpha,\beta)= d\eta(\sharp^{c}\alpha, \sharp^{c} \beta). 
\end{equation}
In the Darboux coordinates, the bivector $\Lambda^c$ is computed to be 
\begin{equation}\label{AP-Bra}
\Lambda^c= \frac{\partial  }{\partial q^{i}}\wedge \frac{\partial  }{\partial p_{i}} + p_i \frac{\partial  }{\partial z} \wedge \frac{\partial  }{\partial p_i}. 
\end{equation} 
Referring to the bivector field $\Lambda^c$ we introduce the following musical mapping 
\begin{equation}\label{Sharp-Delta}
\sharp_{\Lambda^c} (\alpha) :=\Lambda^c(\bullet,\alpha)= \sharp^{c}\alpha - \langle \alpha, \mathcal{R}^\eta \rangle  \mathcal{R}^\eta. 
\end{equation}
 In coordinates, we compute the musical mappings $\sharp^c$ and  $\sharp_{\Lambda^c}$   as
\begin{equation}\label{zap}
\begin{split}
\sharp^c &:\alpha_i dq^i + \alpha^i dp^i + 
\zeta dz\mapsto 
\alpha^i \frac{\partial}{\partial q^i}-(\alpha_i + p_i \zeta)\frac{\partial}{\partial p_i}  + (\zeta + \alpha^ip_i ) \frac{\partial}{\partial z}. 
\\
\sharp_{\Lambda^c} &:\alpha_i dq^i + \alpha^i dp^i + 
\zeta dz\mapsto 
\alpha^i \frac{\partial}{\partial q^i}-(\alpha_i + p_i \zeta)\frac{\partial}{\partial p_i}  + \alpha^ip_i  \frac{\partial}{\partial z}. 
\end{split}
\end{equation}
The pair $(M,\Lambda^c)$ is an almost Poisson manifold since the bracket 
\begin{equation}\label{AP-Bracket}
    \{F,H\}^{(APc)}=\Lambda^c(dF,dH)
\end{equation}
satisfies the Leibniz identity but fails to satisfy the Jacobi identity. Then the almost Poisson bracket \eqref{AP-Bracket} turns out to be 
\begin{equation}\label{AP-Bra-c}
\{F,H\}^{(APc)} (q,p,z) = \frac{\partial F}{\partial q^{i}}\frac{\partial H}{\partial p_{i}} -
\frac{\partial F}{\partial p_{i}}\frac{\partial H}{\partial q^{i}} + p_{i}
\frac{\partial F}{\partial z}
\frac{\partial H}{\partial p_{i}}   - p_{i}\frac{\partial F}{\partial p_{i}}  \frac{\partial H}{\partial z} .
\end{equation}

On the other hand the pair $(\Lambda^c,\mathcal{R}^\eta)$ induces a Jacobi structure \cite{Lichnerowicz-Jacobi,Marle-Jacobi}, since the equalities
\begin{equation}
[\Lambda^c,\Lambda^c]= 2\mathcal{R}^\eta\wedge \Lambda^c, \qquad [\mathcal{R}^\eta,\Lambda^c]=0,
\end{equation}
hold, where the bracket is the Schouten–Nijenhuis bracket(we have employed the definition in  \cite{Marsden1999}). This realization permits us to define a Jacobi (contact) bracket 
\begin{equation}\label{cont-bracket-L} 
\{F,H\}^c =\Lambda^c(dF,dH) +F\mathcal{R}^\eta(H) - H\mathcal{R}^\eta(F).
\end{equation}
Notice that this bracket satisfies the Jacobi identity but not the Leibniz identity. In this local picture, the contact bracket \eqref{cont-bracket-L} is 
\begin{equation}\label{Lag-Bra}
\{F,H\}^c  = \frac{\partial F}{\partial q^{i}}\frac{\partial H}{\partial p_{i}} -
\frac{\partial F}{\partial p_{i}}\frac{\partial H}{\partial q^{i}} + \left(F  - p_{i}\frac{\partial F}{\partial p_{i}} \right)\frac{\partial H}{\partial z} -
\left(H  - p_{i}\frac{\partial H}{\partial p_{i}} \right)\frac{\partial F}{\partial z}.
\end{equation}

\noindent \textbf{Contact Diffeomorphisms.}
We  define the
group of contact diffeomorphisms by \cite{banyaga97}
\begin{equation}\label{Diff-c}
{\rm Diff}_{\mathrm{c}} ( N) =\left\{ \varphi \in {\rm Diff} ( N ) :\varphi ^{\ast }\eta =- h_\varphi   \eta, \quad h_\varphi \in C^\infty%
( N) \right\} .
\end{equation}
Notice that the existence of $h_\varphi$ in the definition manifests the conformal invariance of the contact structure. 
A vector field on the contact
manifold $( N,\eta) $ is a contact vector field (called also infinitesimal conformal contact diffeomorphism) if
it generates a one-parameter group of contact diffeomorphisms. Accordingly, the space of contact vector fields is given by
\begin{equation}\label{algcon}
\mathfrak{X}_{\mathrm{c}} ( N  ) =\left\{ X\in \mathfrak{X} ( 
N ) :\mathcal{L}_{X}\eta  =-h_X \eta ,\quad h_X
\in C^\infty ( N   ) \right\} .  
\end{equation}
Sometimes, a contact vector field is denoted by a pair $(X,h_X)$ to exhibit the conformal factor $h_X$.   

\noindent
\textbf{Contact Hamiltonian Vector Fields.} Consider a contact manifold $(N,\eta)$. 
Referring to its Jacobi structure $(N,\Lambda^c,\mathcal{R}^\eta)$, and the Hamiltonian vector description in \eqref{Ham-v-f-Jac}, in terms of the musical isomorphisms $\sharp_{\Lambda^c}$ and $\sharp^c$,  we define the contact Hamiltonian vector field $X_{H}^c$ as \cite{Br17,BrCrTa17,de2021hamilton,de2019contact,esen2021contact}
\begin{equation}\label{3-def}
X_{H}^c=\sharp_{\Lambda^c}(dH) - H\mathcal{R}^\eta= \sharp^c(dH)-(\langle dH, \mathcal{R}^\eta \rangle +H)\mathcal{R}^\eta.
\end{equation}
The following identities also define the contact Hamiltonian vector field
\begin{equation}\label{Ham-def} 
\iota_{X_{H}}d\eta=dH-\mathcal{R}^\eta(H)\eta,\qquad \iota_{X_{H}}\eta=-H.
\end{equation} 
The relationship between the contact bracket in \eqref{cont-bracket-L} and the contact Hamiltonian vector field is given by 
\begin{equation}\label{3-def-}
X_{H}^c(F)=\{F,H \}^c -  F\mathcal{R}^\eta(H).
\end{equation}
This observation follows the same line of understanding in \eqref{bra-X-J} given for the generic Jacobi framework.
Referring to the Darboux coordinates $(q^{i},p_{i},z)$, for a function $H=H(q,p,z)$, the contact Hamiltonian vector field becomes
\begin{equation}\label{con-dyn-X}
X_{H}^c (q,p,z)=\frac{\partial H}{\partial p_{i}}\frac{\partial}{\partial q^{i}}  - \left(\frac{\partial H}{\partial q^{i}} + \frac{\partial H}{\partial z} p_{i} \right)
\frac{\partial}{\partial p_{i}} + \left(p_{i}\frac{\partial H}{\partial p_{i}} - H\right)\frac{\partial}{\partial z}.
\end{equation}
Thus, we obtain the contact Hamilton's equations for $H$ as
\begin{equation}\label{conham-X}
\dot{q}^i = \frac{\partial H}{\partial p_{i}}, 
\qquad
\dot{p}_i  = -\frac{\partial H}{\partial q^{i}}- 
p_{i}\frac{\partial H}{\partial z}, 
\qquad \dot{z}= p_{i}\frac{\partial H}{\partial p_{i}} - H.
\end{equation}

A direct computation determines a conformal factor for a given contact Hamiltonian vector field via
\begin{equation}\label{L-X-eta}
\mathcal{L}_{X_{H}^c }\eta =
d\iota_{X_{H}^c }\eta+\iota_{X_{H}^c }d\eta= -\mathcal{R}^\eta(H)\eta.
\end{equation}
According to \eqref{L-X-eta}, the flow of a contact Hamiltonian system preserves the contact structure. So, a contact Hamiltonian vector field $X_{H}^c $ belongs to the space $\mathfrak{X}_{\mathrm{c}} ( N  )$ of contact vector fields in \eqref{algcon}.  Still, it does not preserve either the contact one-form or the Hamiltonian function. Instead, we obtain
\begin{equation}
{\mathcal{L}}_{X_{H}^c } \, H = - \mathcal{R}^\eta(H) H.
\end{equation}
Being a non-vanishing top form, we may consider $d\eta^n
\wedge \eta$ as a volume form on $N$.
The Hamiltonian motion does not preserve the volume form since
\begin{equation}
{\mathcal{L}}_{X_{H}^c }  \, (d\eta^n
\wedge \eta) = - (n+1)  \mathcal{R}^\eta(H) d\eta^n
\wedge \eta.
\end{equation}
Assuming the dimension of $N$ to be $2n+1$, we compute the divergence of a contact vector field as
\begin{equation} \label{div-X-H}
\mathrm{div}(X_{H}^c )= -  (n+1)  \mathcal{R}^\eta(H).
\end{equation}

Further, due to the identity \eqref{Jac-integ}, or alternatively, after a direct calculation, we have that 
\begin{equation} \label{integ-contact}
\left[X^c_{F},X^c_{H} \right]=-X^c_{ \{ F,H \}^{c}}  .
\end{equation}
Collecting the observations in \eqref{L-X-eta} and \eqref{integ-contact}, we define the space $\mathfrak{X}_{\mathrm{c-ham}} ( N  )$ of contact Hamiltonian vector fields as a Lie subalgebra  \cite{Br17,BrCrTa17,de2021hamilton,LeLa19}
\begin{equation}\label{X-con}
\mathfrak{X}_{\mathrm{c-ham}} ( N  )=\big\{X_{H}^c \in \mathfrak{X}(N  ): \iota_{X_{H}^c }\eta =-H,~ \iota_{X_{H}^c }d\eta =dH-\mathcal{R}^\eta(H) \eta\big\}
\end{equation}
of the space $\mathfrak{X}_{\mathrm{c}} (N)$ of contact vector fields \eqref{algcon}.   Moreover, one may establish the following
isomorphism from the space of real smooth functions on $N$ to the space of contact Hamiltonian vector fields equipped with minus of the Jacobi-Lie bracket as the Lie algebra structure, that is, 
\begin{equation} \label{iso}
\Psi: ( C^\infty( N ) ,\left\{ \bullet,\bullet \right\} ^{c})\longrightarrow 
\left( \mathfrak{X}_{\mathrm{c-ham}} ( N  ) ,-\left[\bullet ,\bullet
\right] \right) , \qquad H \mapsto X_{H}^c .
\end{equation}

\noindent \textbf{Energy Preserving (Evolution) Dynamics.} In the previous paragraph, we have examined the Hamiltonian dynamics for a given contact manifold $(N,\eta)$ referring to its Jacobi manifold structure $(N,\Lambda^c,\mathcal{R}^\eta)$. Now, we examine Hamiltonian dynamics on a contact manifold referring to its almost Poisson character $(N,\{\bullet,\bullet\}^{APc})$ due to the almost Poisson bracket in \eqref{AP-Bracket}. Accordingly, taking in mind the definition in \eqref{Poisson-ham}, we define energy (or contact evolution) vector field\footnote{We prefer to call this vector field as energy field to avoid a duplication of the term "evolution" which is already coined to time determining vector field \eqref{evolution-H} in cosymplectic category.}   $\mathcal{E}_{H}$ as \cite{Houches}
\begin{equation}
\mathcal{E}_{H}= \sharp_{\Lambda^c}(dH)=\sharp^c(dH)-\langle dH, \mathcal{R}^\eta\rangle \mathcal{R}^\eta,
\end{equation}
where $\sharp_{\Lambda^c}$ is the musical isomorphism induced from the almost Poisson bivector $\sharp^c$. We remark that the energy vector field is preserving the energy that is 
\begin{equation}
\mathcal{E}_{H}(H)= \sharp_{\Lambda^c}(dH)(dH)=0
\end{equation}
as a manifestation of the skew-symmetry of the almost Poisson bivector field $\Lambda^c$. 
It is possible to obtain the following two conditions
\begin{equation}\label{evo-def} 
\iota_{\mathcal{E}_{H}}d\eta=dH-\mathcal{R}^\eta(H)\eta,\qquad \iota_{\mathcal{E}_{H}}\eta=0
\end{equation}  
 are enough to define the energy vector field for a given Hamiltonian function $H$.
A contact energy vector field $\mathcal{E}_{H}$ preserves the Hamiltonian function $H$, but $\mathcal{E}_{H}$
does not preserve the contact structure. In the Darboux coordinates, the contact energy vector field is computed to be
\begin{equation}\label{evo-dyn}	\mathcal{E}_{H}(q,p,z)=\frac{\partial H}{\partial p_i}\frac{\partial}{\partial q^i}  - \left (\frac{\partial H}{\partial q^i} + \frac{\partial H}{\partial z} p_i \right)
	\frac{\partial}{\partial p_i} + p_i\frac{\partial H}{\partial p_i} \frac{\partial}{\partial z},
\end{equation}
so that the integral curves satisfy the energy (evolution) equations
\begin{equation}\label{evo-eq-c}
	\frac{d q^i}{dt}= \frac{\partial H}{\partial p_i}, \qquad \frac{d p_i}{dt}  = -\frac{\partial H}{\partial q^i}- 
	p_i\frac{\partial H}{\partial z}, \quad \frac{d z}{dt} = p_i\frac{\partial H}{\partial p_i}.
\end{equation}
Notice that the divergence of the contact energy vector field $\mathcal{E}_{H}$ is computed to be
\begin{equation} \label{div-E-H-}
\mathrm{div}(\mathcal{E}_{H})= -  n  \mathcal{R}^\eta(H) =  - n \frac{\partial H}{\partial z}.
\end{equation}

 \noindent \textbf{Quantomorphisms.} Let us consider a contact manifold $(N,\eta)$.  By asking the conformal factor $h_\varphi$ to be the unity for a contact diffeomorphism $\varphi$ in \eqref{Diff-c}, one arrives at the conservation of the contact one-form $\varphi^*\eta = \eta$. We call such a mapping as a strict contact diffeomorphism (or a quantomorphism). For a contact manifold $(N,\eta)$, 
we denote the group of all strict contact diffeomorphisms  as
\begin{equation}
{\rm Diff}_{\mathrm{q}} (N)
 =\left\{ \varphi \in {\rm Diff} ( N) :\varphi ^{\ast }\eta =   \eta \right\} \subset {\rm Diff}_{\mathrm{c}} ( N)
.
\end{equation}
The Lie algebra of this Lie group consists of vector fields whose flows are strict contact diffeomorphisms. We call such vector fields as strict contact vector fields (or, infinitesimal quantomorphisms, or infinitesimal strict contact diffeomorphism). We denote the space as 
\begin{equation}\label{X-qua}
\mathfrak{X}_{\mathrm{q}}  ( N) =\left\{ X \in \mathfrak{X} ( N ) :\mathcal{L}_{X}\eta =0\right\} \subset \mathfrak{X}_{\mathrm{\mathrm{c}} }( N ). 
\end{equation}

 \noindent \textbf{Strict Contact Hamiltonian Vector Fields.} For a given Hamiltonian function $H$, the contact Hamiltonian vector field $X^{c}_H$ is a strict contact vector field if and only if $\mathcal{R}^\eta(H)=0$. Accordingly, one has the following space of strict contact Hamiltonian vector fields
\begin{equation}\label{X-qua-ham}
\mathfrak{X}_{\mathrm{q-ham}} ( N )=\big\{\xi_{H}\in \mathfrak{X}(N ): \iota_{\xi_{H}}\eta =-H,~ \iota_{\xi_{H}}d\eta =dH \big\}.
\end{equation}
Notice that one has the following identity in terms of the musical mapping $\flat^{c}$ in \eqref{flat-map} 
\begin{equation}
\xi_{H}= \sharp^c(dH) - H \mathcal{R}^\eta.
\end{equation}

In local coordinates $(q^i,p_i,z)$, the condition $\mathcal{R}^\eta(H)=0$ implies that $H$ must not depend on the $z$ variable.
Referring to the Darboux coordinates $(q^{i},p_{i},z)$, for a Hamiltonian function $H=H(q,p)$ independent of the fiber variable $z$, the strict contact Hamiltonian vector field is 
\begin{equation}\label{con-dyn-xi}
\xi_{H}(q,p)=\frac{\partial H}{\partial p_{i}}\frac{\partial}{\partial q^{i}}  -  \frac{\partial H}{\partial q^{i}} 
\frac{\partial}{\partial p_{i}} + \left(p_{i}\frac{\partial H}{\partial p_{i}} - H\right)\frac{\partial}{\partial z}.
\end{equation}
Thus, we obtain strict contact Hamilton's equations as 
\begin{equation}\label{conham-xi}
\dot{q}^i = \frac{\partial H}{\partial p_{i}}, 
\qquad
\dot{p}_i  = -\frac{\partial H}{\partial q^{i}}, 
\qquad \dot{z}= p_{i}\frac{\partial H}{\partial p_{i}} - H.
\end{equation}
One can see that this flow is divergence-free.

Assuming the contactization bundle $N\mapsto M$, for two functions that are not dependent on the fiber variable $z$, the contact bracket $\{\bullet ,\bullet \}^{c}$ in \eqref{cont-bracket-L} locally possesses a similar form to the canonical Poisson bracket on $M$. Accordingly, a direct calculation reads that 
\begin{equation}
\left[ \xi_{H},\xi_{F}\right]=-\xi_{ \{ H,F \}^{c}}.
\end{equation}
We assume the contactization $\pi:N\mapsto M$ of a contact manifold $(N,\eta)$ over the symplectic manifold $(M,\Omega)$ as displayed in \eqref{contactization}. 
This motivates us to define an isomorphism from the space $C^\infty(M)$ of smooth functions defined on M to the space $\mathfrak{X}_{\mathrm{q-ham}}(N)$ of strict contact Hamiltonian vector fields on the contact manifold $N$. So we have that 
\begin{equation} \label{iso-}
\zeta: ( C^\infty(M) ,\left\{ \bullet,\bullet \right\}^s )\longrightarrow 
\left( \mathfrak{X}_{\mathrm{q-ham}} ( N) ,-\left[\bullet ,\bullet
\right] \right) , \qquad H \mapsto X^c_{\pi^*H}.
\end{equation} 
We present Table \ref{contact-table} for a list of vector fields on contact manifolds.

\subsection{Cocontact Manifolds}

 An even-dimensional, say $(2n+2)$, manifold $ \bar{N} $, is said to be a cocontact manifold if it admits a one form $\eta$  and a closed one-form $\tau$ satisfying the non-degeneracy condition
 \begin{equation}\label{non-int-coco}
     \tau \wedge d\eta ^n \wedge \eta \neq 0. 
 \end{equation}
 We denote a cocontact bundle with a triplet $(\bar{N},\tau,\eta)$. 
The condition \eqref{non-int-coco} reads that the one-form $\eta$ should not be closed for $n\geq 1$. In this case, the top form $\tau \wedge d\eta ^n \wedge \eta$ gives rise to a volume form  
\begin{equation}\label{cocontact-volume}
\overline{\nu} =\frac{(-1)^{n(n-1)/2}}{n!} \tau\wedge  \eta \wedge (d\eta)^n.
\end{equation}
There are two distinguished vector fields on the cocontact manifold. The first is the cosymplectic Reeb field $\mathcal{R}^\tau$ uniquely defined by the conditions 
 \begin{equation}\label{R-coss}
  \iota_{\mathcal{R}^\tau}\tau = 1, \qquad \iota_{\mathcal{R}^\tau}\eta = 0,\qquad \iota_{\mathcal{R}^\tau} d\eta = 0, 
  \end{equation}
and the contact Reeb field $\mathcal{R}^\eta$ defined by the conditions 
 \begin{equation}\label{R-cont}
  \iota_{\mathcal{R}^\eta}\tau = 0, \qquad \iota_{\mathcal{R}^\eta}\eta = 1,\qquad \iota_{\mathcal{R}^\eta} d\eta = 0.
  \end{equation} 

\noindent  \textbf{From Old to New.}
We can define cocontact manifold starting with a contact manifold or an (exact) cosymplectic manifold. Let us examine these constructions one by one. First, we start with a contact manifold $(N,\eta)$ and consider a line bundle $\bar{N}$ over $N$ depicted as 
 \begin{equation}\label{con-coco}
\mathbb{R}\rightsquigarrow (\bar{N}=\mathbb{R}\times N,\tau,\eta) \overset{pr}{\longrightarrow } (N,\eta).
\end{equation}
Here, $\tau$ is assumed to be the exact differential $dt$ where $t$ is the standard coordinate for the extension $\mathbb{R}$. We call this procedure cocontactization from a contact manifold. Secondly, we start with an exact cosymplectic manifold $(\bar{M},\tau,\Omega=-d\Theta)$, and then we consider the extension of 
 \begin{equation}\label{cos-coco}
\mathbb{R}\rightsquigarrow (\bar{N}=\bar{M}\times \mathbb{R},\tau,\eta)  \overset{pr}{\longrightarrow } (\bar{M},\tau,\Omega=-d\Theta).
\end{equation}
Here, $\eta$ is taken to be $dz-\Theta$ where $z$ is the standard coordinate on $\mathbb{R}$. We call this procedure cocontactization of a cosymplectic manifold. Merging these two constructions, it is possible to arrive at a $(2n+2)$-dimensional cocontact manifold starting with a $2n$-dimensional exact symplectic manifold. This needs to define a line bundle structure on $M$ iteratively. More concretely, 
 \begin{equation} \label{sym-coco}
\mathbb{R}^2\rightsquigarrow (\bar{N}=\mathbb{R}\times M \times \mathbb{R},\tau,\eta)  \overset{pr}{\longrightarrow }  (M,\Omega=-d\Theta).
\end{equation}
Here, we consider the coordinates on $\mathbb{R}^2$ as $(z,t)$ and then $\tau$ is assumed to be the exact differential $dt$ while $\eta$ is taken to be $dz-\Theta$. 
The following picture demonstrates all these extensions
\begin{equation}
\xymatrix{
& \bar{M} \times \mathbb{R} =\bar{N} = \mathbb{R}  \times N\ar[dl] \ar[dr]
\\\bar{M}=\mathbb{R}\times M \ar[dr]&&N=M\times \mathbb{R} \ar[dl]
\\&M\ar@/_1pc/[ur]_{\text{\quad contactization}}\ar@/^1pc/[ul]^{\text{cosymplectization\quad }}} \label{T}
\end{equation} 
where each of the downward arrows represents line bundles while upward arrows stand for cosymplectization, contactization, etc.
An exact cocontact manifold $(\bar{N},\tau,\eta)$ (that is the case $\tau=dt$) admits symplectic character as well if it is equipped with the symplectic two-form \eqref{con-to-symp}. So, we can argue that $(\bar{N}=\mathbb{R} \times N,\tau,\eta)$ is an exact cocontact manifold if and only if $(\mathbb{R} \times N ,d(e^t\eta) )$ is a symplectic manifold. 

Now, we start at the bottom of the picture \eqref{T} and go up to the top to obtain the Darboux coordinates on each of these geometries. Accordingly, we start with $2n$-dimensional symplectic manifold $M$ with the Darboux coordinates $(q^i,p_i)$. Then we define the Darboux coordinates $(t,q^i,p_i)$ on the cosymplectic manifold $\bar{M}$ whereas we take $(q^i,p_i,z)$  on the contact manifold $N$. At last, we have $(t,q^i,p_i,z)$ as the coordinates on the cocontact manifold $\bar{N}$. In this realization, we have 
\begin{equation}
\tau=dt, \qquad \eta=dz-p_i dq^i. 
\end{equation} 
The Reeb vector fields read
\begin{equation}
\mathcal{R}^\tau = \frac{\partial }{\partial  t}\ ,\quad \mathcal{R}^\eta = \frac{\partial }{\partial z}. 
\end{equation} 

\noindent
\textbf{Musical Mappings.
} Consider a cocontact manifold $(\bar{N},\tau,\eta)$. The non-degeneracy condition \eqref{non-int-coco} reads the following flat isomorphism from the space of sections of the tangent bundle $T\bar{N}$ to the sections of the cotangent bundle $T^*\bar{N}$ defined by 
\begin{equation}\label{flat-cc}
    \flat^{cc}: \mathfrak{X}(\bar{N} ) \longrightarrow  \Gamma^1(\bar{N}) ,\qquad X \mapsto \langle\tau,X\rangle\tau + \iota_{X} d\eta + \langle\eta,X\rangle\eta.
\end{equation}
The inverse of the flat isomorphism is the sharp isomorphism $\sharp^{cc} = (\flat^{cc})^{-1}$. 
In terms of the musical mappings, we have the Reeb fields as 
   \begin{equation}
 \mathcal{R}^\tau = (\flat^{cc})^{-1}(\tau), \qquad \mathcal{R}^\eta = (\flat^{cc})^{-1}(\eta).
    \end{equation} 
We define a bivector field referring to the sharp map $\sharp^{cc}$ as 
    \begin{equation} \label{Lambda-cc}
\Lambda^{cc}(\alpha,\beta) = d\eta(\sharp^{cc}\alpha,\sharp^{cc}\beta) 
      \end{equation}
In coordinates, we have 
\begin{equation}
\Lambda ^{cc} (t,q,p,z) = \frac{\partial}{\partial q^i}\wedge\frac{\partial}{\partial p_i} +p_i \frac{\partial}{ \partial z} \wedge \frac{\partial}{\partial p_i}.
\end{equation}

 The bivector $\Lambda^{cc}$ induces a musical map
\begin{equation}\label{sharp-L-cc}
 \sharp_{\Lambda^{cc}}: \Gamma^1(\bar{N})\longrightarrow \mathfrak{X}(\bar{N}),\qquad \alpha\mapsto \sharp_{\Lambda^{cc}}(\alpha)= \Lambda^{cc}(\bullet, \alpha).
 \end{equation} 
It can be seen that the one-forms $\eta$ and $\tau$ belong to the kernel of the musical mapping $\sharp_{\Lambda^{cc}}$.  
 A direct calculation provides the following identity
 \begin{equation}
 \sharp_{\Lambda^{cc}} \alpha =\sharp^{cc}\alpha - \langle \alpha, \mathcal{R}^\eta \rangle \mathcal{R}^\eta - \langle \alpha, \mathcal{R}^\tau \rangle \mathcal{R}^\tau. 
\end{equation}
In the Darboux coordinates, the musical maps are computed to be
\begin{equation}\label{zap-cc}
\begin{split}
\sharp^{cc} &:\alpha_i dq^i + \alpha^i dp_i + 
\zeta dz + udt \mapsto 
\alpha^i \frac{\partial}{\partial q^i}-(\alpha_i + p_i \zeta)\frac{\partial}{\partial p_i}  + (\zeta + \alpha^ip_i ) \frac{\partial}{\partial z} + u \frac{\partial}{\partial t}
\\
\sharp_{\Lambda^{cc}} &:\alpha_i dq^i + \alpha^i dp_i + 
\zeta dz  + udt \mapsto 
\alpha^i \frac{\partial}{\partial q^i}-(\alpha_i + p_i \zeta)\frac{\partial}{\partial p_i}  + \alpha^ip_i  \frac{\partial}{\partial z}. 
\end{split}
\end{equation}

The two-tuple $(\bar{N},\Lambda^{cc})$ is an almost Poisson manifold since the bracket 
\begin{equation}\label{AP-Bracket-cc}
    \{F,H\}^{(APcc)}:=\Lambda^{cc}(dF,dH)
\end{equation}
satisfies the Leibniz identity but fails to satisfy the Jacobi identity. In the Darboux coordinates, the almost Poisson bracket \eqref{AP-Bracket-cc}  turns out to be 
\begin{equation}\label{AP-Bra-cont}
\{F,H\}^{(APcc)} (t,q,p,z) = \frac{\partial F}{\partial q^{i}}\frac{\partial H}{\partial p_{i}} -
\frac{\partial F}{\partial p_{i}}\frac{\partial H}{\partial q^{i}} + p_{i}
\frac{\partial F}{\partial z}
\frac{\partial H}{\partial p_{i}}   - p_{i}\frac{\partial F}{\partial p_{i}}  \frac{\partial H}{\partial z} .
\end{equation}

For a cocontact manifold $(\bar{N},\tau,\eta)$ the pair  $(\Lambda^{cc},\mathcal{R}^\eta)$ determines a Jacobi structure since
\begin{equation}
[\Lambda^{cc},\Lambda^{cc}]= 2\mathcal{R}^\eta\wedge \Lambda^{cc}, \qquad [\mathcal{R}^\eta,\Lambda^{cc}]=0,
\end{equation}
where the bracket is the Schouten–Nijenhuis bracket. We have employed the bracket definition in  \cite{Marsden1999}. This realization permits us to define a Jacobi (cocontact) bracket 
\begin{equation}\label{cocont-bracket-L} 
\{F,H\}^{cc} =\Lambda^{cc}(dF,dH) +F\mathcal{R}^\eta(H) - H\mathcal{R}^\eta(F).
\end{equation}
Notice that the bracket satisfies the Jacobi identity but not the Leibniz identity. In this local picture, the cocontact bracket \eqref{cocont-bracket-L} is 
\begin{equation}\label{Lag-Bra-coco}
\{F,H\}^{cc} (t,q,p,z) = \frac{\partial F}{\partial q^{i}}\frac{\partial H}{\partial p_{i}} -
\frac{\partial F}{\partial p_{i}}\frac{\partial H}{\partial q^{i}} + \left(F  - p_{i}\frac{\partial F}{\partial p_{i}} \right)\frac{\partial H}{\partial z} -
\left(H  - p_{i}\frac{\partial H}{\partial p_{i}} \right)\frac{\partial F}{\partial z}.
\end{equation}

\subsection{Cocontact Hamiltonian Dynamics}\label{sec-coco-ham}

Consider a Hamiltonian function $H$ defined on a cocontact manifold $(\bar{N},\tau,\eta)$. We define nine vector fields related to the Hamiltonian function. We begin defining the group of diffeomorphisms.

\noindent \textbf{Cocontact Diffeomorphisms.}
We  define the
group of (weakly) cocontact diffeomorphisms by 
\begin{equation}\label{Diff-cc}
{\rm Diff}_{\mathrm{cc}} ( \bar{N} ) =\left\{ \varphi \in {\rm Diff} ( \bar{N} ) :\varphi ^{\ast }\eta = - h_\varphi\eta - g_\varphi\tau, ~ \varphi^{\ast }\tau=\tau, ~ h_\varphi, g_\varphi \in C^\infty%
( \bar{N})\right\} .
\end{equation}
A vector field on the cocontact
manifold  is called a (weakly) cocontact vector field if it generates a one-parameter group of (weakly) cocontact diffeomorphisms. This definition provides us the following  Lie algebra
\begin{equation}\label{X-cc}
\mathfrak{X}_{\mathrm{cc}} ( \bar{N} ) =\left\{ X\in \mathfrak{X} ( 
N ) :\mathcal{L}_{X}\eta  =-h_X \eta - g_X \tau , ~\mathcal{L}_{X} \tau=0, ~ h_X,g_X 
\in C^\infty ( \bar{N}) \right\} .  
\end{equation}
of cocontact vector fields.

\noindent 
\textbf{Cocontact Hamiltonian Vector Fields.}
For a cocontact manifold $(\bar{N},\tau,\eta)$, consider the associated Jacobi structure. Thus, equation \eqref{Ham-v-f-Jac} determines the definition of cocontact Hamiltonian vector field as
\begin{equation}\label{murat}
\begin{split}
    X_{H}^{cc} &=  \sharp_{\Lambda^{cc}}(dH) -H\mathcal{R}^\eta
    \\
   & = \sharp^{cc}(dH)-\langle dH, \mathcal{R}^\tau \rangle \mathcal{R}^\tau-(\langle dH, \mathcal{R}^\eta \rangle +H)\mathcal{R}^\eta 
\end{split}
\end{equation}
in terms of the musical isomorphism $\sharp_{\Lambda^{cc}}$ in \eqref{sharp-L-cc} and $\sharp^{cc}$. 
The relationship between the cocontact bracket in \eqref{cocont-bracket-L} and the cocontact Hamiltonian vector field is given by 
\begin{equation}\label{3-def--}
X_{H}^{cc}(F)=\{F,H \}^{cc} -  F\mathcal{R}^\eta(H).
\end{equation}
As in the contact case, this is in synchronization with \eqref{bra-X-J} given for the abstract Jacobi framework.

We take the Lie derivative of the differential forms $\tau$, $\eta$, and $d\eta$ in the direction of the cocontact Hamiltonian vector field $X_{H}^{cc}$ to measure the rate of change in the fundamental objects:
\begin{equation}\label{L-X-eta-cc}
\begin{split}
    \mathcal{L}_{X_{H}^{cc}}\tau & = 0,
\\
    \mathcal{L}_{X_{H}^{cc}}\eta & =  -\langle dH, \mathcal{R}^\eta \rangle \eta-\langle dH, \mathcal{R}^\tau \rangle\tau,
\\
    \mathcal{L}_{X_{H}^{cc}}d\eta & = -d \langle dH, \mathcal{R}^\eta \rangle \wedge \eta - \langle dH, \mathcal{R}^\eta \rangle d\eta -d \langle dH, \mathcal{R}^\tau \rangle\wedge  \tau.
\end{split}
\end{equation}
From the first two equations above, we see that a cocontact Hamiltonian vector field preserves the one-form $\tau$ whereas it does not preserve the one-form $\eta$. From the last equation, we see that if a vector field is in the kernel of the one-form $\tau$ and the kernel of the one-form $\eta$ then it belongs to the kernel of $\mathcal{L}_{X_{H}^{cc}}\eta$ as well.
To sum up, for an $(2n+2)$ dimensional cocontact manifold, we compute the rate of change of the top form as 
\begin{equation}\label{L-X-vol-cc}
\begin{split}
    \mathcal{L}_{X_{H}^{cc}} ( \tau \wedge d\eta ^n \wedge \eta ) = -(n+1)\mathcal{R}^\eta(H) \tau \wedge d\eta ^n \wedge \eta .
    \end{split}
\end{equation}
This shows that the Hamiltonian motion does not preserve the volume form. From the calculation \eqref{L-X-vol-cc}, we deduce that the divergence of the cocontact Hamiltonian vector field $X_{H}^{cc}$ is
\begin{equation} \label{div-X-H-cc}
\mathrm{div}(X_{H}^{cc} )=  -(n+1)\mathcal{R}^\eta(H).
\end{equation}
From equation \eqref{murat} and using the skew-symmetry of the almost Poisson bivector field $\Lambda^{cc}$, we compute the rate of change of the Hamiltonian function $H$ along the Hamiltonian motion as follows
\begin{equation}
{\mathcal{L}}_{X_{H}^{cc} } \, H =  -\mathcal{R}^\eta(H)H.
\end{equation}
Therefore, the Hamiltonian function is not preserved along the motion. 

There is an alternatively and equivalent way to define a cocontact Hamiltonian vector field $X_H^{cc}$ defined on a cocontact manifold $(\bar{N},\tau,\eta)$. It is as follows  
	\begin{equation} 
			   \iota_{X_H^{cc}}d \eta = d H-\langle dH, \mathcal{R}^\eta \rangle \eta-\langle dH, \mathcal{R}^\tau \rangle\tau, \qquad 
			\iota_{X_H^{cc}}\eta = -H, \qquad 
			\iota_{X_H^{cc}}\tau = 0. 
	\end{equation}
This realization of cocontact vector field $X_H^{cc}$  is important since it becomes immediate now to see that $X_H^{cc}$ belongs to the Lie algebra \eqref{X-cc} of cocontact vector fields. Notice that for the cocontact Hamiltonian vector field $X_H^{cc}$  one has the conformal functions appearing in the definition \eqref{X-cc}  as
	\begin{equation} 
h_X =\langle dH, \mathcal{R}^\eta \rangle , \qquad g_X= \langle dH, \mathcal{R}^\tau \rangle .
	\end{equation}

Using \eqref{Jac-integ}, it is possible to prove that the space of cocontact Hamiltonian vector fields is closed under the Lie bracket of vector fields as
\begin{equation} \label{jac-cc}
\left[X^{cc}_{F},X^{cc}_{H} \right]=-X^{cc}_{ \{ F,H \}^{cc}}  .
\end{equation}
Thus, the space of cocontact Hamiltonian vector fields is a Lie subalgebra of the space $\mathfrak{X}_{\mathrm{cc}} ( \bar{N} )$ of cocontact vector fields. So, we define the following
\begin{equation}\label{X-cc-ham}
\mathfrak{X}_{\mathrm{cc-ham}} ( \bar{N})=\big\{X_{H}^{cc} \in \mathfrak{X}(\bar{N}):  \iota_{X_H^{cc}}d \eta = d H-\langle dH, \mathcal{R}^\eta \rangle \eta-\langle dH, \mathcal{R}^\tau \rangle\tau, ~
			\iota_{X_H^{cc}}\eta = -H, ~
			\iota_{X_H^{cc}}\tau = 0.\big\}.
\end{equation}
Furthermore, the equality \eqref{jac-cc} implies an 
isomorphism from the space of real smooth functions on $\bar{N}$ to the space of cocontact Hamiltonian vector fields
\begin{equation} \label{bar-psi}
\bar{\Psi}: ( C^\infty( \bar{N} ) ,\left\{ \bullet,\bullet \right\} ^{cc})\longrightarrow 
\left( \mathfrak{X}_{\mathrm{cc-ham}} ( \bar{N}  ) ,-\left[\bullet ,\bullet
\right] \right) , \qquad H \mapsto X_{H}^{cc} .
\end{equation}

Let us close the previous discussion around cocontact Hamiltonian vector fields by writing the vector field and the dynamics in terms of the  the Darboux coordinates $(t,q^{i},p_{i},z)$ on the cocontact manifold $\bar{N}$. So, for a given Hamiltonian function $H=H(t,q,p,z)$, the cocontact Hamiltonian vector field is computed to be
\begin{equation}\label{con-dyn-}
X_{H}^{cc} (t,q,p,z) =\frac{\partial H}{\partial p_{i}}\frac{\partial}{\partial q^{i}}  - \left(\frac{\partial H}{\partial q^{i}} + \frac{\partial H}{\partial z} p_{i} \right)
\frac{\partial}{\partial p_{i}} + \left(p_{i}\frac{\partial H}{\partial p_{i}} - H\right)\frac{\partial}{\partial z}.
\end{equation}
Thus, we obtain the cocontact Hamilton's equations for $H$ as
\begin{equation}\label{conham-}
\dot{q}^i = \frac{\partial H}{\partial p_{i}}, 
\qquad
\dot{p}_i  = -\frac{\partial H}{\partial q^{i}}- 
p_{i}\frac{\partial H}{\partial z}, 
\qquad \dot{z}= p_{i}\frac{\partial H}{\partial p_{i}} - H.
\end{equation}

\noindent \textbf{Energy Preserving Vector Fields.}  In the cocontact Hamiltonian vector field definition, we have used the Jacobi structure $(N,\Lambda^{cc},\mathcal{R}^\eta)$ of a contact manifold $(N,\eta)$. Now, we examine Hamiltonian dynamics on a cocontact manifold referring to its almost Poisson character $(N,\{\bullet,\bullet\}^{APcc})$ due to the almost Poisson bracket defined in \eqref{AP-Bracket-cc}. Accordingly, taking in mind the definition in \eqref{Poisson-ham}, we define the cocontact energy vector field  $\mathcal{E}_{H}^{cc}$ as \cite{Houches}
\begin{equation}\label{E-H-cc}
\mathcal{E}_{H}^{cc}= \sharp_{\Lambda^{cc}}(dH)=\sharp^{cc}(dH) 
- \langle dH, \mathcal{R}^\eta \rangle \mathcal{R}^\eta - \langle dH, \mathcal{R}^\tau \rangle \mathcal{R}^\tau
.
\end{equation}
where $\sharp_{\Lambda^{cc}}$ is the musical isomorphism induced from the almost Poisson bivector $\Lambda^{cc}$.

Let us consider once more a cocontact manifold $(\bar{N},\tau,\eta)$. For a given Hamiltonian function $H$ on the cocontact manifold, a cocontact energy vector field $\mathcal{E}_{H}^{cc}$ is defined to be 
\begin{equation}\label{E-H-cc-alt} 
\iota_{\mathcal{E}_{H}^{cc}}d\eta=dH-\langle dH, \mathcal{R}^\eta \rangle \eta-\langle dH, \mathcal{R}^\tau \rangle\tau,\qquad \iota_{\mathcal{E}_{H}^{cc}}\eta=0, \qquad 
			\iota_{\mathcal{E}_{H}^{cc}}\tau = 0.
   \end{equation}  
Algebraically, we can state that the cocontact energy  vector fields are defined referring to the almost Poisson structure of cocontact manifold $(\bar{N},\tau,\eta)$ equipped with the almost Poisson bivector $\Lambda^{cc}$ in \eqref{Lambda-cc}. Referring to the musical mapping $\sharp_{\Lambda^{cc}}$ in \eqref{sharp-L-cc} induced by this almost Poisson bivector, we write the cocontact energy vector field as 
\begin{equation}
\mathcal{E}_{H}^{cc}= \sharp_{\Lambda^{cc}}(dH)=\sharp^{cc}(dH) 
- \langle dH, \mathcal{R}^\eta \rangle \mathcal{R}^\eta - \langle dH, \mathcal{R}^\tau \rangle \mathcal{R}^\tau
.
\end{equation}
Accordingly, a cocontact energy vector field fails to satisfy the Jacobi identity. In the Darboux coordinates, the cocontact energy vector field is computed to be
\begin{equation}\label{evo-dyn-cc}	\mathcal{E}_{H}^{cc}(t,q,p,z)=\frac{\partial H}{\partial p_i}\frac{\partial}{\partial q^i}  - \left (\frac{\partial H}{\partial q^i} + \frac{\partial H}{\partial z} p_i \right)
	\frac{\partial}{\partial p_i} + p_i\frac{\partial H}{\partial p_i} \frac{\partial}{\partial z},
\end{equation}
so that the dynamics is 
\begin{equation}\label{evo-eq-cc}
	\frac{d q^i}{dt}= \frac{\partial H}{\partial p_i}, \qquad \frac{d p_i}{dt}  = -\frac{\partial H}{\partial q^i}- 
	p_i\frac{\partial H}{\partial z}, \quad \frac{d z}{dt} = p_i\frac{\partial H}{\partial p_i}.
\end{equation}
To see the conservation of the volume form, we first compute the Lie derivatives in the direction of the cocontact energy vector field $\mathcal{E}_{H}^{cc}$ as
\begin{equation}\label{L-E-eta-cc}
\begin{split}
    \mathcal{L}_{\mathcal{E}_{H}^{cc}}\tau & = 0,
\\
    \mathcal{L}_{\mathcal{E}_{H}^{cc}}\eta & = dH - \langle dH,\mathcal{R}^\eta \rangle \eta - \langle dH,\mathcal{R}^\tau \rangle \tau,
\\
    \mathcal{L}_{\mathcal{E}_{H}^{cc}}d\eta & = - d\langle dH,\mathcal{R}^\eta \rangle \wedge \eta - \langle dH,\mathcal{R}^\eta \rangle d\eta - d\langle dH,\mathcal{R}^\tau \rangle \wedge \tau.
\end{split}
\end{equation}
The first identity shows that a cocontact evolution vector preserves the one-form $\tau$ whereas the second one shows that it does not preserve the one-form $\eta$.
Using the identities in \eqref{L-E-eta-cc}, we obtain
\begin{equation}\label{L-E-eta-cc-vol}
\begin{split}
    \mathcal{L}_{\mathcal{E}_{H}^{cc}} ( \tau \wedge d\eta ^n \wedge \eta ) = -n\mathcal{R}^\eta(H)\tau \wedge d\eta ^n \wedge \eta.
    \end{split}
\end{equation}
Therefore we deduce that the cocontact evolution vector does not preserve the volume form. The divergence of a cocontact evolution vector field is determined as
\begin{equation} \label{div-E-H}
\mathrm{div}(\mathcal{E}_{H}^{cc})= -n\mathcal{R}^\eta(H).
\end{equation}
On the other hand, by computing the Lie derivative of the Hamiltonian function $H$ in the direction of the cocontact energy vector field  $\mathcal{E}_{H}^{cc}$ we obtain
\begin{equation}
{\mathcal{L}}_{\mathcal{E}_{H}^{cc}} \, H =  0,
\end{equation}
so, we deduce that the flows of the vector field $\mathcal{E}_{H}^{cc}$ preserves the Hamiltonian function $H$.

\noindent \textbf{Coquantomorphims (Strict Cocontactomorphisms).} Let us consider a cocontact manifold $(\bar{N},\tau,\eta)$.  By asking the conformal factor $h_\varphi$ to be the unity for a cocontact diffeomorphism in \eqref{Diff-cc}, we define strict cocontact diffeomorphism (or a coquantomorphism). For a cocontact manifold $(\bar{N},\tau,\eta)$,
we denote the group of all strict cocontact diffeomorphisms  as 
\begin{equation}\label{Diff-qcc}
{\rm Diff}_{\mathrm{qcc}} ( \bar{N} ) =\left\{ \varphi \in {\rm Diff} ( \bar{N} ) :\varphi ^{\ast }\eta = \eta - g_\varphi\tau, ~ \varphi^{\ast }\tau=\tau, ~  g_\varphi \in C^\infty%
( \bar{N})\right\} .
\end{equation}
Vector fields generating strict cocontact diffeomorphisms are called strict cocontact vector fields which determine the Lie algebra
\begin{equation}\label{X-qcc}
\mathfrak{X}_{\mathrm{qcc}} ( \bar{N} ) =\left\{ X\in \mathfrak{X} ( 
N ) :\mathcal{L}_{X}\eta  = - g_X \tau , ~\mathcal{L}_{X} \tau=0, ~  g_X 
\in C^\infty ( N   ) \right\} .  
\end{equation}

\noindent 
\textbf{Strict Cocontact Hamiltonian Vector Fields.}
For a given Hamiltonian function $H$ on the cocontact manifold, a cocontact Hamiltonian vector field $X^{cc}_H$ is a strict contact vector field if and only if $\mathcal{R}^\eta(H)=0$. This locally proves that a Hamiltonian function must be independent on $z$-variable. So, we have the space of strict cocontact Hamiltonian Vector Fields defined as follows
\begin{equation}\label{X-qcc-ham}
\mathfrak{X}_{\mathrm{qcc-ham}} ( \bar{N})=\big\{\xi_{H}^{cc} \in \mathfrak{X}(\bar{N}):  \iota_{\xi_{H}^{cc} }d \eta = d H-\langle dH, \mathcal{R}^\tau \rangle\tau, ~
			\iota_{\xi_{H}^{cc} }\eta = -H, ~
			\iota_{\xi_{H}^{cc} }\tau = 0\big\}.
\end{equation}  
Notice that one has the following identity in terms of the musical mapping $\flat^{cc}$ in \eqref{flat-cc}  as 
\begin{equation}
\xi_{H}^{cc} = \sharp^{cc}(dH)-\langle dH, \mathcal{R}^\tau \rangle \mathcal{R}^\tau-H \mathcal{R}^\eta .
\end{equation}
Referring to the Darboux coordinates $(t,q^{i},p_{i},z)$, for a Hamiltonian function $H=H(t,q,p)$ independent of the fiber variable $z$, the strict cocontact Hamiltonian vector field is 
\begin{equation}\label{con-dyn}
\xi_{H}^{cc} (t,q,p)=\frac{\partial H}{\partial p_{i}}\frac{\partial}{\partial q^{i}}  -  \frac{\partial H}{\partial q^{i}} 
\frac{\partial}{\partial p_{i}} + \left(p_{i}\frac{\partial H}{\partial p_{i}} - H\right)\frac{\partial}{\partial z}.
\end{equation}
Thus, we obtain strict cocontact Hamilton's equations as  
\begin{equation}\label{conham}
\dot{q}^i = \frac{\partial H}{\partial p_{i}}, 
\qquad
\dot{p}_i  = -\frac{\partial H}{\partial q^{i}}, 
\qquad \dot{z}= p_{i}\frac{\partial H}{\partial p_{i}} - H.
\end{equation}

Let us see the rate of change of the volume form in the direction of the strict contact vector field $\xi_{H}^{cc}$. First we compute the following identities
\begin{equation}\label{L-xi-eta-cc}
\begin{split}
    \mathcal{L}_{\xi_{H}^{cc}}\tau & = 0,
\\
    \mathcal{L}_{\xi_{H}^{cc}}\eta & = - \langle dH,\mathcal{R}^\tau \rangle \tau,
\\
    \mathcal{L}_{\xi_{H}^{cc}}d\eta & = - d\langle dH,\mathcal{R}^\tau \rangle \wedge \tau.
\end{split}
\end{equation}
We see that $\xi_{H}^{cc}$ preserves the one-form $\tau$, but it does not preserve the one-form $\eta$. Moreover, the second equaion implies that if a vector field is in the kernel of the one-form $\tau$, then it belongs to the kernel of $\mathcal{L}_{X_{H}^{cc}}\eta$ as well.
The identities in \eqref{L-xi-eta-cc} allow us to compute the Lie derivative of the volume form
\begin{equation}\label{L-xi-eta-cc-div}
\begin{split}
    \mathcal{L}_{\xi_{H}^{cc}} ( \tau \wedge d\eta ^n \wedge \eta ) = 0.
    \end{split}
\end{equation}
So, we conclude that strict cocontact Hamiltonian vector field is divergence-free.
We compute the Lie derivative of the Hamiltonian function $H$ as
\begin{equation}
{\mathcal{L}}_{\xi_{H}^{cc}} \, H = 0.
\end{equation}
So, the flows of $\xi_{H}^{cc}$ preserve the Hamiltonian function $H$. 

A direct computation shows the following identity 
\begin{equation}
\left[ \xi_{H}^{cc} ,\xi_{F}^{cc}\right]=-\xi_{ \{ H,F \}^{cc}},
\end{equation}
which will be used to learn more about strict cocontact Hamiltonian vector fields. Consider the cocontactization $\bar{\pi}:\bar{N}\mapsto \bar{M}$ of a cocontact manifold $(\bar{N},\tau,\eta)$ over the cosymplectic manifold $(\bar{M},\tau,\Omega)$. 
This motivates us to define an isomorphism from the space $C^\infty(\bar{M})$ of smooth functions defined on $\bar{M}$ to the space $\mathfrak{X}_{\mathrm{qcc-ham}}(\bar{N})$ of strict cocontact Hamiltonian vector fields on the contact manifold $\bar{N}$. So, we have that 
\begin{equation} \label{iso--}
\bar{\zeta}: ( C^\infty(\bar{M}) ,\left\{ \bullet,\bullet \right\}^s )\longrightarrow 
\left( \mathfrak{X}_{\mathrm{qcc-ham}} ( \bar{N}) ,-\left[\bullet ,\bullet
\right] \right) , \qquad H \mapsto X^c_{\bar{\pi}^*H}.
\end{equation}

 \subsection{Cocontact Evolution Dynamics}\label{sec-coco-evo}

Consider a cocontact manifold $(\bar{N},\tau,\eta)$. 
Recall the algebra of cocontact vector fields introduced in \eqref{X-cc}. In the case of cocontact Hamiltonian vector fields \eqref{murat}, one integrates the second condition $\mathcal{L}_{X} \tau=0$ as $\iota_{X} \tau=0$ (see the third condition in \eqref{X-cc-ham}). In this subsection we integrate the condition $\mathcal{L}_{X} \tau=0$  as $\iota_{X} \tau=1$. This assumption relates $t$ and the time variables in an affine way. In this subsection, we shall examine the three vector fields defined in Subsection \ref{sec-coco-ham} one by one, replacing the condition $\iota_{X} \tau=0$ by this other one $\iota_{X} \tau=1$.

\noindent 
\textbf{Cocontact Evolution Vector Fields.} We define a cocontact evolution vector field $\bar{X}_H^{cc}$ as   
	\begin{equation} 
		   \iota_{\bar{X}_H^{cc}}d \eta = d H-\langle dH, \mathcal{R}^\eta \rangle \eta-\langle dH, \mathcal{R}^\tau \rangle\tau, \qquad 
			\iota_{\bar{X}_H^{cc}}\eta = -H, \qquad 
			\iota_{\bar{X}_H^{cc}}\tau = 1. 
	\end{equation}
In terms of the musical mappings $\sharp^{cc}$ and $\sharp_{\Lambda^{cc}}$ in \eqref{sharp-L-cc} we write  the cocontact evolution vector field  as 
\begin{equation}
\begin{split}
\bar{X}_H^{cc} &= \sharp_{\Lambda^{cc}}(dH) + \mathcal{R}^\tau- H  \mathcal{R}^\eta  
\\ & =\sharp^{cc}(dH) 
- \left( H + \langle dH, \mathcal{R}^\eta \rangle  \right) \mathcal{R}^\eta + \left(1 - \langle dH, \mathcal{R}^\tau \rangle \right)\mathcal{R}^\tau.
\end{split}
\end{equation}
In the Darboux coordinates $(t,q^{i},p_{i},z)$, for a given Hamiltonian function $H=H(t,q,p,z)$, the cocontact evolution vector field turns out to be
\begin{equation}\label{con-dyn-Xbar}
\bar{X}_H^{cc}(t,q,p,z) =\frac{\partial }{\partial t}+ \frac{\partial H}{\partial p_{i}}\frac{\partial}{\partial q^{i}}  - \left(\frac{\partial H}{\partial q^{i}} + \frac{\partial H}{\partial z} p_{i} \right)
\frac{\partial}{\partial p_{i}} + \left(p_{i}\frac{\partial H}{\partial p_{i}} - H\right)\frac{\partial}{\partial z}.
\end{equation}
Thus, we obtain the dynamics for $H$ as
\begin{equation}\label{conham-Xbar}
\dot{t} =1, \qquad 
\dot{q}^{i} = \frac{\partial H}{\partial p_{i}}, 
\qquad
\dot{p}_{i} = -\frac{\partial H}{\partial q^{i}}- 
p_{i}\frac{\partial H}{\partial z}, 
\qquad \dot{z}  = p_{i}\frac{\partial H}{\partial p_{i}} - H.
\end{equation}
Given $\dot{t}=1$, we can treat $t$ as a time-parameter up to an affine term. Therefore, one can see that cocontact evolution dynamics offers a nice geometric framework for describing time-dependent dissipative systems.

The Lie derivatives of the forms $\tau,\eta$, and $d\eta$ produce the identities
\begin{equation}\label{L-Xbar-eta-cc}
\begin{split}
    \mathcal{L}_{\bar{X}_H^{cc}}\tau & = 0,
\\
    \mathcal{L}_{\bar{X}_H^{cc}}\eta & = - \langle dH,\mathcal{R}^\eta \rangle \eta - \langle dH,\mathcal{R}^\tau \rangle \tau,
\\
    \mathcal{L}_{\bar{X}_H^{cc}}d\eta & = - d\langle dH,\mathcal{R}^\eta \rangle \wedge \eta - \langle dH,\mathcal{R}^\eta \rangle d\eta - d\langle dH,\mathcal{R}^\tau \rangle \wedge \tau,
\end{split}
\end{equation}
respectively. Thus, although a cocontact evolution vector field preserves the one-form $\tau$, it does not preserve the one-form $\eta$. 

Following the computations in \eqref{L-Xbar-eta-cc}, the Lie derivative of the volume form is 
\begin{equation}\label{L-Xbar-eta-cc-vol}
\begin{split}
    \mathcal{L}_{\bar{X}_H^{cc}} ( \tau \wedge d\eta ^n \wedge \eta ) = -(n+1)\mathcal{R}^\eta(H) \tau \wedge d\eta ^n \wedge \eta. 
    \end{split}
\end{equation}
This formula shows that a cocontact evolution vector field does not preserve the volume form. More concretely, the divergence of the cocontact evolution vector field becomes
\begin{equation} \label{div-Xbar-H}
\mathrm{div}(\bar{X}_H^{cc})= -(n+1)\mathcal{R}^\eta(H).
\end{equation}
Furthermore, by a direct computation of the Lie derivative of the Hamiltonian function $H$ in the direction of the cocontact evolution vector field, we have
\begin{equation}
{\mathcal{L}}_{\bar{X}_H^{cc}} \, H = \mathcal{R}^\tau(H) -\mathcal{R}^\eta(H)H.
\end{equation}
Thus, the cocontact evolution vector field $\bar{X}_H^{cc}$ does not preserve the Hamiltonian function $H$.

\noindent \textbf{Energy Evolution Dynamics.}  Now, consider the almost Poisson structure $(\bar{N},\Lambda^{cc})$ and the determined energy preserving cocontact dynamics in \eqref{E-H-cc-alt}. We replace the condition $\iota_{X} \tau=0$  by $\iota_{X} \tau=1$,
and arrive at the definition of energy evolution vector field $\bar{\mathcal{E}}_{H}^{cc}$ as 
\begin{equation}\label{E-H-cc-alt--} 
\iota_{\bar{\mathcal{E}}_{H}^{cc}}d\eta=dH-\langle dH, \mathcal{R}^\eta \rangle \eta-\langle dH, \mathcal{R}^\tau \rangle\tau,\qquad \iota_{\bar{\mathcal{E}}_{H}^{cc}}\eta=0, \qquad 
			\iota_{\bar{\mathcal{E}}_{H}^{cc}}\tau = 1.
   \end{equation}  
Considering the musical mapping $\sharp_{\Lambda^{cc}}$, we write the energy evolution vector field as 
\begin{equation*}
   \bar{\mathcal{E}}_{H}^{cc}= \sharp_{\Lambda^{cc}}(dH)+ \mathcal{R}^\tau  =\sharp^{cc}(dH) -
\langle dH, \mathcal{R}^\eta \rangle  \mathcal{R}^\eta + \left( 1- \langle dH,  \mathcal{R}^\tau \rangle \right)\mathcal{R}^\tau
. 
\end{equation*}
In the Darboux coordinates, the energy evolution vector field is 
\begin{equation}\label{evo-dyn-coco}
\bar{\mathcal{E}}_{H}^{cc}(t,q,p,z)=\frac{\partial }{\partial t}+\frac{\partial H}{\partial p_i}\frac{\partial}{\partial q^i}  - \left (\frac{\partial H}{\partial q^i} + \frac{\partial H}{\partial z} p_i \right)
	\frac{\partial}{\partial p_i} + p_i\frac{\partial H}{\partial p_i} \frac{\partial}{\partial z},
\end{equation}
so that, the dynamics is
\begin{equation}\label{evo-eq-}
\dot{t} =1,\qquad 
	\frac{d q^i}{dt}= \frac{\partial H}{\partial p_i}, \qquad \frac{d p_i}{dt}  = -\frac{\partial H}{\partial q^i}- 
	p_i\frac{\partial H}{\partial z}, \quad \frac{d z}{dt} = p_i\frac{\partial H}{\partial p_i}.
\end{equation}
To investigate the conservation of the volume form, we compute
\begin{equation}\label{L-Ebar-eta-cc}
\begin{split}
    \mathcal{L}_{\bar{\mathcal{E}}_{H}^{cc}}\tau & = 0,
\\
    \mathcal{L}_{\bar{\mathcal{E}}_{H}^{cc}}\eta & = dH - \langle dH,\mathcal{R}^\eta \rangle \eta - \langle dH,\mathcal{R}^\tau \rangle \tau,
\\
    \mathcal{L}_{\bar{\mathcal{E}}_{H}^{cc}}d\eta & = - d\langle dH,\mathcal{R}^\eta \rangle \wedge \eta - \langle dH,\mathcal{R}^\eta \rangle d\eta - d\langle dH,\mathcal{R}^\tau \rangle \wedge \tau.
\end{split}
\end{equation}
From the first two equations, we see that an energy evolution vector field preserves the one-form $\tau$  but it does not preserve the one-form $\eta$. Using the identities in \eqref{L-Ebar-eta-cc}, we calculate the Lie derivative in the direction of the energy evolution vector field as
\begin{equation}\label{L-Ebar-eta-cc-vol}
\begin{split}
    \mathcal{L}_{\bar{\mathcal{E}}_{H}^{cc}} ( \tau \wedge d\eta ^n \wedge \eta ) = -n\mathcal{R}^\eta(H) \tau \wedge d\eta ^n \wedge \eta
    \end{split}
\end{equation}
and deduce that an energy evolution vector field does not preserve the volume form. This shows that the divergence of the energy evolution vector field is
\begin{equation} \label{div-Ebar-H}
\mathrm{div}(\bar{\mathcal{E}}_{H}^{cc})= -n\mathcal{R}^\eta(H).
\end{equation}

We observe that the Hamiltonian function $H$ is not preserved since
\begin{equation}
{\mathcal{L}}_{\bar{\mathcal{E}}_{H}^{cc}} \, H = \mathcal{R}^\tau(H).
\end{equation}
For contact dynamics, such a vector field preserves the energy as a result of the skew-symmetry of the almost Poisson bivector field $\Lambda^{c}$. In the cocontact case, the energy vector field, given in \eqref{E-H-cc-alt--}, does not preserve energy, but we still refer to it as energy since the global definition resembles that of contact energy vector fields. 

 \noindent 
\textbf{Strict Evolution Vector Fields.}
For a given Hamiltonian function $H$ on the cocontact manifold, a cocontact Hamiltonian vector field $X^{cc}_H$ is a strict contact vector field if and only if $\mathcal{R}^\eta(H)=0$. Therefore, the Hamiltonian function must independent of $z$-variable. So, we have the space of strict cocontact Hamiltonian vector fields defined as 
\begin{equation}\label{X-qcc-ham-evo}
\iota_{\bar{\xi}_{H}^{cc} }d \eta = d H-\langle dH, \mathcal{R}^\tau \rangle\tau, ~
			\iota_{\bar{\xi}_{H}^{cc} }\eta = -H, ~
			\iota_{\bar{\xi}_{H}^{cc} }\tau = 1.
\end{equation}  
Notice that one has the following identity
\begin{equation}
\bar{\xi}_{H}^{cc} (t,q,p,z)= \sharp^{cc}(dH)+ \left(1-\langle dH, \mathcal{R}^\tau \rangle \right) \mathcal{R}^\tau-H \mathcal{R}^\eta 
\end{equation}
 in terms of the musical mapping $\flat^{cc}$ in \eqref{flat-map}. 
In the Darboux coordinates $(t,q^{i},p_{i},z)$, if the Hamiltonian function $H=H(t,q,p)$ is independent of the fiber variable $z$, the strict cocontact Hamiltonian vector field is 
\begin{equation}\label{con-dyn-xibar}
\bar{\xi}_{H}^{cc}(t,q,p,z)=\frac{\partial}{\partial t}+ \frac{\partial H}{\partial p_{i}}\frac{\partial}{\partial q^{i}}  -  \frac{\partial H}{\partial q^{i}} 
\frac{\partial}{\partial p_{i}} + \left(p_{i}\frac{\partial H}{\partial p_{i}} - H\right)\frac{\partial}{\partial z}.
\end{equation}
Thus, we obtain the strict cocontact Hamilton's equations as  
\begin{equation}\label{conham-xibar}
\dot{t}=1 ,\qquad 
\dot{q}^i = \frac{\partial H}{\partial p_{i}}, 
\qquad
\dot{p}_i  = -\frac{\partial H}{\partial q^{i}}, 
\qquad \dot{z}= p_{i}\frac{\partial H}{\partial p_{i}} - H.
\end{equation}
 Let us investigate the conservation of the volume form. By a straightforward calculation we obtain the identities
\begin{equation}\label{L-xibar-eta-cc}
\begin{split}
    \mathcal{L}_{\bar{\xi}_{H}^{cc}}\tau & = 0,
\\
    \mathcal{L}_{\bar{\xi}_{H}^{cc}}\eta & = - \langle dH,\mathcal{R}^\tau \rangle \tau,
\\
    \mathcal{L}_{\bar{\xi}_{H}^{cc}}d\eta & = - d\langle dH,\mathcal{R}^\tau \rangle \wedge \tau.
\end{split}
\end{equation}
We observe that a strict evolution vector field preserves the one-form $\tau$ whereas it does not preserve the one-form $\eta$.
Using \eqref{L-xibar-eta-cc}, we deduce that the Lie derivative of the volume form is 
\begin{equation}\label{L-xibar-eta-cc-vol}
\begin{split}
    \mathcal{L}_{\bar{\xi}_{H}^{cc}} ( \tau \wedge d\eta ^n \wedge \eta ) = 0.
    \end{split}
\end{equation}
From the equation \eqref{L-xibar-eta-cc-vol} we have that $\bar{\xi}_{H}^{cc}$ is divergence-free. 
On the other hand, the Lie derivative of the Hamiltonian function $H$ is
\begin{equation}
{\mathcal{L}}_{\bar{\xi}_{H}^{cc}} \, H = \mathcal{R}^\tau(H).
\end{equation}
So, the strict evolution vector field does not preserve the Hamiltonian function $H$.

 \subsection{Cocontact Gradient  Dynamics}
\label{sec-coco-grad}
 
Consider a cocontact manifold $(\bar{N},\tau,\eta)$, and 
its associated algebra of cocontact vector fields presented in \eqref{X-cc}. In the case of cocontact Hamiltonian vector fields \eqref{murat}, one integrates the second condition $\mathcal{L}_{X} \tau=0$ as $\iota_{X} \tau=0$ (see the third condition in \eqref{X-cc-ham}). Instead, in this subsection we shall examine the three vector fields defined in Subsection \ref{sec-coco-ham} but replacing the condition $\iota_{X} \tau=0$ by $\iota_{X} \tau= \langle dH,\mathcal{R}^\tau  \rangle $.

\noindent 
\textbf{Cocontact Gradient Vector Fields.}  We define a cocontact evolution vector field $\bar{X}_H^{cc}$ as   
	\begin{equation} 
	\iota_{\mathrm{grad}X_H^{cc}}d \eta = d H-\langle dH, \mathcal{R}^\eta \rangle \eta-\langle dH, \mathcal{R}^\tau \rangle\tau, \quad 			\iota_{\mathrm{grad}X_H^{cc}}\eta = -H, \quad  			\iota_{\mathrm{grad}X_H^{cc}} \tau= \langle dH,\mathcal{R}^\tau  \rangle . 
	\end{equation}
In terms of the musical mappings $\sharp^{cc}$ and $\sharp_{\Lambda^{cc}}$ in \eqref{sharp-L-cc} the cocontact gradient vector field becomes 
\begin{equation}
\begin{split}
\mathrm{grad}X_H^{cc}   =\sharp^{cc}(dH) 
- \left( H + \langle dH, \mathcal{R}^\eta \rangle  \right) \mathcal{R}^\eta  .
\end{split}
\end{equation}
In the Darboux coordinates $(t,q^{i},p_{i},z)$, for a given Hamiltonian function $H=H(t,q,p,z)$, the cocontact gradient  vector field turns out to be
\begin{equation}\label{con-dyn-grad}
\mathrm{grad}X_H^{cc} (t,q,p,z) = \frac{\partial H}{\partial t}\frac{\partial}{\partial t} + \frac{\partial H}{\partial p_{i}}\frac{\partial}{\partial q^{i}}  - \left(\frac{\partial H}{\partial q^{i}} + \frac{\partial H}{\partial z} p_{i} \right)
\frac{\partial}{\partial p_{i}} + \left(p_{i}\frac{\partial H}{\partial p_{i}} - H\right)\frac{\partial}{\partial z}.
\end{equation}
Thus, we obtain the cocontact gradient equations for $H$ as
\begin{equation}\label{conham-grad}
\dot{t}= \frac{\partial H}{\partial t} 
,\qquad 
\dot{q}^i = \frac{\partial H}{\partial p_{i}}, 
\qquad
\dot{p}_i  = -\frac{\partial H}{\partial q^{i}}- 
p_{i}\frac{\partial H}{\partial z}, 
\qquad \dot{z}= p_{i}\frac{\partial H}{\partial p_{i}} - H.
\end{equation}

We obtain the following identities 
\begin{equation}\label{L-grad-eta-cc}
\begin{split}
    \mathcal{L}_{\mathrm{grad}X_H^{cc}}\tau & = d\langle dH,\mathcal{R}^\tau \rangle,
\\
    \mathcal{L}_{\mathrm{grad}X_H^{cc}}\eta & = - \langle dH,\mathcal{R}^\eta \rangle \eta - \langle dH,\mathcal{R}^\tau \rangle \tau,
\\
    \mathcal{L}_{\mathrm{grad}X_H^{cc}}d\eta & = - d\langle dH,\mathcal{R}^\eta \rangle \wedge \eta - \langle dH,\mathcal{R}^\eta \rangle d\eta - d\langle dH,\mathcal{R}^\tau \rangle \wedge \tau.
\end{split}
\end{equation}
So, the cocontact gradient vector field preserves neither the one-form $\tau$ nor the one-form $\eta$.  From the identities in \eqref{L-X-eta-cc}, we can compute the Lie derivative of the top form as
\begin{equation}\label{L-grad-eta-cc-vol}
\begin{split}
    \mathcal{L}_{\mathrm{grad}X_H^{cc}} ( \tau \wedge d\eta ^n \wedge \eta ) = \left((\mathcal{R}^\tau)^2(H) - (n+1)\mathcal{R}^\eta(H)\right) \tau \wedge d\eta ^n \wedge \eta.
    \end{split}
\end{equation}
That is, we find out that a cocontact gradient vector field does not preserve the volume form. The divergence of a cocontact Hamiltonian vector field is given by
\begin{equation} \label{div-grad-H}
\mathrm{div}(\mathrm{grad}X_H^{cc})= (\mathcal{R}^\tau)^2(H) - (n+1)\mathcal{R}^\eta(H).
\end{equation}
Moreover, the cocontact gradient vector field does not preserve the Hamiltonian function $H$ since the Lie derivative of the Hamiltonian function is
\begin{equation}
{\mathcal{L}}_{\mathrm{grad}X_H^{cc}} \, H = \left(\mathcal{R}^\tau(H)\right)^2 - \mathcal{R}^\eta(H)H.
\end{equation}

\noindent \textbf{Energy Gradient Vector Fields.}  If we consider the almost Poisson structure $(\bar{N},\Lambda^{cc})$, we determine the energy preserving cocontact dynamics in \eqref{E-H-cc-alt}. We replace the condition $\iota_{X} \tau=0$  by $\iota_{X} \tau=\langle dH,\mathcal{R}^\tau  \rangle$
and arrive at the definition of energy evolution vector field $\mathrm{grad} \mathcal{E}_{H}^{cc}$ as 
\begin{equation}\label{gradE-H-cc-alt--} 
\iota_{\mathrm{grad} \mathcal{E}_{H}^{cc}}d\eta=dH-\langle dH, \mathcal{R}^\eta \rangle \eta-\langle dH, \mathcal{R}^\tau \rangle\tau,\qquad \iota_{\mathrm{grad} \mathcal{E}_{H}^{cc}}\eta=0, \qquad 
	\iota_{\mathrm{grad} \mathcal{E}_{H}^{cc}}\tau = \langle dH,\mathcal{R}^\tau  \rangle.
   \end{equation}  
Using the musical mapping $\sharp_{\Lambda^{cc}}$, we write the energy evolution vector field as 
\begin{equation}
\mathrm{grad} \mathcal{E}_{H}^{cc}=  \sharp^{cc}(dH) - \langle dH, \mathcal{R}^\eta \rangle \mathcal{R}^\eta.
\end{equation}
In the Darboux coordinates, the energy gradient vector field is expressed as follows
\begin{equation}\label{evo-dyn-cc-}	\mathrm{grad} \mathcal{E}_{H}^{cc}(t,q,p,z)=\frac{\partial H}{\partial t}\frac{\partial}{\partial t} +\frac{\partial H}{\partial p_i}\frac{\partial}{\partial q^i}  - \left (\frac{\partial H}{\partial q^i} + \frac{\partial H}{\partial z} p_i \right)
	\frac{\partial}{\partial p_i} + p_i\frac{\partial H}{\partial p_i} \frac{\partial}{\partial z},
\end{equation}
so that, the dynamics is 
\begin{equation}\label{evo-eq-cc-grad}
\dot{t}= \frac{\partial H}{\partial t} , \qquad 	\frac{d q^i}{dt}= \frac{\partial H}{\partial p_i}, \qquad \frac{d p_i}{dt}  = -\frac{\partial H}{\partial q^i}- 
	p_i\frac{\partial H}{\partial z}, \quad \frac{d z}{dt} = p_i\frac{\partial H}{\partial p_i}.
\end{equation}
The energy gradient vector field does not preserve the forms $\tau,\eta$ and $d\eta$ since we have
\begin{equation}\label{L-gradE-eta-cc}
\begin{split}
    \mathcal{L}_{\mathrm{grad} \mathcal{E}_{H}^{cc}}\tau & = d\langle dH,\mathcal{R}^\tau \rangle,
\\
    \mathcal{L}_{\mathrm{grad} \mathcal{E}_{H}^{cc}}\eta & = dH - \langle dH,\mathcal{R}^\eta \rangle \eta - \langle dH,\mathcal{R}^\tau \rangle \tau,
\\
    \mathcal{L}_{\mathrm{grad} \mathcal{E}_{H}^{cc}}d\eta & = - d\langle dH,\mathcal{R}^\eta \rangle \wedge \eta - \langle dH,\mathcal{R}^\eta \rangle d\eta - d\langle dH,\mathcal{R}^\tau \rangle \wedge \tau.
\end{split}
\end{equation}
Moreover, these identities lead us to compute the Lie derivative of the volume form
\begin{equation}\label{L-gradE-eta-cc-vol}
\begin{split}
    \mathcal{L}_{\mathrm{grad} \mathcal{E}_{H}^{cc}} ( \tau \wedge d\eta ^n \wedge \eta ) = \left((\mathcal{R}^\tau)^2(H) - n\mathcal{R}^\eta(H)\right) \tau \wedge d\eta ^n \wedge \eta.
    \end{split}
\end{equation}
We deduce that the energy gradient vector field does not preserve the volume form. We obtain the divergence of energy gradient vector field as
\begin{equation} \label{div-gradE-H}
\mathrm{div}(\mathrm{grad} \mathcal{E}_{H}^{cc})= (\mathcal{R}^\tau)^2(H) - n\mathcal{R}^\eta(H).
\end{equation}

For the conservation of the Hamiltonian function $H$, we calculate the Lie derivative
\begin{equation}\label{res}
{\mathcal{L}}_{\mathrm{grad} \mathcal{E}_{H}^{cc}} \, H = \left(\mathcal{R}^\tau(H)\right)^2
\end{equation}
of the Hamiltonian function $H$ in the direction of the energy gradient vector field $\mathrm{grad} \mathcal{E}_{H}^{cc}$. The result in \eqref{res} reads that the flows of the energy gradient vector field do not preserve the Hamiltonian function.

 \noindent 
\textbf{Strict Gradient Vector Fields.} 
For a given Hamiltonian function $H$ on the cocontact manifold, a cocontact Hamiltonian vector field $X^{cc}_H$ is a strict contact vector field if and only if $\mathcal{R}^\eta(H)=0$. This locally gives that Hamiltonian function must independent of $z$-variable. So we have the space of strict gradient vector field $\mathrm{grad} \xi_{H}^{cc}$ defined as follows 
\begin{equation}\label{X-qcc-ham-grad}
\iota_{\mathrm{grad} \xi_{H}^{cc} }d \eta = d H-\langle dH, \mathcal{R}^\tau \rangle\tau, \qquad 
			\iota_{\mathrm{grad} \xi_{H}^{cc} }\eta = -H,  \qquad 
			\iota_{\mathrm{grad} \xi_{H}^{cc} }\tau = \langle dH,\mathcal{R}^\tau  \rangle.
\end{equation}  
Note that one has the following identity  \eqref{flat-map} 
\begin{equation}
\mathrm{grad} \xi_{H}^{cc}  = \sharp^{cc}(dH) -H \mathcal{R}^\eta .
\end{equation}
In the Darboux coordinates $(t,q^{i},p_{i},z)$, for a Hamiltonian function $H=H(t,q,p,z)$ independent of the fiber variable $z$, the strict gradient vector field is 
\begin{equation}\label{con-dyn-gradxi}
\mathrm{grad} \xi_{H}^{cc} (t,q,p,z)=\frac{\partial H}{\partial t}\frac{\partial}{\partial t} +\frac{\partial H}{\partial p_{i}}\frac{\partial}{\partial q^{i}}  -  \frac{\partial H}{\partial q^{i}} 
\frac{\partial}{\partial p_{i}} + \left(p_{i}\frac{\partial H}{\partial p_{i}} - H\right)\frac{\partial}{\partial z}.
\end{equation}
Thus, we obtain strict gradient equations as  
\begin{equation}\label{conham-gradxi}
\dot{t}= \frac{\partial H}{\partial t}
, \qquad 
\dot{q}^i = \frac{\partial H}{\partial p_{i}}, 
\qquad
\dot{p}_i  = -\frac{\partial H}{\partial q^{i}}, 
\qquad \dot{z}= p_{i}\frac{\partial H}{\partial p_{i}} - H.
\end{equation}

 The following identities show that the one forms $\tau$ and $\eta$ are not preserved under the strict gradient vector field:
\begin{equation}\label{L-gradxi-eta-cc}
\begin{split}
    \mathcal{L}_{\mathrm{grad} \xi_{H}^{cc}}\tau & = d\langle dH,\mathcal{R}^\tau \rangle,
\\
    \mathcal{L}_{\mathrm{grad} \xi_{H}^{cc}}\eta & = - \langle dH,\mathcal{R}^\tau \rangle \tau,
\\
    \mathcal{L}_{\mathrm{grad} \xi_{H}^{cc}}d\eta & = - d\langle dH,\mathcal{R}^\tau \rangle \wedge \tau,
\end{split}
\end{equation}
Further, we have that
\begin{equation}\label{L-gradxi-eta-cc-vol}
\begin{split}
    \mathcal{L}_{\mathrm{grad} \xi_{H}^{cc}} ( \tau \wedge d\eta ^n \wedge \eta ) = (\mathcal{R}^\tau)^2(H) \tau \wedge d\eta ^n \wedge \eta.
    \end{split}
\end{equation}
Thus, we deduce that the strict gradient vector field does not preserve the volume form. We write the divergence of the strict gradient vector field as
\begin{equation} \label{div-gradxi-H}
\mathrm{div}(\mathrm{grad} \xi_{H}^{cc})= (\mathcal{R}^\tau)^2(H).
\end{equation}

\subsection{Kinetic Lift of Cocontact Hamiltonian Dynamics: Momentum Formulation}\label{sec-kin-coco-mom}

 Consider a cocontact manifold $(\bar{N},\tau,\eta)$.We now determine the dual space $\mathfrak{X}_{\mathrm{cc}}^* (\bar{N})$ of the space of contact Hamiltonian vector fields $\mathfrak{X}_{\mathrm{cc-ham}}(\bar{N})$ given in \eqref{X-con}. We first note that   $\mathfrak{X}_{\mathrm{cc-ham}}^*(\bar{N})$ is a subspace of the space $\Gamma^1(\bar{N})\otimes \mathrm{Den} (\bar{N})$ of one-form densities. To be more precise, we compute the $L_2$-pairing (simply multiply-and-integrate) of an arbitrary contact vector field $X_{H}^{cc} $ with a one-form density $\bar{\Upsilon} \otimes \overline{\nu}$. Making use of the identities of the Cartan calculus, we obtain
\begin{equation}\label{cc-dual-comp}
\begin{split}
 \langle \bar{\Upsilon} \otimes \overline{\nu},& X_{H}^{cc}   \rangle _{L_2} = \int_{\bar{N}} \langle \bar{\Upsilon} , X_{H}^{cc} 
 \rangle   \overline{\nu} =\int_{\bar{N}} \langle \bar{\Upsilon} , \sharp_{\Lambda^{cc}}(dH) -H\mathcal{R}^\eta    \rangle\overline{\nu} 
\\
&= - \int_{\bar{N}} \langle {\sharp_{\Lambda^{cc}}(\bar{\Upsilon} )}, dH  \rangle  \overline{\nu}  -  \int_{\bar{N}} H  \langle \bar{\Upsilon} ,\mathcal{R}^\eta   \rangle  \overline{\nu}  \\
&= - \int_{\bar{N}} \langle \sharp^{cc}\bar{\Upsilon} - \langle \bar{\Upsilon} , \mathcal{R}^\eta \rangle \mathcal{R}^\eta - \langle \bar{\Upsilon} , \mathcal{R}^\tau \rangle \mathcal{R}^\tau, dH  \rangle   \overline{\nu}  -  \int_{\bar{N}} H  \langle \bar{\Upsilon} ,\mathcal{R}^\eta   \rangle  \overline{\nu} 
\\
&= - \int_{\bar{N}} \left(\iota_{{\sharp^{cc}\bar{\Upsilon} }}  dH \right)  \overline{\nu}  + \int_{\bar{N}}  \langle \bar{\Upsilon} , \mathcal{R}^\eta \rangle  \langle \mathcal{R}^\eta,dH \rangle  \overline{\nu}  
+
\int_{\bar{N}}  \langle \bar{\Upsilon} , \mathcal{R}^\tau \rangle  \langle \mathcal{R}^\tau,dH \rangle  \overline{\nu} \\&
  \qquad   \qquad 
-  \int_{\bar{N}} H  \langle \bar{\Upsilon} ,\mathcal{R}^\eta   \rangle  \overline{\nu} 
\\&=  -
 \int_{\bar{N}}  dH \wedge \iota_{{\sharp^{cc}\bar{\Upsilon} }} \overline{\nu} + 
 \int_{\bar{N}}  \iota_{\mathcal{R}^\eta}dH  \langle \bar{\Upsilon} , \mathcal{R}^\eta \rangle  
 \overline{\nu} +  
 \int_{\bar{N}}  \iota_{\mathcal{R}^\tau }dH  \langle \bar{\Upsilon} , \mathcal{R}^\tau \rangle  
 \overline{\nu} -  \int_{\bar{N}} H  \langle \bar{\Upsilon} ,\mathcal{R}^\eta   \rangle  \overline{\nu} 
  \\&=  \int_{\bar{N}}  H d \iota_{{\sharp^{cc}\bar{\Upsilon} }} \overline{\nu}
  + \int_{\bar{N}} \langle \bar{\Upsilon} , \mathcal{R}^\eta \rangle dH \wedge  \iota_{\mathcal{R}^\eta} \overline{\nu} 
  +\int_{\bar{N}} \langle \bar{\Upsilon} , \mathcal{R}^\tau \rangle dH \wedge  \iota_{\mathcal{R}^\tau} \overline{\nu} 
  - \int_{\bar{N}}   H   \langle \bar{\Upsilon} , \mathcal{R}^\eta \rangle  \overline{\nu} 
   \\&= \int_{\bar{N}}  H d \iota_{{\sharp^{cc}\bar{\Upsilon} }} \overline{\nu}
  - \int_{\bar{N}} H d\langle \bar{\Upsilon} , \mathcal{R}^\eta \rangle \wedge  \iota_{\mathcal{R}^\eta} \overline{\nu} 
  -\int_{\bar{N}} H  \langle \bar{\Upsilon} , \mathcal{R}^\eta \rangle    d\iota_{\mathcal{R}^\eta} \overline{\nu}
  \\&
  \qquad   \qquad 
  - \int_{\bar{N}} H d\langle \bar{\Upsilon} , \mathcal{R}^\tau\rangle \wedge  \iota_{\mathcal{R}^\tau} \overline{\nu} 
  -\int_{\bar{N}} H  \langle \bar{\Upsilon} , \mathcal{R}^\tau \rangle    d\iota_{\mathcal{R}^\tau} \overline{\nu}
- \int_{\bar{N}}   H   \langle \bar{\Upsilon} , \mathcal{R}^\eta \rangle  \overline{\nu}    
  \\ &= 
  \int_{\bar{N}}  H \left( \mathrm{div}({\sharp^{cc}\bar{\Upsilon} })  
  -   \iota_{\mathcal{R}^\eta}d\langle \bar{\Upsilon} , \mathcal{R}^\eta \rangle  - \langle \bar{\Upsilon} , \mathcal{R}^\eta \rangle    \mathrm{div}(\mathcal{R}^\eta)   \right)  \overline{\nu}   
   \\&
  \qquad   \qquad 
  - \int_{\bar{N}}  H \left( 
  \iota_{\mathcal{R}^\tau}d\langle \bar{\Upsilon} , \mathcal{R}^\tau \rangle  + \langle \bar{\Upsilon} , \mathcal{R}^\tau \rangle    \mathrm{div}(\mathcal{R}^\tau) + \langle \bar{\Upsilon} , \mathcal{R}^\eta \rangle   \right)  \overline{\nu} 
 \\ &= 
  \int_{\bar{N}}  H \left( \mathrm{div}({\sharp^{cc}\bar{\Upsilon} }) 
  -   \mathcal{L}_{\mathcal{R}^\eta} \langle \bar{\Upsilon} , \mathcal{R}^\eta \rangle  -    \mathcal{L}_{\mathcal{R}^\tau} \langle \bar{\Upsilon} , \mathcal{R}^\tau \rangle  - \langle \bar{\Upsilon} , \mathcal{R}^\eta \rangle \right)  \overline{\nu} ,
\end{split}
\end{equation}
where $\mathrm{div}$ stands for the divergence with respect to the cocontact volume $\overline{\nu} $ in \eqref{cocontact-volume}. Notice that the divergence of the Reeb vector fields are zero, and we employed these observations in the last line of the calculation. If the volume form is fixed, the non-degeneracy of the pairing gives the dual space as
\begin{equation} \label{cc-Ham-alg}
\mathfrak{X}_{\mathrm{cc-ham}}^* (\bar{N}) = \big \{\bar{\Upsilon}  \in 
\Gamma^1 (\bar{N}) : \mathrm{div}({\sharp^{cc}\bar{\Upsilon} }) 
  -   \mathcal{L}_{\mathcal{R}^\eta} \langle \bar{\Upsilon} , \mathcal{R}^\eta \rangle  -    \mathcal{L}_{\mathcal{R}^\tau} \langle \bar{\Upsilon} , \mathcal{R}^\tau \rangle  - \langle \bar{\Upsilon} , \mathcal{R}^\eta \rangle \neq 0 \big
\}  \cup \{0\}.
\end{equation}
We next recall the Lie algebra isomorphism $\bar{\Psi}:H\mapsto X_{H}^{cc} $ of \eqref{bar-psi}. Identifying the dual $C^\infty(\bar{N})$ with the space of densities on the cocontact manifold, and fixing the contact volume form, \eqref{cc-dual-comp} determines the dual mapping of \eqref{X-cc-ham} as
\begin{equation}\label{Delta-*}
\begin{split}
\bar{\Psi}^*:  \mathfrak{X}_{\mathrm{cc-ham}}^*(\bar{N}) & \longrightarrow (C^\infty)^*(\bar{N}),\\
&\bar{\Upsilon} \mapsto f= \mathrm{div}({\sharp^{cc}\bar{\Upsilon} }) 
  -   \mathcal{L}_{\mathcal{R}^\eta} \langle \bar{\Upsilon} , \mathcal{R}^\eta \rangle  -    \mathcal{L}_{\mathcal{R}^\tau} \langle \bar{\Upsilon} , \mathcal{R}^\tau \rangle  - \langle \bar{\Upsilon} , \mathcal{R}^\eta \rangle.
  \end{split}
\end{equation}
In terms of the Darboux coordinates $(t,q^i,p_i,z)$, we can compute the density function $f$ from \eqref{Delta-*} as follows. Consider a one-form section 
\begin{equation}
\bar{\Upsilon}  =\bar{\Upsilon} _i dq^i + \bar{\Upsilon} ^i dp_i + \bar{\Upsilon} _t dt+\bar{\Upsilon} _z dz.
\end{equation}
Let us present the three terms on the right-hand side of \eqref{Delta-*} one by one. 
A direct calculation reads the contact divergence of $\bar{\Upsilon} $ as 
\begin{equation}
 \mathrm{div}\left({\sharp^{cc}\bar{\Upsilon} }\right)  = \frac{\partial \bar{\Upsilon}^i}{\partial q^i} - \frac{\partial \bar{\Upsilon}_i}{\partial p_i}  -p_i\left(\frac{\partial \bar{\Upsilon}_z}{\partial p_i}-\frac{\partial \bar{\Upsilon}^i}{\partial z} \right) - n\bar{\Upsilon}_z + \frac{\partial \bar{\Upsilon}_z}{\partial z} + \frac{\partial \bar{\Upsilon}_t}{\partial t}.
\end{equation}
Then the second and third Lie derivative term is locally computed to be 
\begin{equation} 
  \mathcal{L}_{\mathcal{R}^\eta} \langle \bar{\Upsilon} , \mathcal{R}^\eta \rangle   = \mathcal{R}^\eta(\bar{\Upsilon} _z) = \frac{\partial \bar{\Upsilon} _z}{\partial z},\qquad     \mathcal{L}_{\mathcal{R}^\tau} \langle \bar{\Upsilon} , \mathcal{R}^\tau \rangle   = \mathcal{R}^\tau(\bar{\Upsilon} _t) = \frac{\partial \bar{\Upsilon} _t}{\partial t} .
\end{equation}
The fourth term is simply $ \langle \bar{\Upsilon} , \mathcal{R}^\eta \rangle = \bar{\Upsilon} _z$. Adding all of these terms we arrive at the definition of the time-dependent dissipative density function as 
\begin{equation}
   f(t,q,p,z)=  \frac{\partial \Upsilon^i}{\partial q^i} - \frac{\partial \Upsilon_i}{\partial p_i}  -p_i(\frac{\partial \Upsilon_z}{\partial p_i}-\frac{\partial \Upsilon^i}{\partial z} ) -(n+1)\Upsilon_z.
    \end{equation}

A cocontact Hamiltonian vector field $X^{cc}_H$ has non-zero divergence as one easily deduces from \eqref{div-X-H}; moreover, referring to \eqref{coad-gen}, the coadjoint action of $X^{cc}_H $ on a dual element $\bar{\Upsilon}$ is just
\begin{equation}
ad_{X_{H}^{cc} }^\ast \bar{\Upsilon} = \mathcal{L}_{X_{H}^{cc} }  \bar{\Upsilon} -  (n+1)  \mathcal{R}^\eta(H) \bar{\Upsilon} .
\end{equation} 
So that the Lie-Poisson dynamics is what we call cocontact momentum-Vlasov dynamics
\begin{equation}\label{LP-cont-mom}
\dot{ \bar{\Upsilon} }= - \mathcal{L}_{X_{H}^{cc} }  \bar{\Upsilon} +  (n+1)  \mathcal{R}^\eta(H) \bar{\Upsilon} .
\end{equation}

\noindent \textbf{Contact (Autonomous) Kinetic Dynamics in Momentum Formulation.} Let us study the particular situation of no explicit time dependence. This geometrically refers to a contact manifold
 $(N,\eta)$. In this case, the associated Lie algebra is the space $\mathfrak{X}_{\mathrm{c-ham}}(N)$ of contact Hamiltonian vector fields given in \eqref{X-con}. So we have no (time) Reeb vector field $\mathcal{R}^\tau$ since there is no time dependence. We shall show the Lie-Poisson dynamics on the dual of the space of contact Hamiltonian vector fields $\mathfrak{X}_{\mathrm{c-ham}}^*(N)$. For more explicit discussion and calculations, one may consult a recent exposition  \cite{esen2023conformal} on contact kinetic theory. Once the contact volume is fixed, an examination of the definition \eqref{cc-Ham-alg} for a time-independent case implies that the dual space is 
\begin{equation} \label{c-Ham-alg}
\mathfrak{X}_{\mathrm{c-ham}}^* (N) = \big \{\Upsilon  \in 
\Gamma^1 (N) :\mathrm{div}({\sharp^{c}\Upsilon}) 
  -   \mathcal{L}_{\mathcal{R}^\eta} \langle \Upsilon, \mathcal{R}^\eta \rangle   - \langle \Upsilon, \mathcal{R}^\eta \rangle \neq 0 \big
\}  \cup \{0\}.
\end{equation}
We next recall the Lie algebra isomorphism $H\mapsto X_{H}^c $ of \eqref{iso}. Identifying the dual $C^\infty(N)$ with the space of densities on the contact manifold, and fixing the contact volume form, the definition in \eqref{Delta-*} reduces to the following way to determine the dual of \eqref{iso} as
\begin{equation}\label{Psi-*}
\Psi^*:  \mathfrak{X}_{\mathrm{c-ham}}^*(N) \longrightarrow (C^\infty)^*(N),\qquad \Upsilon \mapsto f=\mathrm{div}({\sharp^{c}\Upsilon}) 
  -   \mathcal{L}_{\mathcal{R}^\eta} \langle \Upsilon, \mathcal{R}^\eta \rangle   - \langle \Upsilon, \mathcal{R}^\eta \rangle .
\end{equation}
In terms of the Darboux coordinates $(q^i,p_i,z)$, we can compute the density function $f$ from \eqref{Psi-*} as follows. Consider a one-form section 
\begin{equation}
\Upsilon  =\Upsilon_i dq^i + \Upsilon^i dp_i + \Upsilon_z dz \mapsto 
    f(q,p,z)=  \frac{\partial \Upsilon^i}{\partial q^i} - \frac{\partial \Upsilon_i}{\partial p_i}  -p_i(\frac{\partial \Upsilon_z}{\partial p_i}-\frac{\partial \Upsilon^i}{\partial z} ) -(n+1)\Upsilon_z.
    \end{equation}
A contact Hamiltonian vector field $X^{c}_H$ has non-zero divergence since \eqref{div-X-H}, and using \eqref{coad-gen}, we calculate the coadjoint action of $X^{c}_H $ on a dual element $\Upsilon$ 
\begin{equation}
ad_{X_{H}^c }^\ast \Upsilon= \mathcal{L}_{X_{H}^c } \Upsilon -  (n+1)  \mathcal{R}^\eta(H)\Upsilon.
\end{equation} 
So that the Lie-Poisson dynamics is what we call contact momentum-Vlasov dynamics
\begin{equation}\label{LP-cont-mom-}
\dot{\Upsilon}= - \mathcal{L}_{X_{H}^c } \Upsilon +  (n+1)  \mathcal{R}^\eta(H)\Upsilon.
\end{equation}

\subsection{Kinetic Lift of Cocontact Hamiltonian Dynamics: Density Formulation}\label{sec-kin-coco-den}

Let $(\bar{N},\tau,\eta)$
 be a  cocontact manifold.  
 We first start with the following observation 
 \begin{equation}\label{flip-flow-}
\int \{F,H\}^{cc}K\overline{\nu}  = \int F \{H,K\}^{cc} \overline{\nu}    +  (n+3) \int  FK \mathcal{R}^\eta(H)  \overline{\nu}  
\end{equation}
for all smooth functions $F$, $H$, and $K$ defined on the cocontact manifold $\bar{N}$. Indeed, we have
\begin{equation}\label{cont-bra-comt}
\begin{split}
\int \{F,H\}^{cc}K\overline{\nu}  &= \int  \left( X_H^{cc}(F)+F \mathcal{R}^\eta(H)  \right)  K\overline{\nu}  
\\ &=\int  K\left(\iota_{X_H^{cc}} dF \right)      \overline{\nu}  + \int KF \mathcal{R}^\eta(H)  \overline{\nu}  
\\ &=\int  KdF    \wedge  \iota_{X_H^{cc}} \overline{\nu}  + \int KF \mathcal{R}^\eta(H)  \overline{\nu}  
\\ &= - \int  FdK\wedge  \iota_{X_H^{cc}} \overline{\nu} 
- \int  FK d\iota_{X_H^{cc}} \overline{\nu}  
 + \int KF \mathcal{R}^\eta(H)  \overline{\nu}  
 \\ &= 
 - \int  F(\iota_{X_H^{cc}} dK) \overline{\nu} 
- \int  FK \mathrm{div}(X_H^{cc}) \overline{\nu}  
 + \int KF \mathcal{R}^\eta(H)  \overline{\nu}  
 \\ &= 
 - \int  FX_H^{cc}(K) \overline{\nu} 
+ \int  FK (n+1)  \mathcal{R}^\eta(H)  \overline{\nu}  
 + \int KF \mathcal{R}^\eta(H)  \overline{\nu}  
 \\
  &= 
 - \int  FX_H^{cc}(K) \overline{\nu} 
+ \int  FK (n+1)  \mathcal{R}^\eta(H)  \overline{\nu}  
 + \int KF \mathcal{R}^\eta(H)  \overline{\nu}  
  \\
  &= - \int F \left(\{K,H\}^{cc}-K\mathcal{R}^\eta(H) \right)   + (n+2) \int  FK  \mathcal{R}^\eta(H)  \overline{\nu}  
    \\
  &=  \int F \{H,K\}^{cc} \overline{\nu}    +  (n+3) \int  FK \mathcal{R}^\eta(H)  \overline{\nu} .
\end{split}
\end{equation}
From \eqref{flip-flow-}, we obtain that 
 \begin{equation}\label{flip-flow--}
\begin{split}
\int_N X_{H}^{cc} (F)K\overline{\nu}  &=  \int_N \{F,H \}^{cc} K-  F\mathcal{R}^\eta(H) K \overline{\nu}  \\
&= \int_N F \{H,K\}^{cc} \overline{\nu}    +  (n+2) \int_N  FK \mathcal{R}^\eta(H) \overline{\nu}  .
\end{split}
\end{equation}
This calculation leads us to obtain the dual mapping of the Lie derivative action 
 \begin{equation}\label{Psi-mapping}
\bar{\Xi}_{X_{H}^{cc} }:C^\infty(\bar{N}) \longrightarrow C^\infty(\bar{N}),\qquad (X_{H}^{cc} ,K)\mapsto X_{H}^{cc} (K) 
 \end{equation}
 of a cocontact Hamiltonian vector field $X_{H}^{cc} $. 
Acordingly, the dual of \eqref{Psi-mapping} is given by 
 \begin{equation}
\bar{\Xi}^*_{X_{H}^{cc} }: (C^\infty)^*(\bar{N}) \longrightarrow (C^\infty)^*(\bar{N}),\qquad f\mapsto \Psi^*_{X_{H}^{cc} }(f) = \{H,f\}^{cc} - (n+2) f.
 \end{equation}
We then define the dynamics generated by the dual action as
    \begin{equation}\label{LP-con-den-VF}
\dot{f} = - \bar{\Xi}^*_{X_{H}^{cc} }(f) = - \{H,F\}^{cc} + (n+2) f \mathcal{R}^\eta(H). 
\end{equation} 
In terms of the Darboux coordinates $(t,q^i,p_i,z)$, the cocontact kinetic dynamics in density form reads
\begin{equation} 
    \frac{\partial f}{\partial t} +\frac{\partial H}{\partial p_i} \frac{\partial f}{\partial q^i} -\frac{\partial H}{\partial q^i} \frac{\partial f}{\partial p_i}
 =p_i \left(\frac{\partial f}{\partial p_i} \frac{\partial H}{\partial z}- \frac{\partial H}{\partial p_i} \frac{\partial f}{\partial z}\right) +(n+1)f \frac{\partial H}{\partial z} + \frac{\partial f}{\partial z} H. 
\end{equation} \label{LP-cont-Psi}

\noindent \textbf{Contact (Autonomous) Kinetic Dynamics in Density Formulation.} Let $(N,\eta)$ be a contact manifold; then, observe that   \begin{equation}\label{flip-flow-1}
\int \{F,H\}^{c}K\nu = \int F \{H,K\}^{c} \nu   +  (n+3) \int  FK \mathcal{R}^\eta(H)  \nu 
\end{equation}
for all smooth functions $F$, $H$, and $K$ defined on the contact manifold $N$. Therefore, 
\begin{equation}\label{flip-flow--1} 
\int_N X_{H}^c (F)K\nu  = \int_N F \{H,K\}^{c} \nu   +  (n+2) \int_N  FK \mathcal{R}^\eta(H) \nu . 
\end{equation}
This calculation leads us to obtain the dual mapping of the Lie derivative action 
 \begin{equation}\label{Psi-mapping-}
\Xi_{X_{H}^c }:C^\infty(N) \longrightarrow C^\infty(N),\qquad (X_{H}^c ,K)\mapsto X_{H}^c (K) 
 \end{equation}
 of a contact vector field $X_{H}^c $. 
Acordingly, the dual of \eqref{Psi-mapping-} is given by 
 \begin{equation}
\Xi^*_{X_{H}^c }: (C^\infty)^*(N) \longrightarrow (C^\infty)^*(N),\qquad f\mapsto \Psi^*_{X_{H}^c }(f) = \{H,f\}^{c} - (n+2) f.
 \end{equation}
We then define the dynamics generated by the dual action as
    \begin{equation}\label{LP-con-den-VF-}
\dot{f} = - \Xi^*_{X_{H}^c }(f) = - \{H,F\}^{c} + (n+2) f \mathcal{R}^\eta(H). 
\end{equation} 
In terms of the Darboux coordinates $(q^i,p_i,z)$, the kinetic dynamics turns out to be the contact kinetic dynamics in density form 
\begin{equation} 
    \frac{\partial f}{\partial t} +\frac{\partial H}{\partial p_i} \frac{\partial f}{\partial q^i} -\frac{\partial H}{\partial q^i} \frac{\partial f}{\partial p_i}
 =p_i \left(\frac{\partial f}{\partial p_i} \frac{\partial H}{\partial z}- \frac{\partial H}{\partial p_i} \frac{\partial f}{\partial z}\right)  \quad+(n+1)f \frac{\partial H}{\partial z} + \frac{\partial f}{\partial z} H.
\end{equation} \label{LP-cont-Psi-}

\section{Summary and Discussion}

In this work, we have examined Lie-Poisson formulation of time-dependent kinetic theory in many aspects. To be more precise about our novel vector field definitions on cocontact manifolds, we first examine dynamical vector fields defined on cosymplectic manifolds (see Table \ref{cosymplectic-table}), and on contact manifolds (see Table \ref{contact-table}). We have proposed cosymplectic geometric analysis of time dependent Vlasov equation in Section \ref{sec-kin-cos-mom} in terms of the momentum variables, and in Section  \ref{sec-kin-cos-den} in terms of the plasma density function. These formulations have been linked by a Poisson/momentum map. The dual character (contact-cosymplectic) of cocontact manifolds lead us to determine various different vector fields in Section \ref{sec-coco-ham}, Section \ref{sec-coco-evo}, and Section \ref{sec-coco-grad}. We list all of them in Table \ref{cocontact-table}. 
Then we have introduced the time-dependent contact generalization of the Vlasov equation. This has been done in Section \ref{sec-kin-coco-mom} for momentum variables and \ref{sec-kin-coco-den} in terms of the plasma density function. Once again, a Poisson/momentum map has been constructed in order to display the relationship between these two formulations. This investigation involves the determination of proper Lie algebras and dual spaces. To list the results of the paper we first 
determine proper Lie algebras and the dual spaces in Table \ref{kinetic-algebra-table}. Later we present the kinetic equations in terms of the momentum and density functions in symplectic, cosymplectic, contact and cocontact manifolds. See Table \ref{kinetic-table}.

Now we list here some of the problems we shall focus on in upcoming papers: 
\begin{itemize}
   \item The Legendre transformation takes Lagrangian dynamics to Hamiltonian dynamics. For non-degenerate theories, defining the transformation is pretty straightforward. For degenerate theories, the transformation is not immediate. Tulczyjew triple is a geometric model formulating the Legendre transformation  
even for degenerate cases \cite{tul77,TuUr99}. Tulczyjew triple has been constructed for many different physical theories. For contact theory, in \cite{esen2021contact,esen2021implicit} the Tulczyjew triple has been constructed for coorientable contact setting whereas it is given in \cite{grabowska2022contact} for a more general setting without a need for coorientability. We plan to construct a Tulczyjew triple for cocontact geometry.

\item Each local chart of a symplectic manifold is symplectic, but the inverse of this assertion does not hold in general. For example, locally conformally symplectic manifolds, where local symplectic characteristics cannot be glued up to a global symplectic character \cite{Vaisman85}. In our previous work \cite{EsLeSaZa21}, we introduced the Hamilton-Jacobi formalism on locally conformally symplectic manifolds. We further explored this local-global interplay in the context of field theory in \cite{EsLeSaZa-k-sympl,EsLeSaZa-Cauchy}. More recently, we proposed the Hamilton-Jacobi manifold \cite{AtEsLe23} in the context of locally conformally cosymplectic manifolds \cite{ChLeMa91}. Following these approach, we plan to investigate the local-global interplay for cocontact manifolds, aiming to establish a geometric structure that we may call \emph{locally conformally cocontact manifolds} as well as to examine Hamiltonian dynamics and the Hamilton-Jacobi formalism within this framework.

\end{itemize}

\begin{landscape}
\begin{center}  
\begin{tabular}{ ||c|c||c|c|c|c|c|| }
 \hline
 \multicolumn{7}{||c||}{\textbf{Cosymplectic Manifold} \qquad $(\bar{M},\tau,\Omega)$, \qquad $\sharp_{\Lambda^{cs}}=\sharp^{cs}-\langle \bullet, \mathcal{R}^\tau \rangle \mathcal{R}^\tau$}\\
 \hline \hline
V-Field & Notation & $\iota_X\Omega$ & $\iota_X\tau$ &Music & $\mathrm{div} X$ & $X(H)$  \\
 \hline  \hline Hamiltonian &$X^{cs}$& 
 $dH- \langle dH,\mathcal{R}^{\tau} \rangle\tau$  & $0$  &  $\sharp^{cs}dH-\langle dH, \mathcal{R}^\tau \rangle \mathcal{R}^\tau$ & $0$  &  $0$  \\  \hline 
  Gradient &$\operatorname{grad} H$& 
 $dH- \langle dH,\mathcal{R}^{\tau} \rangle\tau$ & $\langle dH,\mathcal{R}^{\tau} \rangle$ & $\sharp^{cs} dH $  &   $(\mathcal{R}^\tau)^2(H)$  &  $\left(\mathcal{R}^\tau(H)\right)^2$  \\ \hline  
 Evolution &$E_H$& 
 $dH- \langle dH,\mathcal{R}^{\tau} \rangle\tau$ & $ 1$  &  $\sharp^{cs}dH + 
 (1- \langle dH,\mathcal{R}^{\tau} \rangle)\mathcal{R}^{\tau}$ &  $0$  &  $\mathcal{R}^\tau(H)$ \\  
 \hline
\end{tabular} 
\captionof{table}{Dynamics on Cosymplectic Manifolds \label{cosymplectic-table}}
\end{center}

\begin{center}
\begin{tabular}{ ||c|c||c|c|c|c|c||  }
 \hline
 \multicolumn{7}{||c||}{\textbf{Contact  Manifold} \qquad $(N,\eta)$, \qquad $\sharp_{\Lambda^c} =\sharp^c-\langle \bullet, \mathcal{R}^\eta\rangle \mathcal{R}^\eta$} \\
 \hline \hline
V-Field &  Notation & $\iota_Xd\eta$ & $\iota_X\eta$ &Music& $\mathrm{div} X$ & $X(H)$  \\
 \hline  \hline Contact & $X_{H}^c$& 
 $dH-\mathcal{R}^\eta(H) \eta$  & $-H$ & $\sharp^{c}dH-(H+\mathcal{R}^\eta(H))\mathcal{R}^\eta  $ & $-(n+1)\mathcal{R}^\eta(H)$ &  $-\mathcal{R}^\eta(H)H$ \\ Hamiltonian &&&&&&\\  \hline Contact & $\mathcal{E}_{H}$& 
 $dH-\mathcal{R}^\eta(H)\eta$ & $0$ &  $\sharp^{c}dH-\mathcal{R}^\eta(H)\mathcal{R}^\eta $ & $-n\mathcal{R}^\eta(H)$ &  $0$\\ Energy (evolution) &&&&&& \\  \hline  Strict Contact& $\xi_{H}$& 
 $dH$ & $-H$ &  $\sharp^{c}(dH)- H\mathcal{R}^\eta$ & $0$ &  $0$ \\  Hamiltonian  &&&&&& \\  
 \hline
\end{tabular}
\captionof{table}{Dynamics on Contact Manifolds \label{contact-table}}
\end{center}

\begin{center}
\begin{tabular}{ ||c|c||c|c|c|c|c||  }
 \hline
 \multicolumn{7}{||c||}{\textbf{Cocontact  Manifold} \qquad $(\bar{N},\tau,\eta)$, \qquad $ \sharp_{\Lambda^{cc}}  =\sharp^{cc} - \langle \bullet, \mathcal{R}^\eta \rangle \mathcal{R}^\eta - \langle \bullet, \mathcal{R}^\tau \rangle \mathcal{R}^\tau. $} \\
 \hline \hline
V-Field & Notation & $\iota_X d\eta$ & $\iota_X \eta$ & $\iota_X \tau$ & $\mathrm{div} X$ & $X(H)$    \\
 \hline  \hline
 Cocontact  &$X_{H}^{cc}$& 
 $	  d H-\langle dH, \mathcal{R}^\eta \rangle \eta-\langle dH, \mathcal{R}^\tau \rangle\tau$ & $			  -H$ & $			  0$ & $-(n+1)\mathcal{R}^\eta(H)$  &  $-\mathcal{R}^\eta(H)H$  
 \\ Hamiltonian &&&&&&\\
 \hline Cocontact&$\mathcal{E}_{H}^{cc}$& 
 $ dH-\langle dH, \mathcal{R}^\eta \rangle \eta-\langle dH, \mathcal{R}^\tau \rangle\tau$ & $  0$ & $   0$ & $-n\mathcal{R}^\eta(H)$  &  $0$  
  \\  Energy &&&&&&
 \\
 \hline Strict Cocontact &$\xi_{H}^{cc} $&   
 $ d H-\langle dH, \mathcal{R}^\tau \rangle\tau$ & $ 
			-H
			$ & $ 0$ & $0$  &  $0$  
    \\ Hamiltonian &&&&&&\\
 \hline 
 Cocontact&$\bar{X}_H^{cc}$&
 $ d H-\langle dH, \mathcal{R}^\eta \rangle \eta-\langle dH, \mathcal{R}^\tau \rangle\tau$ & $ 
	  -H
			$ & $  1$ & $-(n+1)\mathcal{R}^\eta(H)$  &  $\mathcal{R}^\tau(H) $  
 \\ Evolution &&&&&&$- \mathcal{R}^\eta(H)H$\\
 \hline 
 Cocontact  & $\bar{\mathcal{E}}_{H}^{cc}$ & 
$ dH-\langle dH, \mathcal{R}^\eta \rangle \eta-\langle dH, \mathcal{R}^\tau \rangle\tau$ & $ 		 0$ & $  1$ & $-n\mathcal{R}^\eta(H)$  &  $\mathcal{R}^\tau(H)$  
 \\ Energy Evolution &&&&&&\\
  \hline 
Strict  Cocontact   & $\bar{\xi}_{H}^{cc} $ &
$ d H-\langle dH, \mathcal{R}^\tau \rangle\tau$
	 & $ 		 -H$ & $			 1$ & $0$  &  $\mathcal{R}^\tau(H)$  
 \\Evolution &&&&&&\\
 \hline 
 Cocontact  & $\mathrm{grad}X_H^{cc}$ & 
 $  d H-\langle dH, \mathcal{R}^\eta \rangle \eta-\langle dH, \mathcal{R}^\tau \rangle\tau $ & $ 		  -H	$ & $ \langle dH,\mathcal{R}^\tau  \rangle$ & $(\mathcal{R}^\tau)^2(H) $  &  $(\mathcal{R}^\tau(H))^2$  
 \\Gradient &&&&&$- (n+1)\mathcal{R}^\eta(H)$&$ - \mathcal{R}^\eta(H)H$\\
 \hline 
Energy    & $\mathrm{grad} \mathcal{E}_{H}^{cc}$ & 
$ dH-\langle dH, \mathcal{R}^\eta \rangle \eta-\langle dH, \mathcal{R}^\tau \rangle\tau 
	$ & $ 		 0$ & $  \langle dH,\mathcal{R}^\tau  \rangle$ & $(\mathcal{R}^\tau)^2(H) $  &  $(\mathcal{R}^\tau(H))^2$  
 \\Gradient &&&&&$- n\mathcal{R}^\eta(H)$&\\
  \hline 
 Strict  & $\mathrm{grad} \xi_{H}^{cc} $ & 
$ d H-\langle dH, \mathcal{R}^\tau \rangle\tau	$
	 & $    -H$ & $		 \langle dH,\mathcal{R}^\tau  \rangle$ & $(\mathcal{R}^\tau)^2(H)$  &  $(\mathcal{R}^\tau(H))^2$  
 \\Gradient &&&&&&\\
 \hline
\end{tabular}
\captionof{table}{Dynamics on Cocontact Manifolds \label{cocontact-table}}
\end{center}

\begin{center}
\begin{tabular}{ ||c||c|c||  }
 \hline
 \multicolumn{3}{||c||}{\textbf{Hamiltonian Vector Fields and Dual Spaces}} \\
 \hline \hline
Manifold & Lie Algebra & Dual Space  \\
 \hline  \hline 
 Symplectic $M$& $\mathfrak{X}_{\mathrm{ham}}(M)$ in \eqref{alg-sym-ham}  &  $\mathfrak{X}_{\mathrm{ham}}^*(M)$ in \eqref{dual-symp-ham}   \\  \hline Cosymplectic  $\bar{M}$&$\mathfrak{X}_{\mathrm{cs-ham}}(\bar{M})$ in \eqref{X-cs} & 
$\mathfrak{X}_{\mathrm{cs-ham}}^*(\bar{M})$ in \eqref{X-cs-ham}    \\  \hline 
Contact $N$& $\mathfrak{X}_{\mathrm{c-ham}}(N)$ in \eqref{X-con} & 
$\mathfrak{X}_{\mathrm{c-ham}}^*(N)$ in \eqref{c-Ham-alg}  \\  \hline 
Cocontact $\bar{N}$& $\mathfrak{X}_{\mathrm{cc-ham}}(\bar{N})$ in \eqref{X-cc-ham} & 
$\mathfrak{X}_{\mathrm{cc-ham}}^*(\bar{N})$ in \eqref{cc-Ham-alg}  \\ 
 \hline
\end{tabular}
\captionof{table}{Lie Algebras and Dual Spaces\label{kinetic-algebra-table}}
\end{center}

\begin{center}
\begin{tabular}{ ||c||c|c||  }
 \hline
 \multicolumn{3}{||c||}{\textbf{Kinetic Dynamics}} \\
 \hline \hline
Manifold / Coord. & Lie-Poisson in Momenta & Density and Lie-Poisson Equation\\
 \hline  \hline 
 Symplectic $M$  &  $\Pi =
 \mathcal{L}_{X_H}\Pi$&  $f =\mathrm{div}\sharp^s(\Pi)$\\ $(q^i,p_i)$ &&  $\dot{f} =\{H,f\}^s 
 $  \\  \hline Cosymplectic  $\bar{M}$  &$\dot{\bar{\Pi}}=-\mathcal{L}_{X^{cs}_H }\bar{\Pi}$ & $f=\mathrm{div} (\sharp^{cs} \bar{\Pi})   + \mathcal{L}_{\mathcal{R}^\tau }  \langle \bar{\Pi},\mathcal{R}^\tau\rangle $ \\   $(t,q^i,p_i)$ && $\dot{f} =\{H,f\}^{cs}$ \\
 \hline 
Contact $N$  &  $\dot{\Upsilon}= - \mathcal{L}_{X_{H}^{c} } \Upsilon +  (n+1)  \mathcal{R}^\eta(H)\Upsilon$&  $f=\mathrm{div}({\sharp^{c}\Upsilon}) 
  -   \mathcal{L}_{\mathcal{R}^\eta} \langle \Upsilon, \mathcal{R}^\eta \rangle   - \langle \Upsilon, \mathcal{R}^\eta \rangle$
   \\ $(q^i,p_i,z)$ & & $\dot{f} = - \{H,f\}^{c} + (n+2) f \mathcal{R}^\eta(H)$ \\
 \hline 
Cocontact $\bar{N}$ & $\dot{\bar{\Upsilon} }= - \mathcal{L}_{X_{H}^{cc} } \bar{\Upsilon} +  (n+1)  \mathcal{R}^\eta(H)\bar{\Upsilon} $& $f= \mathrm{div}({\sharp^{cc}\bar{\Upsilon} }) 
  -   \mathcal{L}_{\mathcal{R}^\eta} \langle \bar{\Upsilon} , \mathcal{R}^\eta \rangle  -    \mathcal{L}_{\mathcal{R}^\tau} \langle \bar{\Upsilon} , \mathcal{R}^\tau \rangle  - \langle \bar{\Upsilon} , \mathcal{R}^\eta \rangle$
    \\ $(t,q^i,p_i,z)$ & & $\dot{f} = - \{H,f\}^{cc} + (n+2) f \mathcal{R}^\eta(H)$\\
 \hline
\end{tabular}
\captionof{table}{Kinetic Dynamics \label{kinetic-table}}
\end{center}

\end{landscape}

\appendix 

\section{Fundamental Constructions}\label{Appapp}

   \subsection{Jacobi Manifolds and Dynamics}
    A Jacobi manifold \cite{Kirillov,Lich,Marle-Jacobi,vaisman2002jacobi} is a triplet $(P,\Lambda^J,\epsilon)$, where $\Lambda^J$ is a bivector field (a skew-symmetric contravariant 2-tensor field) and $\epsilon$ is a vector field satisfying
\begin{equation}\label{ident-Jac}
     [\Lambda^J,\Lambda^J] = 2 \epsilon \wedge \Lambda^J, \qquad 
         [\epsilon,\Lambda^J] = 0,
\end{equation} 
    where the bracket is the Schouten–Nijenhuis bracket. We can define a Jacobi bracket on the space $C^\infty(P)$ of smooth functions on $M$ by 
\begin{equation}\label{Jacobi-bracket}
   \{F,H\}^{J}= \Lambda^J(dF, dH)+ F \epsilon(H) - H \epsilon (F).
\end{equation} 
The bracket \eqref{Jacobi-bracket} is skew-symmetric and satisfies  the Jacobi identity but not the Leibniz identity. Instead, it is a first-order differential operator since
\begin{equation}\label{Jacobi-bracket-}
\{F, KH\}^{J} = K \{F, H\}^{J}+ H \{F, K\}^{J} + KH \epsilon(F),
\end{equation}
for all $F$, $K$ and $H$. This is equivalently written as the weak Leibniz identity
\begin{equation}
        \operatorname{supp}(\{H,F\}^{J}) \subseteq \operatorname{supp} (H) \cap \text{supp} (F).
\end{equation}

For a Hamiltonian function $H$ on a Jacobi manifold $(P,\Lambda^J,\epsilon)$, the Hamiltonian vector field is defined to be
\begin{equation}\label{Ham-v-f-Jac}
X_H^J=\sharp_{\Lambda^J}(dH)-H\epsilon.
\end{equation}
The relationship between the Hamiltonian vector field and the Jacobi bracket is given by the equation
\begin{equation}\label{bra-X-J}
X_H^J(F)= \{F,H\}^{J}- F \epsilon(H).
\end{equation} 
The space of Hamiltonian vector fields is closed; moreover, we have that
 \begin{equation}\label{Jac-integ}
[X_F,X_H]=-X_{\{F,H\}^{J}},
\end{equation}
where the bracket on the left-hand side is the Jacobi-Lie bracket of vector fields. 
This result leads us to determine the Lie algebra 
 \begin{equation}
\mathfrak{X}_{\mathrm{J-ham}}(P)= \big\{X_{H}^J \in \mathfrak{X}(M ): X_{H}^J =\sharp_{\Lambda^J}(dH)-H\epsilon \big\}
\end{equation}
of Hamiltonian vector fields on the Jacobi manifold. Further, by taking the Lie algebra bracket on the space $\mathfrak{X}_{\mathrm{J-ham}}(P)$ as the minus of the Jacobi-Lie bracket of vector fields, the mapping 
 \begin{equation}
 \vartheta:C^\infty(P)\longrightarrow \mathfrak{X}_{\mathrm{J-ham}}(P),\qquad H\mapsto X_H^J
\end{equation}
turns out to be a Lie algebra homomorphism.

\subsection{Almost Poisson Manifolds and Dynamics}

Consider a manifold $P$ equipped with a skew-symmetric $\mathbb{R}$-bilinear bracket  \begin{equation}\label{PoissonBracket}
  \{\bullet,\bullet\}:C^{\infty}(P) \times C^{\infty}(P)\longrightarrow C^{\infty}(P)
\end{equation}
defined on the space $C^{\infty}(P)$ of smooth functions. The bracket is called an almost Poisson bracket if it satisfies the Leibniz identity. 
In this case, the pair $(P,\{\bullet,\bullet\})$ is called an almost Poisson manifold \cite{BhaskaraViswanath,Dufour,Lichnerowicz-Poi,Vaisman94,weinstein1983local,Weinstein98}. An almost Poison bracket is said to be Poisson if it satisfies the Jacobi identity. In this case, the pair $(P,\{\bullet,\bullet\})$ is called a Poisson manifold. 

One may identify an almost Poisson bracket $\{\bullet,\bullet\}$ with a bivector field $\Lambda$ on $P$ according to the following definition 
\begin{equation} \label{bivec-PoissonBra}
\Lambda(dF,dH):=\{F,H\}
\end{equation}
for all $F$ and $H$ in $C^{\infty}(P)$. One can extend this definition to all one-forms by linearity. So, we may represent an almost Poisson manifold by a pair $(P,\Lambda)$ as 
well. This geometry enables us to define a musical mapping $\sharp_{\Lambda}$  induced by the bivector field $\Lambda$ as
\begin{equation} \label{sharp2}
\sharp_{\Lambda}: \Gamma^1(P)\longrightarrow \mathfrak{X}(P),  \qquad \sharp_{\Lambda}(\alpha)( \beta) := \Lambda(\beta,\alpha),
\end{equation}
where the pairing on the left-hand side is the one between the space $\Gamma^1(P)$ of one-form sections and the space $\mathfrak{X}(P)$ of vector fields. The bracket $\{\bullet,\bullet\}$ satisfies the Jacobi identity if and only if the bivector field $\Lambda$ commutes with itself under the Schouten-Nijenhuis bracket, that is, 
\begin{equation} \label{Poisson-cond}
[\Lambda,\Lambda]=0.
\end{equation} 
A bivector field is a Poisson bivector field if the condition \eqref{Poisson-cond} holds. So, we can denote a Poisson manifold by a pair $(P,\Lambda)$ with  $\Lambda$ satisfying the Jacobi identity \eqref{Poisson-cond}. Accordingly, a Poisson manifold is a particular instance of a Jacobi manifold where the vector field vanishes.

Consider an (almost) Poisson manifold $(P,\{\bullet,\bullet\})$. Given a Hamiltonian function $H$ (any real-valued function on $P$), the Hamiltonian vector field $X_H$ is defined to be, for all $F$ in $C^{\infty}(P)$,
\begin{equation}\label{Hamvf-}
X_H(F):=\{F,H\},
\end{equation}
where the notation $X_H(F)$ stands for the directional derivative of $F$ in the direction of $X_H$. In terms of $\sharp_{\Lambda}$, the Hamiltonian vector field $X_H$ can be defined as 
\begin{equation}\label{Poisson-ham}
\sharp_{\Lambda}(dH)= X_H.
\end{equation}
The Hamiltonian dynamics generated by a Hamiltonian function $H$ is then  
\begin{equation} \label{HamEq}
\dot{x}=\{x,H\},
\end{equation} 
for $x$ in $P$.
The skew-symmetry of the bracket implies that the Hamiltonian function $H$ is conserved all along the motion that is $\dot{H}=0$. This physically corresponds to the conservation of energy.

 \subsection{Lie-Poisson Dynamics of Kinetic Theories} \label{App-LP-Kinetic}
   
Let $G$ be a Lie group, $\mathfrak{g}$ be its Lie algebra with skew-symmetric bilinear bracket
\begin{equation}
  [\bullet,\bullet]:\mathfrak{g}\times \mathfrak{g}\mapsto \mathfrak{g}
\end{equation}
satisfying the Jacobi identity \cite{GF,vara84} . Referring to this Lie bracket, one can define (left) adjoint and (right) coadjoint actions of the Lie algebra $\mathfrak{g}$ on itself and the dual space as 
\begin{equation}\label{coad-gen}
\begin{split}
ad &:\mathfrak{g}\times \mathfrak{g}\mapsto \mathfrak{g},\qquad ad_\xi\eta:=[\xi,\eta],
\\
ad^* &:\mathfrak{g}\times \mathfrak{g}^*\mapsto \mathfrak{g}^*,\qquad \langle ad^{\ast }_\xi \rho, \eta \rangle = 
 \langle \rho , ad_\xi \eta \rangle,
\end{split}
\end{equation}
respectively.
The dual space $\mathfrak{g}^*$ is equipped with a Lie Poisson bracket \cite{Marsden1999,holm11,holm2009geometric}, defined as follows: for two functionals $\mathcal{F}$ and $\mathcal{H}$ defined on $\mathfrak{g}^*$, we put
\begin{equation}
\left\{ \mathcal{F},\mathcal{H}\right\} _{\mathfrak{g}^{\ast }} (
\rho  ) = \Big\langle \rho ,\left[ \frac{\delta \mathcal{F}}{%
\delta \rho},\frac{\delta \mathcal{H}}{\delta \rho }\right] \Big\rangle.  \label{LPBr}
\end{equation}%
One may consider the minus Lie-Poisson bracket as well. In this work, we prefer the positive one. For a Hamiltonian functional $\mathcal{H}$, the dynamics is governed by the Lie-Poisson
equations computed in terms of the coadjoint action as
\begin{equation}
\dot{\rho}= - ad_{\delta \mathcal{H} / \delta \rho }^{\ast }\rho. 
\label{LPEq}
\end{equation}

\noindent
\textbf{Lie-Poisson Dynamics of Kinetic Theories.}
Consider a bunch of particles rest in a manifold $M\subset \mathbb{R}^3$. $\mathrm{Diff}(M)$  is the group of all diffeomorphisms on $M$, \cite{ArKh98,banyaga97,ebin70,marsden82}. The motion of the particles is determined by the left action of $\mathrm{Diff}(M)$ on the particle space $M$. The right action commutes with the particle motion and constitutes
an infinite dimensional symmetry group 
called the particle relabelling symmetry. The Lie algebra of $\mathrm{Diff}(M)$ is the space $\mathfrak{X}(M)$ of vector fields. 
We define the dual space $\mathfrak{X}^{\ast }(M) $ of the Lie algebra as space of one-form densities $\Lambda ^{1}(M) \otimes \mathrm{Den}(M) $ on $M$. Here, the pairing between a vector field $X$ and a dual element $\Pi\otimes \mu$ is defined to be $L_2$ (simply multiply-and-integrate form) pairing
\begin{equation}\label{L-2}
\langle \bullet, \bullet \rangle_{L_2}: \Gamma^1(M)\otimes \mathrm{Den}(M)\times \mathfrak{X}(M)\longrightarrow \mathbb{R}, \qquad (\Pi\otimes \mu,X)\mapsto \int_{M} \langle\Pi , X \rangle   \mu.   
\end{equation}
Here, the pairing inside the integral is between the one-form $\Pi$ and the vector field $X$,   $\mu$ is a density, that is, a volume form, on $M$.
The adjoint and coadjoint actions are 
\begin{equation} \label{Lie-brr}
\begin{split}
ad_X Y & = [X,Y]_{\mathfrak{X}\left( M%
\right)}=-[X,Y]_{JL},
\\ad_X^\ast (\Pi\otimes \mu) & = (\mathcal{L}_X \Pi + \mathrm{div}(X)\Pi) \otimes \mu,
\end{split}
\end{equation}
respectively. Here, $\mathrm{div}(X)$ is the divergence of the vector field with respect to the volume form $\mu$. 

If a particle in a medium is flowing, following the dynamics generated by the vector field $X$ defined on the manifold $M$, then the kinetic dynamics of the medium can be understood as the lift of this motion. To describe it, we introduce a linear Hamiltonian functional defined as follows
\begin{equation}
\mathcal{H}(\Pi \otimes \mu) = \int_M \langle \Pi, X \rangle \mu,
\end{equation}
This functional operates on the space of one-form densities $\Gamma^1(M) \otimes \mathrm{Den}(M)$. We can prove that $\delta H/\delta \Pi = X$. When the volume is fixed, the Lie-Poisson equation takes the form
\begin{equation}\label{LP-gen}
\dot{\Pi} = -\mathcal{L}_{X}\Pi - \mathrm{div}(X)\Pi.
\end{equation}
However, if the dynamics is generated by a divergence-free vector field (such as in the case of incompressible fluid flow or Vlasov flow), then the second term on the right-hand side of \eqref{LP-gen} vanishes, and we obtain
\begin{equation}
\dot{\Pi} = -\mathcal{L}_{X}\Pi. \label{LP}
\end{equation}
This intermediate level represents the momentum formulation of kinetic dynamics. To express the dynamics in terms of density functions, one must establish a Poisson mapping from the space of one-forms to the space of densities. In the present work, this has been accomplished for Hamiltonian dynamics through the dualization of the Lie algebra homomorphism $H\mapsto X_H.$

\section*{Acknowledgements}
We acknowledge the financial support of Grant PID2019-106715GBC21, and the Severo Ochoa Programme for Centres of Excellence in R\&D (CEX2019-000904-S).

\bibliographystyle{abbrv}
\bibliography{references}

\begin{thebibliography}{10}

\bibitem{AbMa78}
R.~Abraham and J.~E. Marsden.
\newblock {\em Foundations of mechanics}.
\newblock Benjamin/Cummings Publishing Co., Inc., Advanced Book Program,
  Reading, Mass., 1978.

\bibitem{Arn}
V.~I. Arnold.
\newblock {\em Mathematical methods of classical mechanics}, volume~60 of {\em
  Graduate Texts in Mathematics}.
\newblock Springer-Verlag, New York, second edition, 1989.
\newblock Translated from the Russian by K. Vogtmann and A. Weinstein.

\bibitem{ArKh98}
V.~I. Arnold and B.~A. Khesin.
\newblock {\em Topological methods in hydrodynamics}, volume 125 of {\em
  Applied Mathematical Sciences}.
\newblock Springer-Verlag, New York, 1998.

\bibitem{AtEsLe23}
B.~Ate\c{s}li, O.~Esen, M.~de~Le\'{o}n, and C.~Sard\'{o}n.
\newblock On locally conformally cosymplectic {H}amiltonian dynamics and
  {H}amilton-{J}acobi theory.
\newblock {\em J. Phys. A}, 56(1):Paper No. 015204, 40, 2023.

\bibitem{azuaje2023canonical}
R.~Azuaje and A.~Escobar-Ruiz.
\newblock Canonical and canonoid transformations for hamiltonian systems on
  (co) symplectic and (co) contact manifolds.
\newblock {\em Journal of Mathematical Physics}, 64(3), 2023.

\bibitem{banyaga97}
A.~Banyaga.
\newblock {\em The structure of classical diffeomorphism groups}, volume 400 of
  {\em Mathematics and its Applications}.
\newblock Kluwer Academic Publishers Group, Dordrecht, 1997.

\bibitem{BazzGoer15}
G.~Bazzoni and O.~Goertsches.
\newblock K-cosymplectic manifolds.
\newblock {\em Annals of Global Analysis and Geometry}, 47:239--270, 2015.

\bibitem{BhaskaraViswanath}
K.~H. Bhaskara and K.~Viswanath.
\newblock {\em Poisson algebras and {P}oisson manifolds}, volume 174 of {\em
  Pitman Research Notes in Mathematics Series}.
\newblock Longman Scientific \& Technical, Harlow; copublished in the United
  States with John Wiley \& Sons, Inc., New York, 1988.

\bibitem{Br17}
A.~Bravetti.
\newblock Contact {H}amiltonian dynamics: the concept and its use.
\newblock {\em Entropy}, 19(10):Paper No. 535, 12, 2017.

\bibitem{Bravettithermo}
A.~Bravetti.
\newblock Contact geometry and thermodynamics.
\newblock {\em Int. J. Geom. Methods Mod. Phys.}, 16:1940003, 51, 2019.

\bibitem{BrCrTa17}
A.~Bravetti, H.~Cruz, and D.~Tapias.
\newblock Contact {H}amiltonian mechanics.
\newblock {\em Ann. Physics}, 376:17--39, 2017.

\bibitem{Leon-GradOnCosymp}
F.~Cantrijn, M.~de~Le{\'o}n, and E.~Lacomba.
\newblock Gradient vector fields on cosymplectic manifolds.
\newblock {\em Journal of Physics A: Mathematical and General}, 25(1):175,
  1992.

\bibitem{Cape}
B.~Cappelletti-Montano, A.~De~Nicola, and I.~Yudin.
\newblock A survey on cosymplectic geometry.
\newblock {\em Rev. Math. Phys.}, 25(10):1343002, 55, 2013.

\bibitem{ChLeMa91}
D.~Chinea, M.~de~Le\'{o}n, and J.~C. Marrero.
\newblock Locally conformal cosymplectic manifolds and time-dependent
  {H}amiltonian systems.
\newblock {\em Comment. Math. Univ. Carolin.}, 32(2):383--387, 1991.

\bibitem{LeGaGrMuRi23}
M.~de~Le\'{o}n, J.~Gaset, X.~Gr\`acia, M.~C. Mu\~{n}oz Lecanda, and X.~Rivas.
\newblock Time-dependent contact mechanics.
\newblock {\em Monatsh. Math.}, 201(4):1149--1183, 2023.

\bibitem{LeGaMuRiRo23}
M.~de~Le\'{o}n, J.~Gaset, M.~C. Mu\~{n}oz Lecanda, X.~Rivas, and
  N.~Rom\'{a}n-Roy.
\newblock Multicontact formulation for non-conservative field theories.
\newblock {\em J. Phys. A}, 56(2):Paper No. 025201, 44, 2023.

\bibitem{LeLa19}
M.~de~Le\'{o}n and M.~Lainz.
\newblock Contact {H}amiltonian systems.
\newblock {\em J. Math. Phys.}, 60(10):102902, 18, 2019.

\bibitem{Leon-cocontact23}
M.~de~Le\'{o}n, M.~Lainz, A.~L\'{o}pez-Gord\'{o}n, and X.~Rivas.
\newblock Hamilton-{J}acobi theory and integrability for autonomous and
  non-autonomous contact systems.
\newblock {\em J. Geom. Phys.}, 187:Paper No. 104787, 22, 2023.

\bibitem{de2021hamilton}
M.~de~Le{\'o}n, M.~Lainz, and A.~Mu\~niz Brea.
\newblock {The Hamilton--Jacobi theory for contact Hamiltonian systems}.
\newblock {\em Mathematics}, 9(16):1993, 2021.

\bibitem{de2019contact}
M.~de~Le{\'o}n and M.~Lainz~Valc{\'a}zar.
\newblock Contact {H}amiltonian systems.
\newblock {\em Journal of Mathematical Physics}, 60(10):102902, 2019.

\bibitem{leon89}
M.~de~Le\'{o}n and P.~R. Rodrigues.
\newblock {\em Methods of differential geometry in analytical mechanics},
  volume 158 of {\em North-Holland Mathematics Studies}.
\newblock North-Holland Publishing Co., Amsterdam, 1989.

\bibitem{Leon-CosympReduction}
M.~de~Le{\'o}n and M.~Saralegi.
\newblock Cosymplectic reduction for singular momentum maps.
\newblock {\em Journal of Physics A: Mathematical and General}, 26(19):5033,
  1993.

\bibitem{LeonSar2}
M.~de~Le\'{o}n and C.~Sard\'{o}n.
\newblock Cosymplectic and contact structures for time-dependent and
  dissipative {H}amiltonian systems.
\newblock {\em J. Phys. A}, 50(25):255205, 23, 2017.

\bibitem{LeonTuyn}
M.~de~Le\'{o}n and G.~M. Tuynman.
\newblock A universal model for cosymplectic manifolds.
\newblock {\em J. Geom. Phys.}, 20(1):77--86, 1996.

\bibitem{Dufour}
J.~Dufour and N.~T. Zung.
\newblock {\em Poisson structures and their normal forms}, volume 242 of {\em
  Progress in Mathematics}.
\newblock Birkh\"{a}user Verlag, Basel, 2005.

\bibitem{ebin70}
D.~G. Ebin and J.~Marsden.
\newblock Groups of diffeomorphisms and the motion of an incompressible fluid.
\newblock {\em Ann. of Math. (2)}, 92:102--163, 1970.

\bibitem{EsLeSaZa-Cauchy}
O.~Esen, M.~de~Le\'{o}n, C.~Sard\'{o}n, and M.~Zajac.
\newblock Cauchy data space and multisymplectic formulation of conformal
  classical field theories.
\newblock {\em Ann. Physics}, 434:Paper No. 168616, 26, 2021.

\bibitem{EsLeSaZa-k-sympl}
O.~Esen, M.~de~Le\'{o}n, C.~Sard\'{o}n, and M.~Zajac.
\newblock The globalization problem of the {H}amilton-{D}edonder-{W}eyl
  equations on a local {$k$}-symplectic framework.
\newblock {\em Mediterr. J. Math.}, 18(1):Paper No. 26, 25, 2021.

\bibitem{EsLeSaZa21}
O.~Esen, M.~de~Le\'{o}n, C.~Sard\'{o}n, and M.~Zajac.
\newblock Hamilton-{J}acobi formalism on locally conformally symplectic
  manifolds.
\newblock {\em J. Math. Phys.}, 62(3):Paper No. 033506, 15, 2021.

\bibitem{esen2023conformal}
O.~Esen, A.~Gezici, M.~Grmela, H.~G{\"u}mral, M.~Pavelka, and S.~S{\"u}tl{\"u}.
\newblock Conformal and contact kinetic dynamics and their geometrization.
\newblock {\em arXiv preprint arXiv:2307.06080}, 2023.

\bibitem{EsGrMiGu19}
O.~Esen, M.~Grmela, H.~G\"{u}mral, and M.~Pavelka.
\newblock Lifts of symmetric tensors: fluids, plasma, and {G}rad hierarchy.
\newblock {\em Entropy}, 21(9):Paper No. 907, 33, 2019.

\bibitem{EsGu11}
O.~Esen and H.~G\"{u}mral.
\newblock Lifts, jets and reduced dynamics.
\newblock {\em Int. J. Geom. Methods Mod. Phys.}, 8(2):331--344, 2011.

\bibitem{EsGu12}
O.~Esen and H.~G\"{u}mral.
\newblock Geometry of plasma dynamics {II}: {L}ie algebra of {H}amiltonian
  vector fields.
\newblock {\em J. Geom. Mech.}, 4(3):239--269, 2012.

\bibitem{esen2021contact}
O.~Esen, M.~Lainz, M.~de~Le{\'o}n, and J.~C. Marrero.
\newblock Contact dynamics: {L}egendrian and {L}agrangian submanifolds.
\newblock {\em Mathematics}, 9(21):2704, 2021.

\bibitem{esen2021implicit}
O.~Esen, M.~Lainz~Valc\'{a}zar, M.~de~Le\'{o}n, and C.~Sard\'{o}n.
\newblock Implicit contact dynamics and {H}amilton-{J}acobi theory.
\newblock {\em Differential Geom. Appl.}, 90:Paper No. 102030, 2023.

\bibitem{EsSu21}
O.~Esen and S.~S\"{u}tl\"{u}.
\newblock Matched pair analysis of the {V}lasov plasma.
\newblock {\em J. Geom. Mech.}, 13(2):209--246, 2021.

\bibitem{GF}
M.~Fecko.
\newblock {\em Differential geometry and Lie groups for physicists}.
\newblock Cambridge university press, 2006.

\bibitem{Goto15}
S.-i. Goto.
\newblock Legendre submanifolds in contact manifolds as attractors and
  geometric nonequilibrium thermodynamics.
\newblock {\em J. Math. Phys.}, 56(7):073301, 30, 2015.

\bibitem{grabowska2022contact}
K.~Grabowska and J.~Grabowski.
\newblock Contact geometric mechanics: the tulczyjew triples.
\newblock {\em arXiv preprint arXiv:2209.03154}, 2022.

\bibitem{Grmela14}
M.~Grmela.
\newblock Contact geometry of mesoscopic thermodynamics and dynamics.
\newblock {\em Entropy}, 16(3):1652--1686, 2014.

\bibitem{guil90}
V.~Guillemin and S.~Sternberg.
\newblock {\em Symplectic techniques in physics}.
\newblock Cambridge University Press, Cambridge, second edition, 1990.

\bibitem{Gu10}
H.~G\"{u}mral.
\newblock Geometry of plasma dynamics. {I}. {G}roup of canonical
  diffeomorphisms.
\newblock {\em J. Math. Phys.}, 51(8):083501, 23, 2010.

\bibitem{holm11}
D.~D. Holm.
\newblock {\em Geometric mechanics. {P}art {I}}.
\newblock Imperial College Press, London, second edition, 2011.
\newblock Dynamics and symmetry.

\bibitem{holm2009geometric}
D.~D. Holm, T.~Schmah, and C.~Stoica.
\newblock {\em Geometric mechanics and symmetry}, volume~12 of {\em Oxford
  Texts in Applied and Engineering Mathematics}.
\newblock Oxford University Press, Oxford, 2009.
\newblock From finite to infinite dimensions, With solutions to selected
  exercises by David C. P. Ellis.

\bibitem{Kirillov}
A.~A. Kirillov.
\newblock Local {L}ie algebras.
\newblock {\em Uspehi Mat. Nauk}, 31(4(190)):57--76, 1976.

\bibitem{Liber87}
P.~Libermann and C.~M. Marle.
\newblock {\em Symplectic geometry and analytical mechanics}, volume~35 of {\em
  Mathematics and its Applications}.
\newblock D. Reidel Publishing Co., Dordrecht, 1987.

\bibitem{Lich}
A.~Lichnerowicz.
\newblock Sur la reductivite de certaines algebres d’atomorphismes.
\newblock {\em C.R. Acad. Sci. Paris}, 253:1302--1304, 1961.

\bibitem{Lichnerowicz-Poi}
A.~Lichnerowicz.
\newblock Les vari\'{e}t\'{e}s de {P}oisson et leurs alg\`ebres de {L}ie
  associ\'{e}es.
\newblock {\em J. Differential Geometry}, 12(2):253--300, 1977.

\bibitem{Lichnerowicz-Jacobi}
A.~Lichnerowicz.
\newblock Les vari\'{e}t\'{e}s de {J}acobi et leurs alg\`ebres de {L}ie
  associ\'{e}es.
\newblock {\em J. Math. Pures Appl. (9)}, 57(4):453--488, 1978.

\bibitem{Marle-Jacobi}
C.~Marle.
\newblock On {J}acobi manifolds and {J}acobi bundles.
\newblock In {\em Symplectic geometry, groupoids, and integrable systems
  ({B}erkeley, {CA}, 1989)}, volume~20 of {\em Math. Sci. Res. Inst. Publ.},
  pages 227--246. Springer, New York, 1991.

\bibitem{marsden82}
J.~E. Marsden.
\newblock A group theoretic approach to the equations of plasma physics.
\newblock {\em Canad. Math. Bull.}, 25(2):129--142, 1982.

\bibitem{marsden1983hamiltonian}
J.~E. Marsden, T.~Ratiu, R.~Schmid, R.~Spencer, and A.~J. Weinstein.
\newblock Hamiltonian systems with symmetry, coadjoint orbits and plasma
  physics.
\newblock {\em Atti della Accademia delle scienze di Torino}, 117(1):289--340,
  1983.

\bibitem{marsden83b}
J.~E. Marsden, T.~Ratiu, and A.~Weinstein.
\newblock Reduction and {H}amiltonian structures on duals of semidirect product
  {L}ie algebras.
\newblock In {\em Fluids and plasmas: geometry and dynamics ({B}oulder,
  {C}olo., 1983)}, volume~28 of {\em Contemp. Math.}, pages 55--100. Amer.
  Math. Soc., Providence, RI, 1984.

\bibitem{Marsden1999}
J.~E. Marsden and T.~S. Ratiu.
\newblock {\em Introduction to mechanics and symmetry}, volume~17 of {\em Texts
  in Applied Mathematics}.
\newblock Springer-Verlag, New York, second edition, 1999.
\newblock A basic exposition of classical mechanical systems.

\bibitem{marsden1982hamiltonian}
J.~E. Marsden and A.~Weinstein.
\newblock The {H}amiltonian structure of the {M}axwell-{V}lasov equations.
\newblock {\em Phys. D}, 4(3):394--406, 1981/82.

\bibitem{MaNo17}
A.~Mastromartino and Y.~Nogier.
\newblock On the contact and co-contact of higher order.
\newblock {\em Mat. Vesnik}, 69(2):89--100, 2017.

\bibitem{McSa17}
D.~McDuff and D.~Salamon.
\newblock {\em Introduction to symplectic topology}.
\newblock Oxford Graduate Texts in Mathematics. Oxford University Press,
  Oxford, third edition, 2017.

\bibitem{morrison80}
P.~Morrison.
\newblock The {M}axwell-{V}lasov equations as a continuous {H}amiltonian
  system.
\newblock {\em Phys. Lett. A}, 80(5-6):383--386, 1980.

\bibitem{morrison1981hamiltonian}
P.~J. Morrison.
\newblock Hamiltonian field description of the one-dimensional poisson-vlasov
  equations.
\newblock Technical report, Princeton Univ., NJ (USA). Plasma Physics Lab.,
  1981.

\bibitem{morrison1982poisson}
P.~J. Morrison.
\newblock Poisson brackets for fluids and plasmas.
\newblock In {\em AIP Conference proceedings}, volume~88, pages 13--46.
  American Institute of Physics, 1982.

\bibitem{Mrugala}
R.~Mrugala, J.~D. Nulton, J.~C. Sch\"{o}n, and P.~Salamon.
\newblock Contact structure in thermodynamic theory.
\newblock {\em Rep. Math. Phys.}, 29(1):109--121, 1991.

\bibitem{RaSc81}
T.~S. Ra\c{t}iu and R.~Schmid.
\newblock The differentiable structure of three remarkable diffeomorphism
  groups.
\newblock {\em Math. Z.}, 177(1):81--100, 1981.

\bibitem{Ri23}
X.~Rivas.
\newblock Nonautonomous {$k$}-contact field theories.
\newblock {\em J. Math. Phys.}, 64(3):Paper No. 033507, 20, 2023.

\bibitem{RiTo23}
X.~Rivas and D.~Torres.
\newblock Lagrangian-{H}amiltonian formalism for cocontact systems.
\newblock {\em J. Geom. Mech.}, 15(1):1--26, 2023.

\bibitem{Houches}
A.~A. Simoes, D.~M. de~Diego, M.~L. Valc\'{a}zar, and M.~de~Le\'{o}n.
\newblock The geometry of some thermodynamic systems.
\newblock In {\em Geometric structures of statistical physics, information
  geometry, and learning}, volume 361 of {\em Springer Proc. Math. Stat.},
  pages 247--275. Springer, Cham, 2021.

\bibitem{tul77}
W.~M. Tulczyjew.
\newblock The {L}egendre transformation.
\newblock {\em Ann. Inst. H. Poincar\'{e} Sect. A (N.S.)}, 27(1):101--114,
  1977.

\bibitem{TuUr99}
W.~o.~M. Tulczyjew and P.~Urba\'{n}ski.
\newblock A slow and careful {L}egendre transformation for singular
  {L}agrangians.
\newblock {\em Acta Phys. Polon. B}, 30(10):2909--2978, 1999.
\newblock The Infeld Centennial Meeting (Warsaw, 1998).

\bibitem{Vaisman85}
I.~Vaisman.
\newblock Locally conformal symplectic manifolds.
\newblock {\em Internat. J. Math. Math. Sci.}, 8(3):521--536, 1985.

\bibitem{Vaisman94}
I.~Vaisman.
\newblock {\em Lectures on the geometry of {P}oisson manifolds}, volume 118 of
  {\em Progress in Mathematics}.
\newblock Birkh\"{a}user Verlag, Basel, 1994.

\bibitem{vaisman2002jacobi}
I.~Vaisman.
\newblock Jacobi manifolds.
\newblock {\em Selected topics in Geometry and Mathematical Physics},
  1:81--100, 2002.

\bibitem{vara84}
V.~S. Varadarajan.
\newblock {\em Lie groups, {L}ie algebras, and their representations}, volume
  102 of {\em Graduate Texts in Mathematics}.
\newblock Springer-Verlag, New York, 1984.
\newblock Reprint of the 1974 edition.

\bibitem{weinstein1983local}
A.~Weinstein.
\newblock The local structure of {P}oisson manifolds.
\newblock {\em J. Differential Geom.}, 18(3):523--557, 1983.

\bibitem{Weinstein98}
A.~Weinstein.
\newblock Poisson geometry.
\newblock {\em Differential Geom. Appl.}, 9(1-2):213--238, 1998.
\newblock Symplectic geometry.

\end{thebibliography}

\end{document}